\documentclass[12pt]{amsart}

\usepackage[matrix,arrow,curve,frame]{xy}
\usepackage{amsmath,amsthm,amssymb,enumerate}
\usepackage{latexsym}
\usepackage{amscd}
\usepackage[colorlinks=true]{hyperref}
\usepackage{euscript}
\usepackage{todonotes}
\usepackage{mathrsfs}
\usepackage{xypic,calc}

\newtheorem*{thm*}{Theorem}
\newtheorem{thm}[subsection]{Theorem}
\newtheorem{defn}[subsection]{Definition}

\newtheorem{claim}[subsection]{Claim}

\newtheorem{Const}[subsection]{Construction}
\newtheorem{lemma}[subsection]{Lemma}

\newtheorem{remark}{Remark}

\newtheorem{Error}{Error}

\theoremstyle{definition}
\newtheorem{example}[subsection]{Example}

\newcommand{\cat}{\mathcal}
\newcommand{\lra}{\longrightarrow}
\newcommand{\lla}{\longleftarrow}
\newcommand{\llra}[1]{\stackrel{#1}{\lra}}
\newcommand{\llla}[1]{\stackrel{#1}{\lla}}
\newcommand{\sSet}{s\cat Set}

\newcommand{\Q}{\mathbb Q}
\newcommand{\Z}{\mathbb Z}
\newcommand{\F}{\mathbb F}
\newcommand{\N}{\mathbb N}
\newcommand{\W}{\mathbb W}

\newcommand{\C}{\mathbb C}

\DeclareMathOperator{\colim}{colim}
\DeclareMathOperator{\hocolim}{hocolim}
\DeclareMathOperator{\holim}{holim}

\DeclareMathOperator*{\invlim}{\varprojlim}
\DeclareMathOperator*{\dirlim}{\varinjlim}

\DeclareMathOperator{\BG}{BG}

\DeclareMathOperator{\pic}{pic}
\DeclareMathOperator{\Pic}{Pic}
\DeclareMathOperator{\B}{B}
\DeclareMathOperator{\G}{G}

\DeclareMathOperator{\U}{U}

\DeclareMathOperator{\Id}{Id}
\DeclareMathOperator{\BUn}{BU}

\DeclareMathOperator{\Ne}{Nr}
\DeclareMathOperator{\BC}{B \mathscr{C}}
\DeclareMathOperator{\BA}{B \mathscr{A}}

\DeclareMathOperator{\HZ}{H\Z}
\DeclareMathOperator{\Hh}{H}

\DeclareMathOperator{\BGx}{\mathscr{B}G}
\DeclareMathOperator{\BSm}{\mathscr{B}\Sigma}

\DeclareMathOperator{\BSi}{B\Sigma}
\DeclareMathOperator{\rBGx}{r\mathscr{B}G}
\DeclareMathOperator{\rB}{r\mathscr{B}}
\DeclareMathOperator{\rE}{r\mathscr{E}}
\DeclareMathOperator{\rtB}{r\tilde{\mathscr{B}}}
\DeclareMathOperator{\rtE}{r\tilde{\mathscr{E}}}

\DeclareMathOperator{\BGL}{\mathscr{B}GL}
\DeclareMathOperator{\BGLL}{BGL}
\DeclareMathOperator{\BGLs}{\mathscr{B}\mathcal{G}L}
\DeclareMathOperator{\GLs}{\mathcal{G}L}
\DeclareMathOperator{\ENo}{\mbox{EN}}
\DeclareMathOperator{\EGL}{\mbox{EGL}}
\DeclareMathOperator{\EGg}{\mbox{EG}}
\DeclareMathOperator{\BGb}{BGL}
\DeclareMathOperator{\GL}{GL}

\DeclareMathOperator{\Fr}{Fb}

\DeclareMathOperator{\BN}{\mathscr{B}N}
\DeclareMathOperator{\BS}{B\Sigma}

\DeclareMathOperator{\No}{N}

\DeclareMathOperator{\Cq}{hCE(\psi^q)}

\DeclareMathOperator{\Cqp}{hCE(\nabla  \psi^q)}
\DeclareMathOperator{\Eqp}{hE(\nabla \psi^q)}

\DeclareMathOperator{\ku}{kU}

\DeclareMathOperator{\Map}{Map}

\DeclareMathOperator{\E}{E}

\DeclareMathOperator{\Pm}{P}
\DeclareMathOperator{\Ho}{Ho}
\DeclareMathOperator{\Spe}{Sp}

\newfont{\german}{eufm10}

\newcommand\qu{/\kern-.7ex/}

\newcommand{\iso} {\cong}

\title
{The $p$-local stable Adams conjecture: Erratum}
\author{Prasit Bhattacharya}
\address{Department of Mathematical Science, New Mexico State University, Las Cruces, USA}
\email{prasit@nmsu.edu}

\author{Nitu Kitchloo}
\address{Department of Mathematics, Johns Hopkins University, Baltimore, USA}
\email{nitu@math.jhu.edu}

\date{\today}
\begin{document}


{\abstract 
In this article we prove the stable $p$-local Adams conjecture which has remained persistently open. As part of this article, we correct the sequelae of incorrect results in the literature: \cite{F}, \cite{BhK}. This article is self-contained and serves as an efficient proxy for the topological and algebraic methods pertaining to the circle of ideas surrounding the stable Adams conjecture.}

\maketitle
\begin{center}
{\it In memory of Jack Morava}
\end{center}
\tableofcontents

\section{Discussion of errors and Introduction.}

\medskip
\noindent
Let us briefly recalling the stable Adams conjecture. Let $\underline{\ku}_{(p)}$ be defined via the pullback diagram of spectra
\[ \xymatrix{
	\underline{\ku}_{(p)} \ar[d]\ar[r] & \Hh(\Z) \ar[d] \\
	\ku_{(p)} \ar[r]& \Hh(\Z_{(p)})}
\]
where the bottom horizontal map is the Postnikov approximation of $p$-local connective K-theory $\ku_{(p)}$, and the right vertical map is induced by the canonical inclusion. Similarly, let $\pic^{ev}  {\bf S}_{(p)}$ be the even picard spectrum of the $p$-local sphere ${\bf S}_{(p)}$, defined via the pullback 
\[ \xymatrix{
	\pic^{ev} {\bf S}_{(p)} \ar[d]\ar[r] & \Hh(\Z) \ar[d]^{\times 2} \\
	\pic  {\bf S}_{(p)}  \ar[r]& \Hh(\Z)}
\]
where, as before, the bottom horizontal map is the Postnikov approximation of the Picard spectrum of the $p$-local sphere spectrum and the right vertical map is induced by multiplication by $2$. The following theorem is what we mean by the $p$-local stable Adams conjecture:

\bigskip
\begin{thm*} (see Theorem \ref{PLAC})
	
	\noindent
	Fix $p,q$ to be any primes such that $p \neq q$. Recall the classical $p$-local $J$-homomorphism: 
	\[ J : \Z \times \BUn_{(p)} \longrightarrow \Pic^{ev} {\bf S}_{(p)}, \]
	where $\Pic^{ev} {\bf S}_{(p)}$ is the even Picard space defined as $\Omega^{\infty} \pic^{ev} {\bf S}_{(p)}$. Then $J$ lifts to a stable map 
	\[ \underline{J} : \underline{\ku}_{(p)} \longrightarrow \pic^{ev}  {\bf S}_{(p)} \]
	such that $\underline{J}$ is invariant under precomposition with the Adams operation $\psi^q$. More precisely, $\underline{J}$ admits a  canonical factorization through a map $\underline{J}_q$: 
	\[
	\xymatrix{
		& \ar[dl]  \underline{\ku}_{(p)}   \ar[d]^{\underline{J}}  \\
		\underline{\ku}_{h\psi^q} \ar[r]^{\underline{J}_q} & \pic^{ev} {\bf S}_{(p)}}
	\]
	where $\underline{\ku}_{h\psi^q}$ denotes the homotopy co-equalizer (or homotopy orbits) of the $\psi^q$-action on $\underline{\ku}_{(p)}$. Furthermore, the map $\underline{J}_q$ sends a generator of $\pi_1(\underline{\ku}_{h\psi^q}) = \Z$ to the $q \in \Z_{(p)}^{\times} = \pi_1(\pic^{ev} {\bf S}_{(p)})$. 
\end{thm*}

\medskip
\noindent
A $p$-complete version of this conjecture was claimed by E. Friedlander in \cite{F}, which was the main reference for \cite{BhK}. The project \cite{BhK} was intended to shed light on higher associative structures on Moore spectra, which was the topic of Bhattacharya's thesis \cite{Bh}. So it came to us as a surprise when a referee pointed out that the $p$-completed stable Adams conjecture, as stated in \cite{F}, was false (see \cite{BhK} for details). We have subsequently learnt from E. Friedlander that his error may be traced back to the oversight regarding the fact that the Adams operations do not preserve the orientation of the universal spherical fibrations. 

\medskip
\noindent
Following the referee report for \cite{BhK}, the authors made an attempt to fill the gap left by the error in \cite{F} by supplying a correct proof of the $p$-local stable Adams conjecture. Unfortunately, the proof still had two gaps including one of the authors' own making (see \cite{BhK2} for details), that were later flagged by E. Friedlander in \cite{F3}. 
In this article we remedy both these errors by elaborating on an alternate method of proof that had already been sketched in the appendix to \cite{BhK}. This alternate method relies on a classification result by Friedlander (\cite{F}, Theorem 6.1). In order to make this article independent of \cite{F}, we have included an appendix in this article that supplies a shorter proof of Friedlander's classification result after making appropriate changes to the original statement of the classification result as indicated in the appendix. Given this erratum, the application to higher associative structures on Moore spectra proceeds as described in \cite{BhK}. 

\medskip
\noindent
We start by outlining the errors in \cite{BhK} and describing the proposed changes.

\bigskip
\noindent
\begin{Error}
	In section 4 (Construction 4.1) of \cite{BhK} we define a simplicial scheme, $\mbox{S}^{2i}_\Z$ whose complex points (with the analytic topology) is equivalent to the $2i$-sphere. The basepoints were incorrectly stated as the initial scheme $\mbox{Spec}(0)$. They should have been stated as the terminal scheme $\mbox{Spec}(\Z)$ to make this object well defined. Since all subsequent constructions use the correct definition, this error is not consequential. 
	However, in Remark 4.4, we posit the existence of a map
	\[ \tau_{i,j} : \mbox{S}^{2i}_{\Z} \times \mbox{S}^{2j}_{\Z} \longrightarrow \mbox{S}^{2i+2j}_\Z, \]
	which descends to the smash product map on complex points. Unfortunately, on closer study, we see that such a map is not well-defined. 
\end{Error}

\medskip
\noindent
The map $\tau_{i,j}$ first appears in past work by Friedlander (see for instance \cite{F}, page 139). However, the fact that $\tau_{i,j}$ is not well-defined had been overlooked in the literature. We too failed to notice this in \cite{BhK}, until the error was finally flagged by Friedlander in \cite{F3}. This is a consequential issue, since the map $\tau_{i,j}$ forms the basis of the monoidal structure one requires of all constructions that follow. 

\bigskip
\noindent
{\bf Solution to Error 1}: A obvious way to make $\tau_{i,j}$ manifestly well-defined is to replace the scheme $\mbox{S}^{2n}_\Z$ by the scheme $\mathcal{S}^n$ defined as the $n$-fold smash product of the scheme $\mbox{S}^{2}_\Z$. However, $\mathcal{S}^n$ will have a smaller group of symmetries than we require. More precisely, instead of supporting the action of the full general linear group $\mbox{GL}_n$, $\mathcal{S}^n$ will only admit an action by the subgroup: $\Sigma_n \ltimes \mbox{GL}_1^{\times n} \subseteq \mbox{GL}_n$ (i.e. the normalizer of the standard maximal torus in $\mbox{GL}_n$) . To preserve all the symmetries, one is required to alter the definition of $\mbox{S}^{2n}_\Z$ without changing the homotopy type of its $\C$-points in the analytic topology. This new definition, with the full group of symmetries, has been provided by E. Friedlander in his recent update \cite{F3}. Below we will define two simplicial schemes $\mathcal{D}^n$ and $\mathcal{S}^n$ both of which model the $2n$-sphere (see definition \ref{SGlscheme} and remark \ref{base}). $\mathcal{D}^n$ is the model recently introduced by Friedlander in \cite{F3}, while $\mathcal{S}^n$ is the $n$-fold smash product of the scheme $\mbox{S}^{2}_\Z$ as before. We then relate $\mathcal{D}^n$ and $\mathcal{S}^n$ by a zigzag interpolated by the $n$-fold smash product of $\mathcal{D}^1$. These zigzags are fed into our construction, which is also described in detail below. 

\bigskip
\noindent
\begin{Error}
	Having corrected the map $\tau_{i,j}$ above, one may extend it to a $\mbox{GL}_i \times \mbox{GL}_j$-equivariant map of schemes, where $\mbox{GL}_i$ denotes the general linear group scheme. This allows us to construct a map of bar constructions of the $\mbox{GL}_i$-action on the above simplicial schemes (see \cite{BhK}, diagram 4.26):
	\[ \omega_{i,j} : \mbox{SBGL}_i \times \mbox{SBGL}_j \longrightarrow \mbox{SBGL}_{i+j}. \]
	Changing base to $\C$, and applying the (p-completed) \'etale homotopy type functor to $\omega_{i,j}$ one obtains a map 
	\[ \mbox{\'Et}_p (\mbox{SBGL}_i \times \mbox{SBGL}_j) \longrightarrow \mbox{\'Et}_p (\mbox{SBGL}_{i+j}). \]
	Using the classical comparison theorem (\cite{F2}, Theorem 8.4) one has a weak equivalence
	\[ F : \mbox{\'Et}_p (\mbox{SBGL}_i \times \mbox{SBGL}_j) \llra{\sim} (\mbox{SBGL}^{top}(\C)_i)\hat{\, }_p \times (\mbox{SBGL}^{top}(\C)_j)\hat{\, }_p, \]
	where $\mbox{SBGL}^{top}(\C)_i$ denotes the complex points of $\mbox{SBGL}_i$ with the analytic topology. In \cite{BhK} (see diagram 4.25), we implicitly use the inverse equivalence of $F$ and the map $\omega_{i,j}$ to construct a multicategory. However, the above equivalence $F$ does not have a functorial inverse! This is a consequential error since one requires functoriality in subsequent categorical constructions. 
\end{Error}

\bigskip
\noindent
{\bf Solution to Error 2}: Since $F$ fails to have a functorial inverse,  our method as it currently stands of using multicategories to prove the $p$-local stable Adams conjecture (see \cite{BhK}, Section 5) is not viable. Instead, we are led naturally to the structure of a $\Gamma$-space, since such directed equivalences are precisely the ``specialty" condition in the theory of $\Gamma$-spaces. In fact, a proof of the $p$-local stable Adams conjecture using $\Gamma$-spaces was already sketched as an alternate method of proof in the appendix to \cite{BhK}, and aligns with the method used in \cite{F}. We have therefore chosen to fill in the details in the alternate proof using $\Gamma$-spaces, instead of trying to salvage our earlier proof using multicategories. Our current line of reasoning requires a classification result by Friedlander for sectioned $\Gamma$-fibrations, which is reproved in the appendix.

\medskip
\noindent
Even though we have taken pains to make the content of this article as transparent and readable as possible, some parts remain unavoidably dense. It is for that reason, a general description of the key ideas and steps would be helpful before we begin. We summarize these steps below.

\bigskip
\noindent
{\bf{Summary of the key ideas in this article:}}

\noindent
The key step in this article is the construction (theorem \ref{main}) of a $p$-local spherical fibration of $\Gamma$-spaces $\pi : \mathscr{E} \longrightarrow \mathscr{B}$, a notion that is defined in \ref{fib}. This fibration is built from $p$-localizations of the fiberwise compactificatified universal $\GL_n(\C)$-vector bundle. Moreover, by construction, our fibration $\pi$ comes endowed with an endomorphism that represents the $p$-local Adams operation $\psi^q$ for $q$ any prime coprime to $p$. Friedlander's classification theorem \ref{infinityfib} now allows us to classify $\pi$ via a map of $\Gamma$-spaces 
\[ \mathfrak{J} : \mathscr{B} \longrightarrow \BGx_{\mathbb{S}_{(p)}}. \]
The classificatiion theorem implies that the map $\mathfrak{J}$ is invariant under $\psi^q$. Now, by a result of Segal (theorem \ref{Segal}), the category of $\Gamma$-spaces admits a functor to the category of spectra under which the map $\mathfrak{J}$ becomes the $p$-localized stable $J$ homomorphism. Applying this result to $\mathfrak{J}$ allows us to conclude the $p$-local stable Adams conjecture (theorem \ref{PLAC}). 

\smallskip
\noindent
Hence, the crucial step in this article is the construction of the $\psi^q$-equivariant fibration $\pi : \mathscr{E} \longrightarrow \mathscr{B}$ given in theorem \ref{main}. As in \cite{BhK}, this fibration is constructed by first deconstructing the $\Gamma$-spaces $\mathscr{E}$ and $\mathscr{B}$ along an arithmetic fracture square (definition \ref{AFS}), and glueing the $p$-completed form of $\pi$ to its rational form, along the adelic form. This glueing is done in a $\psi^q$-equivariant fashion. The details are summarized as follows:

\smallskip
\noindent
The $p$-completed form of $\mathscr{E}$ is a $\Gamma$-space built from the \'etale homotopy type (definition \ref{EHTb}) of a simplicial complex scheme $\mathcal{D}^n$ (definition \ref{SGlscheme}) which comes endowed with Galois symmetries (definition \ref{Galois}). These symmetries allow us to construct the $p$-completed form of $\psi^q$ as the automorphism of $\C$ whose inverse extends the Frobenius automorphism of the Witt vectors on $\overline{\F}_q$. We define the rational form of $\mathscr{E}$ to be the $\Gamma$-space built from actions of the Lie group $\Sigma_n \ltimes \GL_1(\C)^{\times n}$ acting in the canonical fashion on the $n$-fold smash product $S\C^{\times} \wedge \ldots \wedge S\C^\times$, where $S\C^{\times}$ denotes the unreduced suspension of $\C^\times$. This smash product is the space of complex points of a simplicial scheme $\mathcal{S}^n$ (definition \ref{SGlscheme}) in the analytic topology.  The rational form of $\psi^q$ acting on $S\C^{\times} \wedge \ldots \wedge S\C^\times$ is induced but the $q$-power self map of $\C^\times$ on each factor. 
At this point, a subtle matter arises of making these two (seemingly distinct) forms of $\psi^q$ described above compatible along the adelic form. This compatibility comes down to relating the two $\Gamma$-spaces described above through a $\psi^q$-equivariant zigzag of $\Gamma$-spaces (equation \ref{zigzag}). This zigzag is obtained by glueing two smaller zigzags. This glueing process requires an important theorem \ref{ptstoet} that shows that the functor of points maps naturally to the \'etale homotopy type of a simplicial complex scheme. Having glued the two smaller zigzags, the arithmetic fracture square is then constructed by functorially straightening the consolidated zigzag (equation \ref{kappa}). The two smaller zigzags are described below:

\smallskip
\noindent
In the first of the two zigzags used in equation \ref{zigzag}, we relate the $\Gamma$-space built from the complex points of $\mathcal{S}^n$ in the analytic topology with the $\Gamma$-space built from the complex points of $\mathcal{S}^n$ in the discrete topology in a $\psi^q$-equivariant fashion (theorem \ref{zig-zagpts}). This zigzag interpolates the $\overline{\F}_q$-points of $\mathcal{S}^n$, on which the two forms of $\psi^q$ agree with the Frobenius automorphism (remark \ref{FL}). 
In the second of the two zigzags used in equation \ref{zigzag}, we relate the $\Gamma$-space built from the \'etale homotopy type of $\mathcal{D}^n$ to the $\Gamma$-space built from the complex points of $\mathcal{S}^n$ in the discrete topology in a fashion that is compatible with respect to $\psi^q$. In this zigzag, we first pass through the $\Gamma$-space built from the \'etale homotopy type of $\mathcal{S}^n$ (theorem \ref{special}) which is then related to the $\Gamma$-space built from the (discrete) complex points of $\mathcal{S}$. 

\smallskip
\noindent
{\bf{General remarks about the article:}}

\noindent
The general structure of this article is as follows. 
The first four sections of this paper provide the relevant background information. In Section \ref{BFS} we review the theory of $\Gamma$-spaces. In Section \ref{ToPC} we introduce two important definitions: $\Gamma$-spaces with rigid sections, and fibrations of $\Gamma$-spaces.  Sections \ref{ToPC} and \ref{APC} also offer important examples including ones arising from \'etale homotopy theory, and from the functor of points. Finally, in Section \ref{FS} we glue several of these examples together along a fracture square to construct a sectioned $p$-local spherical fibration, which we classify, leading to a proof of the $p$-local Adams conjecture. For completeness, we reprove Friedlander's classification result in the appendix. 

\medskip
\noindent
We would also like to point out that in this article we have avoided the language of profinite completions of our \'etale types that was used in \cite{BhK}. Avoiding this language allows us to relate our constructions directly to the constructions made in the canonical reference on \'etale homotopy theory \cite{F2}. Moreover, in \cite{BhK} the precise relation between the Bousfield-Kan $p$-completion of the Artin-Mazur \'etale homotopy type of a simplicial scheme and that of its profinite completion was left unclear. 

\smallskip
\noindent
{\bf{Acknowledgements:}}

\noindent
As should be evident to the reader of this article, all our methods are directly or indirectly motivated by the work of Eric Friedlander and we owe him a profound debt of gratitude, not just for the elegance of his ideas, but also for identifying the oversights in \cite{F} and \cite{BhK} and offering helpful suggestions in \cite{F3}. His work was influential in proving the classical Adams conjecture \cite{F4}, even before his work on the stable Adams conjecture. All of this is to say that the authors continue to believe that the stable Adams conjecture ($p$-local or otherwise) should be attributed to Eric Friedlander.  

\smallskip
\noindent
The second author is grateful to David Gepner and Emily Riehl for patiently helping him understand various categorical aspects pertaining to the material presented in the appendix. He would also like to thank 
Gereon Quick for sharing his quick thoughts on the \'etale homotopy type. But most importantly, we thank the referees of this article for their thoughtful suggestions that helped make this article so much more transparent. 

\section{Review of $\Gamma$-spaces after Bousfield, Friedlander and Segal.} \label{BFS}

\medskip
\noindent
We now define the objects of interest beginning with the definition of $\Gamma$-spaces. Details can be found in \cite{BF}. In staying faithful to the context used in \cite{BF} and \cite{F}, we will work in the classical model category of simplicial sets, which we refer to as ``spaces". We reserve the term ``topological spaces" for objects in the category of compactly generated weak Hausdorff spaces (CGWH) endowed with the classical monoidal model structure.  All categorical constructions like taking products or mapping spaces of topological spaces is understood to be done in CGWH. 

\smallskip
\noindent
At various points in this article  (that are clearly indicated),  we will use the geometric realization functor denoted as $|\quad|$ to move from the category to spaces to the category of topological spaces. To return to the original category of spaces, we will use the singular simplices functor which is denoted as  $\mbox{Sing}$.

\medskip
\begin{defn} \label{Gamma} ($\Gamma$-space)
	
	\noindent
	A $\Gamma$-space is a functor $\mathscr{B} : \mathscr{F} \longrightarrow \sSet_*$
	from the category $\mathscr{F}$ of finite pointed sets to the category $\sSet_*$ of pointed simplicial sets. Moreover, we demand that $\mathscr{B}$ sends each singleton to the constant point simplicial set. We denote the category of $\Gamma$-spaces $\mathscr{F}[\sSet_*]$ to be the category of such functors. Note that it is sufficient to define a $\Gamma$-space on the full sub-category generated by the objects ${\bf{n} }= \{0,1,\ldots, n\}$ (with the element $0$  the basepoint). Henceforth, we shall consider all $\Gamma$-spaces restricted to this sub-category, unless evident otherwise. The category $\mathscr{F}[\sSet_*]$ has the structure of a closed simplicial model category (\cite{BF}, Theorem 3.5), with weak equivalences  defined object wise\footnote{In \cite{BF}, two different model structures are described: the ``strict" and the ``stable" model structures. We use the strict one in this article.}.  
	
	\smallskip
	\noindent
	A $\Gamma$-space $\mathscr{B}$ is called``special " if it satisfies the property that the following map of simplicial sets is a weak equivalence: 
	\[ \prod_{i=1}^n d^{\mathscr{B}}_i : \mathscr{B}({\bf n}) \longrightarrow \prod_{i=1}^n \mathscr{B}({\bf 1}), \]
	where $d_i : {\bf n } \longrightarrow {\bf 1}$ denote the various projections that send $i \in {\bf n}$ to $1$, and $j \neq i$ to $0$. 
	
	\smallskip
	\noindent
	Note: In practice one may allow (as several authors do), the value of a $\Gamma$-space $\mathscr{F}$ on the singleton $\ast$ to be a contractible simplicial set. Then, normalizing the value $\mathscr{F}(S)$ for any pointed $S$ by taking the quotient $\mathscr{F}(S)/\mathscr{F}(\ast)$, reduces to the above definition.
\end{defn}

\medskip
\begin{remark} \label{functoriality}
	Given a functor $\mbox{F} : \sSet_* \longrightarrow \sSet_*$, we may (objectwise) compose a $\Gamma$-space $\mathscr{B}$ with $\mbox{F}$ to get a new $\Gamma$-space $\mbox{F}\mathscr{B}$. The functor of particular interest to us is Bousfield-Kan $p$-completion \cite{BK}. We will also work in the algebraic context (see Section \ref{APC}), where the functors of interest will be the \'etale homotopy type and the functor of points. 
\end{remark}

\medskip
\begin{example} \label{exampleN} (The $\Gamma$-space $\mathscr{N}$) 
	
	\noindent
	Let the set of non-negative numbers $\N$ be seen as a discrete pointed simplicial set, with $0$  the basepoint. Define a $\Gamma$-space $\mathscr{N}$ by demanding that $\mathscr{N}({\bf n}) = \N^{\times n}$ for $n>0$ with the zero $n$-tuple  the basepoint. Given $d : {\bf m} \longrightarrow {\bf n}$ in $\mathscr{F}$, we define:
	\[ d(I) = d(i_1, \ldots, i_m) = (j_1, \ldots, j_n) = J, \quad \mbox{where} \quad j_s = \sum_{t \in d^{-1}(s)} i_t, \quad \mbox{if} \quad d^{-1}(s) \neq \emptyset,\]
	and the element $j_s $ is defined to be $0$ if $d^{-1}(s)$ is empty. 
\end{example}

\begin{defn} \label{F/N} (The category $\mathscr{F}[\sSet_*]/\mathscr{N}$)
	
	\noindent
	The $\Gamma$-spaces we consider in this article are naturally objects in $\mathscr{F}[\sSet_*]/\mathscr{N}$, the category over $\mathscr{N}$. Given an n-tuple of non-negative integers $I = (i_1, \ldots, i_n) \in \N^{\times n}$, we may define $\mathscr{B}_I({\bf n})$ to be the pre-image of $I$ under the map $\mathscr{B}({\bf n}) \longrightarrow \mathscr{N}({\bf n})$. In particular, we have a decomposition:
	\[ \mathscr{B}({\bf n}) = \coprod_{I \in \N^{\times n}} \mathscr{B}_I({\bf n}). \]
	Since ${\bf n}$ is implicit in the index $I$ , we will often denote $\mathscr{B}_I({\bf n})$ by just $\mathscr{B}_I$. 
\end{defn}

\noindent
We now recall some basic definitions from \cite{BF}:

\medskip
\begin{defn} \label{sSpectra} (Simplicial spectrum: \cite{BF}, Definition 2.1)
	
	\noindent
	Let $\Sigma$ be the pointed simplicial set representing the circle, endowed with one zero simplex and one non-degenerate one simplex. A spectrum valued in simplicial sets is a collection of pointed simplicial sets $\mathscr{S} := \{ \mathscr{S}_n, n \geq 0 \}$ endowed with structure maps $\Sigma \wedge \mathscr{S}_n \longrightarrow \mathscr{S}_{n+1}$. Such spectra form a closed simplicial model category (\cite{BF}, Theorem 2.3). 
\end{defn}

\medskip
\begin{defn} \label{Q} ($\Omega$-spectra: \cite{BF}, Section 2)
	
	\noindent
	Note that on applying the geometric realization functor, $| \quad |$, to the structure maps of a spectrum $\{ \mathscr{S}_n  \}$ gives rise to a (topological) pre-spectrum $\{ | \mathscr{S}_n | \}$. We define $\mathscr{S}$ to be an $\Omega$-spectrum if the underlying topological pre-spectrum is an $\Omega$-spectrum.

	\medskip
	\noindent
	There exists a functor $\mbox{Q}$ on spectra taking values in $\Omega$-spectra, and receiving a natural weak equivalence from the identity functor. Let $\mbox{Sing}$ denote the singular simplices functor. Then the value of the functor $\mbox{Q}$ on a spectrum $\mathscr{S} = \{ \mathscr{S}_n \}$ is defined to be the spectrum $\mbox{Q}(\mathscr{S})$ with 
	\[ \mbox{Q}(\mathscr{S})_n := \colim_i \mbox{Sing} \,  \Omega^i |\mathscr{S}_{n+i}|. \]
\end{defn}

\medskip
\begin{defn} \label{nabla} (The functor $\nabla$: \cite{BF}, Sections 4,5\footnote{$\nabla$ is our notation and not the notation used in \cite{BF}.})
	
	\noindent
	Let $\mathscr{B}$ be a $\Gamma$-space, and let $\Sigma$ be as in definition \ref{sSpectra}. Set $\Sigma^0 := \Delta[0] \cup \ast$ and $\Sigma^n := \Sigma^{\wedge n}$. Define a simplicial set $\mathscr{B}(\Sigma^n)$ as follows. First observe that for a fixed $k$, evaluating the $\Gamma$-space $\mathscr{B}$ on the finite set $\Sigma^n_k$ of $k$-simplices of $\Sigma^n$ gives rise to simplicial set $\mathscr{B}(\Sigma^n_k)$. As $k$ varies, the collection $\mathscr{B}(\Sigma^n_k)$ defines a bi-simplicial set; one defines $\mathscr{B}(\Sigma^n)$ as the diagonal of this bi-simplicial set. 
	
	\medskip
	\noindent
	Now define a (simplicial) spectrum $\nabla(\mathscr{B})$ with $\{ \nabla(\mathscr{B})_n := \mathscr{B}(\Sigma^n), \, \,  n \geq 0 \}$.  
	The structure map $\Sigma \wedge \mathscr{B}(\Sigma^n) \longrightarrow \mathscr{B}(\Sigma^{n+1})$ is induced by applying $\mathscr{B}$ to the maps of finite sets parametrized by the simplices of $\Sigma$, namely: $\{ x \}  \times \Sigma^n_k \longrightarrow \Sigma^{n+1}_k$ for $x \in \Sigma_k$ and $k \geq 0$. 
\end{defn}

\smallskip
\noindent
For any $k \geq 0$, let us identify $\Sigma_k^0$ with the pointed set ${\bf 1}$. Given a $\Gamma$-space $\mathscr{B}$, we have a zigzag
\[ \mathscr{B}(\Sigma^0) \times \mathscr{B}(\Sigma^0) \llla{\pi_1 \times \pi_2} \mathscr{B}({\bf 2}) \llra{\mu} \mathscr{B}(\Sigma^0),\]
where $\pi_1$ (resp. $\pi_2$) denotes the map induced by the projection ${\bf 2} \rightarrow {\bf 1}$, that sends the element $2$ (resp. $1$) to the basepoint and $1$ (resp. $2$) to $1$. Similarly, $\mu$ denotes the map induced by the map ${\bf 2} \rightarrow {\bf 1}$ that sends both nonzero elements to $1$. If $\mathscr{B}$ is special (see definition \ref{Gamma}), then the map $\pi_1 \times \pi_2$ is a weak equivalence, making $|\nabla(\mathscr{B})_0|$ an abelian topological monoid up to coherent homotopy. 

\medskip
\noindent
The following fundamental result is a restatement of a result due to Graeme Segal \cite{Se}:

\medskip
\begin{thm} (\cite{BF}, Theorem 4.4; \cite{Se}, Proposition 1.4) \label{Segal} The identification $\mathscr{B} \mapsto \nabla(\mathscr{B})$ defined above defines a functor 
	\[ \nabla : \Ho \mathscr{F}[\sSet_*] \longrightarrow \Ho (\Spe) \]
	from the homotopy category of $\Gamma$-spaces, to the homotopy category of Spectra. Furthermore, if $\mathscr{B}$ is special, then the natural map
	\[ \nabla(\mathscr{B})_n \longrightarrow \mbox{Q}(\nabla(\mathscr{B}))_n\]
	is a weak equivalence of spaces for $n > 0$. Moreover, $ |\nabla(\mathscr{B})_1|$ is equivalent to the classifying space \cite{M1} of the topological monoid  $|\nabla(\mathscr{B})_0|$. 
\end{thm}

\section{Graded topological permutative categories and  $\Gamma$-spaces with rigid sections.} \label{ToPC}

\medskip
\noindent
In this section we describe a functorial construction by P. May \cite{M} of a $\Gamma$-space, starting with a topological permutative category. We begin by recalling the definition.

\medskip
\begin{defn} \label{PC} (Topological Permutative Category: \cite{M})

	\noindent
	\noindent
	A Topological Permutative Category $(\mathscr{A}, \otimes, \ast, \sigma)$ is the following data. 
	
	\smallskip
	\noindent
	{\bf (1)} A pointed topological object space $\mbox{Ob}(\mathscr{A})$ with non-degenerate basepoint $\ast$, and a morphism space $\mbox{Mor}(\mathscr{A})$. Moreover, the identity morphism $\mbox{Ob}(\mathscr{A}) \longrightarrow \mbox{Mor}(\mathscr{A})$ is a cofibration. 
	
	\smallskip
	\noindent
	{\bf (2)} A strictly associative monoidal structure $\otimes$ that is represented by a continuous map of pointed object spaces, and compatible maps on morphisms spaces, and so that $\ast$ is a strict unit. 
	
	\smallskip
	\noindent
	{\bf (3)} For any pair of objects $(a,b)$, a continuous, natural isomorphism $\sigma(a,b) \in \mbox{Mor}(a \otimes b, b \otimes a)$ that satisfies the relations (our compositions are read left to right): 
	\[ \sigma(\ast, a) = \sigma(a,\ast) = \mbox{Id}_a, \quad \sigma(a,b) \sigma(b,a) = \mbox{Id}_{a\otimes b}, \quad \sigma(a \otimes b, c) = (\mbox{Id}_a \otimes \sigma(b,c)) \, (\sigma(a,c) \otimes \mbox{Id}_b). \]
	A functor of topological permutative categories is a continuous functor defined in the obvious way. 
	
	\medskip
	\noindent
	Note: The continuity of $\sigma$ can be expressed as a continuous lift:
	\[
	\xymatrix{
		& \mbox{Mor}(\mathscr{A}) \ar[d]^{s,t} \\
		\mbox{Ob}(\mathscr{A}) \times \mbox{Ob}(\mathscr{A}) \ar[r]^{\, \, \otimes \times \otimes(\tau) \, \, \,\,  } \ar[ur]^{\sigma} & \mbox{Ob}(\mathscr{A}) \times \mbox{Ob}(\mathscr{A}) \,}
	\]
	where $s,t$ denote the source and target maps and $\otimes(\tau)$ denotes monoidal structure applied to the twist map $\tau$.
	In terms of the above lift, naturality of $\sigma$ can then be expressed by a commutative diagram:
	\[ 
	\xymatrix{  \mbox{Mor}(\mathscr{A}) \times \mbox{Mor}(\mathscr{A}) \ar[d]^{s_1 \times s_2 \times \otimes(\tau)} \ar[r]^{\otimes \times t_1 \times t_2 \quad \quad } &  \mbox{Mor}(\mathscr{A}) \times \mbox{Ob}(\mathscr{A}) \times \mbox{Ob}(\mathscr{A})  \ar[d]^{ (\underline{\, \, \,}) \, \sigma} \\
		\mbox{Ob}(\mathscr{A}) \times \mbox{Ob}(\mathscr{A}) \times \mbox{Mor}(\mathscr{A})	\ar[r]^{\quad \quad \sigma \, (\underline{\,\, \,})} & \mbox{Mor}(\mathscr{A}),}
	\]
	where $s_1,s_2$ denote the source maps of respective factors of $\mbox{Mor}(\mathscr{A}) \times \mbox{Mor}(\mathscr{A})$ (and simililarly for $t_1,t_2$). The relations satisfied by $\sigma$ given in point {\bf (3)} may be similarly expressed as appropriate commuting diagrams. 
\end{defn}

\noindent
In analogy with a $\Gamma$-space, we may define a $\Gamma$-topological permutative category to be a functor ${\bf n} \mapsto \mathscr{A}(\bf n)$ from the category of pointed finite sets to topological permutative categories, so that $\mathscr{A}(\bf 0)$ is the trivial category. 
In \cite{M}, May defines a functorial construction described below that takes a topological permutative category $\mathscr{A}$ and gives rise to a $\Gamma$-topological permutative category ${\bf n} \mapsto \mathscr{A}(\bf n)$.

\medskip
\begin{thm} \cite{M} \label{perm} There exists a functor $\mathscr{A} \mapsto \{ \bf n \mapsto \mathscr{A}(\bf n)\}$ from the category of topological permutative categories to the category of $\Gamma$-topological permutative categories. Moreover, $\mathscr{A}(\bf 1)$ is the original category $\mathscr{A}$, and for $n > 1$, the topological category $\mathscr{A}(\bf n)$ is equivalent to the category $\mathscr{A}^{\times n}$ (see remark \ref{split} below). Taking the singular nerve $\BA(\bf n)$ (see definition \ref{CS} below) of each category $\mathscr{A}(\bf n)$ gives rise to a special $\Gamma$-space $\BA$, defined as $\bf n \mapsto \BA(\bf n)$. 
\end{thm}
\medskip

\begin{remark} \label{E-infinity}
	Given a permutative category $\mathscr{A}$, we may consider the underlying
	symmetric monoidal category $\mathscr{A}$. It is well known that the nerve of $\mathscr{A}$ (see definition \ref{CS}) is an $E_{\infty}$-space (\cite{Se}, Section 2). By \cite{M} (Theorem 3), this $E_\infty$-structure is naturally equivalent to the one obtained via theorem \ref{Segal} applied to the $\Gamma$-space obtained in theorem \ref{perm} above. This observation will allow us to identify the $E_\infty$-structure of several examples that are to follow. 
	
\end{remark}

\medskip
\begin{Const}\cite{M} \label{ConstM}
	We briefly recall the construction of the $\Gamma$-topological permutative category given in \cite{M}. The object space of the topological permutative category $\mathscr{A}(\bf n)$ has elements given by a collection of elements $a_S \in \mbox{Ob}(\mathscr{A})$ indexed by subsets $S \subseteq {\bf n}$ containing $0$, and isomorphisms $i_{S,T} : a_{S \cup T} \longrightarrow a_S \otimes a_T$ where $S, T \subseteq {\bf n}$ are any pointed subsets so that $S \cap T = \{ 0 \}$. Moreover, one demands that $a_{\{ 0 \}} = \ast$ and that the following diagrams commute:
	\[
	\xymatrix{
		& a_S \otimes \ast \ar[d]^{=} &  a_{R \cup S \cup T} \ar[d]^{i_{R \cup S, T}} \ar[r]^{i_{R, S \cup T}} & a_R \otimes a_{S \cup T} \ar[d]^{i_{S,T}} & a_{R \cup S} \ar[d]^{\Id} \ar[r]^{i_{R, S}} & a_R \otimes a_S \ar[d]^{\sigma} \\
		a_S \ar[r]^{\Id} \ar[ur]^{i_{S,\{0\}}} & a_S & a_{R \cup S} \otimes a_T \ar[r]^{i_{R \cup S}}  & a_R \otimes a_S \otimes a_T & a_{S \cup R} \ar[r]^{i_{S, R}} & a_S \otimes a_R.}
	\]
	Morphisms in $\mathscr{A}(\bf n)$ between $(a_S, i_{S,T})$ and $(b_S, j_{S,T})$ are given by a collection of morphisms in $\mathscr{A}$ $\{ \alpha_S : a_S \longrightarrow b_S \}$ with $\alpha_{\{ 0 \}} = \Id_\ast$, making the following diagrams commute
	\[ 
	\xymatrix{a_{S \cup T} \ar[d]^{\alpha_{S \cup T}} \ar[r]^{i_{S, T}} & a_S \otimes a_T \ar[d]^{\alpha_S \otimes \alpha_T} \\
		b_{S \cup T} \ar[r]^{j_{S, T}} & b_S \otimes b_T.}
	\]
	The monoidal structure $\otimes$ and commutative isomorphism $\sigma$ for $\mathscr{A}(\bf n)$ is defined pointwise for indexing collection $(S, \{S,T \})$. The unit in $\mathscr{A}(\bf n)$ is defined to be $(\ast_S, \Id_{S,T})$. The object and morphism spaces of $\mathscr{A}(\bf n)$ are topologized as respective subspaces of 
	\[ \prod_{S} \mbox{Ob}(\mathscr{A}) \quad \times \prod_{\substack{S,T \\ S \cap T = \{ 0 \}}} \mbox{Mor}(\mathscr{A}), \quad \mbox{and} \quad \quad \prod_{S} \mbox{Mor}(\mathscr{A}). \]
	Given a pointed map $\varphi : {\bf m} \longrightarrow {\bf n}$, and subsets $S, T \subseteq \bf n$ containing $0$, let $U, V \subseteq \bf m$ be given by $U = \varphi^{-1}(S-\{0\}) \cup \{0 \}$ and $V = \varphi^{-1}(T-\{0\}) \cup \{0\}$. Define a functor
	\[ \varphi^{\mathscr{A}} : \mathscr{A}({\bf m}) \longrightarrow \mathscr{A}({\bf n}), \quad (\varphi^{\mathscr{A}}(a)_{S}, \varphi^{\mathscr{A}}(i)_{S,T}) = (a_U, i_{U,V}), \quad \varphi^{\mathscr{A}}(\alpha)_S = \alpha_U. \]
	One checks that $\varphi^{\mathscr{A}}$ is a functor of topological permutative categories, making $\bf n \mapsto \mathscr{A}(\bf n)$ a $\Gamma$-permutative category. 
\end{Const}

\begin{remark} \label{split}
	Notice that for an object $\{ a_S, i_{S,T} \}$ of $\mathscr{A}(\bf n)$, the object $a_S$ of $\mathscr{A}$ is isomorphic to the object $a_{s_1} \otimes \ldots \otimes a_{s_r}$ by the iterated use of the morphisms $i_{S,T}$, where $s_1 < s_2 < \cdots < s_r$ are the nonzero elements of $S$ with ordering induced from $\bf n$. Using this identification, one obtains a functor $\mathscr{A}({\bf n}) \longrightarrow \mathscr{A}^{\times n}$. It is shown in \cite{M} that this functor is an equivalence, with the inverse given by the section that sends an $n$-tuple of objects of $\mathscr{A}$: $\{ a_1, \ldots, a_n \}$ to the family $\{ a_S, i_{S,T} \}$ with $a_S = a_{s_1} \otimes \cdots \otimes a_{s_r}$ with $\{ s_j \}$ as above, and $i_{S,T}$ given by a shuffling of $a_S \otimes a_T$ using $\sigma$. 
\end{remark}

\bigskip
\begin{defn} \label{CS} (The Singular Nerve $\BC$)
	
	\noindent
	Let $\mathscr{C}$ be an arbitrary small topological category with a distinguished base-point object $\ast$. Then the pointed simplicial set $\BC$ (the singular nerve) is defined to be the diagonal of the bi-simplicial set obtained by taking level wise singular simplcies on the nerve $\Ne \mathscr{C}$ of $\mathscr{C}$: 
	\[ \BC := \Delta \circ \mbox{Sing} \Ne \mathscr{C}. \]
	As an example, consider the topological category with objects as a topological space $X$ with morphisms generated by the right action of a topological monoid $\G$ on $X$. Then $\Ne (X, \G)$ is equivalent to the Bar construction whose $n$-simplicies are given by the space $X \times \G^{\times n}$. Applying the singular simplices functor gives rise to a bisimplicial set $\mbox{Sing} \Ne (X, \G)$ whose $(m, n)$ simplicies are given by $\mbox{Sing}_m(X \times \G^{\times n})$. This gives $\B(X, \G)$ a concrete description as the simplicial set obtained by taking the diagonal of the singular bar construction. If the unit element of $\G$ is non-degenerate, then the weak homotopy type of the geometric realization of $\B(X, \G)$ is equivalent to the homotopy orbits $X \times_{\G} \EGg$, for the $\G$-action on $X$. If in addition, we are also given a space $Y$ with a left $\G$-action, then $\B(X, \G)$ generalizes to $\B(X, \G, Y)$ obtained as the diagonal of the singular simplices functor applied to the two sided bar construction $\Ne(X, \G, Y)$ with $n$-simplices $X \times \G^{\times n} \times Y$. 
\end{defn}

\medskip
\begin{remark} \label{trisimp}
	A permutative category internal to simplicial topological spaces is the same as a simplicial object in topological permutative categories. For such an $\mathscr{A}$, May's construction (being functorial) gives rise to a simplicial object in $\Gamma$-topological permutative categories $\bf n \mapsto \mathscr{A}(\bf n)$. The nerve $\Ne \mathscr{A}(\bf n)$ is now a bisimplicial topological space. Hence the singular simplices $\mbox{Sing} \Ne \mathscr{A}(\bf n)$ applied levelwise is a trisimplicial set. We define the singular nerve $\BA(\bf n)$ as the (triple) diagonal: 
	\[ \BA({\bf n}) := \Delta \circ \mbox{Sing} \Ne \mathscr{A}(\bf n). \]
	As in the case of topological permutative categories, the result is a special $\Gamma$-space $\BA$. 
\end{remark}

\bigskip
\noindent
All the topological permutative categories we encounter in this article have a special form which we define as graded topological permutative categories. 

\medskip
\begin{defn} \label{TPC} (Graded Topological Permutative Category)
	
	\noindent
	A graded topological permutative category $(\mathscr{A}, \otimes, \ast, \sigma)$ is the following data: 
	
	\medskip
	\noindent
	{\bf (1)} A collection of topological spaces ${\bf{X}}^i$ with non-degenerate basepoints, so that one has a permutative category $\mathscr{A}$ with
	\[ \mbox{Ob}(\mathscr{A}) = \coprod_{i \in \N} {\bf{X}}^i, \quad \mbox{with} \quad {\bf{X}}^0 := \ast. \]
	For $i \geq 0$, there is a topological monoid $\G_i$ with a non-degenerate unit and with $\G_0$ the trivial monoid. Morphisms of $\mathscr{A}$ are generated by the right action of $\G_i$ on the pointed topological spaces ${\bf{X}}^i$. Hence:
	\[ \mbox{Mor}(\mathscr{A}) = \coprod_{i \in \N} ({\bf{X}}^i \times \G_i). \]
	The source and target maps on the space of morphisms are given by the projection map and the (right) action map on ${\bf{X}}^i \times \G_i$ respectively. 
	
	\medskip
	\noindent
	{\bf (2)} One is also given a strictly associative monoidal structure $\otimes$ that is represented by a collection of maps of pointed topological spaces, and compatible continuous homomorphisms on morphisms:
	\[ \otimes : {\bf{X}}^i \times {\bf{X}}^j \longrightarrow {\bf{X}}^{i+j}, \quad \quad \otimes : \G_i \times \G_j \longrightarrow \G_{i+j}, \]
	so that $\ast$ is a strict unit and that the monoidal structure on the objects  collapses the subspaces $\ast_i \times {\bf{X}}^j$ and ${\bf{X}}^i \times \ast_j$ to the basepoint of ${\bf{X}}^{i+j}$, thereby factoring through the smash product ${\bf{X}}^i \wedge {\bf{X}}^j$ where $i,j > 0$ and $\ast_i$ and $\ast_j$ denotes the respective basepoints of ${\bf{X}}^i$ and ${\bf{X}}^j$. 
	
	\medskip
	\noindent
	{\bf (3)} Furthermore, one has a isomorphism $\sigma$ represented by a lift, and with $s,t$ denoting source and target maps defined as projection and action respectively:
	\[
	\xymatrix{
		& {\bf{X}}^{i+j} \times \mbox{G}_{i+j} \ar[d]^{s,t} \\
		{\bf{X}}^i \times {\bf{X}}^j \ar[r]^{\otimes \times \otimes(\tau) \, \, \, \, \,  } \ar[ur]^{\sigma} & {\bf{X}}^{i+j} \times {\bf{X}}^{i+j}\,}
	\]
	so that $\sigma$ is natural (expressed as a commutative diagram as in definition \ref{PC}), and satisfies the relations (our compositions are read left to right): 
	\[ \sigma(\ast, a) = \sigma(a,\ast) = \Id_a, \quad \sigma(a,b) \sigma(b,a) = \Id_{a\otimes b}, \quad \sigma(a \otimes b, c) = (\Id_a \otimes \sigma(b,c)) \, (\sigma(a,c) \otimes \Id_b). \]
\end{defn}

\medskip
\begin{remark}\label{smash}
	In several examples that follow, ${\bf{X}}^i$ for $i > 0$ will in fact be the $i$-fold smash product of some fixed (possibly simplicial) topological space (see examples, \ref{NormalizerPCX}, \ref{CXX}). However, we work in the above generality because the examples we construct in section \ref{APC} are not canonically equal to the $i$-fold smash products. 
\end{remark}

\medskip
\noindent
The structure of the $\Gamma$-topological permutative category $\mathscr{A}(\bf n)$ can be explicitly described in the special case when $\mathscr{A}$ is a graded topological permutative category. We begin with: 

\medskip
\begin{defn}\label{G_I^+}
	
	\noindent
	Given $I = (i_1, \ldots, i_n)$, and a subset $T \subseteq I$, let $\mu_I(T) = \sum_{t \in T} i_t$, and let $\G_{\mu_I(T)}^{\times} \subseteq \G_{\mu_I(T)}$ denote the subgroup of invertible elements. Now define 
	\[ \G_I =  \G^I \times \G_I^{+} 
	\quad \mbox{where} \quad \G^I = \prod_{k=1}^n \G_{i_k} \quad \mbox{and} \quad \G_I^{+} = \prod_{T \subseteq \bf n} \G^{\times}_{\mu_I(T)}  \]
	where the subsets $T$ in the definition of $\G_I^{\times}$ are assumed to belong to the collection of subsets of ${\bf n}$ with at least two nonzero elements. Notice that while $\G_I$ need only be a monoid, $\G_I^{+}$ is a group. 
\end{defn}

\medskip
\begin{lemma} \label{TPC1}
	Given a graded topological permutative category $\mathscr{A}$ as in definition \ref{TPC}, the category $\mathscr{A}(\bf n)$ is a coproduct of topological categories $\mathscr{A}_I(\bf n)$, where $I \in \N^{\times n}$. The object space of $\mathscr{A}_I(\bf n)$ is:
	\[ \mbox{Ob}(\mathscr{A}_I({\bf n})) =  {\bf{X}}_I \quad \mbox{where} \quad {\bf{X}}_I = {\bf{X}}^I \times \G_I^{+}, \quad \mbox{with} \quad {\bf{X}}^I = \prod_{k=1}^n {\bf{X}}^{i_k}. \]
	Furthermore, the morphisms of $\mathscr{A}_I({\bf n})$ are defined by the canonical right action of $\G_I$ on ${\bf{X}}_I$:
	\[ \mbox{Mor}(\mathscr{A}_I({\bf n})) = {\bf{X}}_I \times \G_I,\]
	with source and target map given by the projection on ${\bf{X}}_I$ and the right action of $\G_I$ on ${\bf{X}}_I$ resp. 
\end{lemma}
\begin{proof}
	By the construction \ref{ConstM}, $\mbox{Ob}(\mathscr{A}(\bf n))$ are given by a collection of objects $a_S$ of $\mathscr{A}$ where $S \subseteq {\bf n}$ is a pointed subset, as well as isomorphisms $i_{S,T} : a_{S \cup T} \longrightarrow a_S \otimes a_T$, where $S, T \subseteq {\bf n}$ are any non-trivial pointed subsets so that $S \cap T = \{ 0 \}$. This collection is required to satisfy compatibility conditions as displayed by the diagrams in construction \ref{ConstM}. 
	
	\medskip
	\noindent
	Given a collection $(a_S, i_{S,T})$ as above, and a nonzero element $s \in {\bf n}$, let $x_s \in {\bf{X}}^{i_s}$ denote the value of $a_{\{ 0, s \}}$. Let $I \in \N^{\times n}$ be defined as $I = (i_1, \ldots, i_n)$. 
	Given a pointed subset $S \subseteq {\bf n}$ with at least two nonzero elements, let $\gamma_S : a_S \longrightarrow x_{s_1} \otimes x_{s_2} \otimes \ldots \otimes x_{s_r}$ denote the isomorphism obtained by the repeated use of maps $i_{S,T}$ where $s_1 < s_2 < \cdots < s_r$ are the nonzero elements of $S \subseteq \bf n$. We represent the elements $\gamma^{-1}_S$ in the subgroup of invertible elements $\G_{\mu_I(S)}^{\times}$, as in definition \ref{G_I^+}. Hence we obtain $x_I := (x_i, x^{\times}_S) \in {\bf{X}}^I \times \G_I^{+} = {\bf{X}}_I$, where $x^{\times}_S := \gamma_S^{-1}$. As $I$ varies over $\N^{\times n}$, we obtain a continuous map 
	\[ \mbox{Ob}(\mathscr{A}({\bf n})) \longrightarrow \coprod_{I \in \N^{\times n}} {\bf{X}}_I. \]
	We now observe that one may uniquely recover the collection $(a_S, i_{S,T})$ from $(x_i, x_S^{\times})$. To see this, notice that $a_S = (x_{s_1} \otimes x_{s_2} \otimes \ldots \otimes x_{s_r}) \, \gamma_S^{-1} = (x_{s_1} \otimes x_{s_2} \otimes \ldots \otimes x_{s_r}) \, x_S^{\times}$. Similarly, we may express $i_{S,T}$ in terms of $\gamma_S = (x_S^{\times})^{-1}$ using the following commutative diagram: 
	
	\[ 
	\xymatrix{a_{S \cup T} \ar[d]^{\gamma_{S \cup T}} \ar[r]^{i_{S, T}} & a_S \otimes a_T \ar[d]^{\gamma_S \otimes \gamma_T} \\
		(x_{u_1} \otimes \ldots \otimes x_{u_{r''}} )\ar[r]^{\sigma \quad \quad \quad \quad } & (x_{s_1} \otimes \ldots \otimes x_{s_r}) \otimes (x_{t_1} \otimes \ldots \otimes x_{t_{r'}}) }
	\]
	where $\sigma$ is the permutation isomorphism that rearranges the monoidal product of elements $\{ x_{u_1}, \ldots, x_{u_{r''}} \}$ where $ u_1 < \cdots < u_{r''} $ are the nonzero elements of $S \cup T$.
	
	\smallskip
	\noindent
	We now describe how to identify the morphisms of $\mathscr{A}(\bf n)$. Recall from construction \ref{ConstM} that morphisms between $(a_S, i_{S,T})$ and $(b_S, j_{S,T})$ in $\mbox{Mor}(\mathscr{A}(\bf n))$ are given by a collection of morphismsms $\alpha_S : a_S \longrightarrow b_S$ that are commute with the isomorphisms $i_{S,T}$ as follows
	
	\[ 
	\xymatrix{a_{S \cup T} \ar[d]^{\alpha_{S \cup T}} \ar[r]^{i_{S, T}} & a_S \otimes a_T \ar[d]^{\alpha_S \otimes \alpha_T} \\
		b_{S \cup T} \ar[r]^{j_{S, T}} & b_S \otimes b_T.}
	\]
	The above diagram can be iterated to obtain a diagram in terms of $\gamma_S$:
	\[ 
	\xymatrix{a_S \ar[d]^{\alpha_S} \ar[r]^{\gamma_S \quad \quad} & x_{s_1} \otimes \ldots \otimes x_{s_r} \ar[d]^{\alpha_{s_1} \otimes \cdots \otimes \alpha_{s_r}} \\
		b_S \ar[r]^{\gamma'_S \quad \quad } & y_{s_1} \otimes \ldots \otimes y_{s_r}  .}
	\]
	From the above diagram we notice that the morphisms $\alpha_1, \ldots, \alpha_n$, and the isomorphism $\gamma'_S$ can be chosen arbitrarily to define a morphism $\{\alpha_S\}$ in $\mbox{Mor}(\mathscr{A}(\bf n))$. For any nonzero element $s \in {\bf n}$, we may identify the morphisms $\alpha_s$ with  the action of an element of $g_s \in \G_{i_s}$. Similarly, for any $S \subseteq {\bf n}$ with at least two nonzero elements, $g^{\times}_S := \gamma_S {\gamma'_S}^{-1}$ yields an element of $\G^{\times}_{\mu_I(S)}$. These elements assemble into an element $g_I := (g_i, g^{\times}_S) \in \G^I \times \G_I^{+} = \G_I$ (see definition \ref{G_I^+}). From this it is straightforward to see that if $(a_S, i_{S,T})$ is identified with point $x_I$ in ${\bf{X}}_I$ as above, then $\{\alpha_S\}$ may be identified with the pair $(x_I,g_I) \in {\bf{X}}_I \times \G_I$. 
	
	\medskip
	\noindent
	It follows that $\mathscr{A}(\bf n)$ is a coproduct of topological categories $\mathscr{A}_I(\bf n)$, having objects given by the space ${\bf{X}}_I$, and morphisms given by the action of the monoid $\G_I$. This completes the proof of the lemma. 
\end{proof}

\smallskip
\noindent
We now study the functoriality of the above description under morphisms in $\mathscr{F}$. 

\smallskip
\begin{lemma} \label{retracts} 
	Let $\mathscr{A}$ be a graded topological permutative category. Given $d : \bf m \longrightarrow \bf n$, the functor $d^{\mathscr{A}} : \mathscr{A}({\bf m}) \longrightarrow \mathscr{A}({\bf n})$ restricts to a functor $d^{\mathscr{A}} : \mathscr{A}_I({\bf m}) \longrightarrow \mathscr{A}_J({\bf{n}})$ with $I \in \N^{\times m}$, and $J := d(I) \in \N^{\times n}$ defined as in example \ref{exampleN}. We introduce helpful notation in the next two paragraphs, after which we describe the behavior of the functor $d^{\mathscr{A}} : \mathscr{A}_I({\bf m}) \longrightarrow \mathscr{A}_J({\bf{n}})$. 
	
	\medskip
	\noindent
	{\bf Start of Notation}:
	Let $x_I = (x_t, x^{\times}_T) \in {\bf{X}}_I = {\bf X}^I \times \G_I^{+}$ be an object of $\mathscr{A}_I({\bf m})$, where $x_t$ denotes the components in ${\bf X}^I = \prod_{k=1}^m {\bf X}^{i_k}$ indexed by some nonzero $t \in {\bf m}$, and $x^{\times}_T$ denotes the components in $\G_I^{+} = \prod_{T \subseteq \bf m} \G^{\times}_{\mu_I(T)}$. Given $d : {\bf m} \longrightarrow {\bf n}$, and any $0 \neq s \in \bf n$, define $x_d(s)$ as the object of $\mathscr{A}$ given by the monoidal product $x_{t_1} \otimes \cdots \otimes x_{t_k}$, where $t_1 < t_2 < \cdots < t_k$ are the nonzero elements of $d^{-1}\{ s \}$ (and as the unit $\ast$ if $d^{-1}\{ s\} = \emptyset$). 
	Similarly, given $g_I = (g_t, g^{\times}_T) \in \G_I$, with  $g_I = (g_t, g^{\times}_T)$  indexed as in $x_I$, define $g_d(s)$ as the monoidal product $g_{t_1} \otimes \cdots \otimes g_{t_k}$.
	
	\smallskip
	\noindent
	Given $d$ as above, and a pointed subset $S \subseteq {\bf n}$, define $d^\ast S \subseteq {\bf m}$ as $d^\ast S := d^{-1}(S-\{0\}) \cup \{0 \}$.  Let $x_d^{\times}(S)$ denote the group element given by the monoidal product $x^{\times}_{d^\ast\{0,s_1\}} \otimes \cdots \otimes x^{\times}_{d^\ast\{ 0,s_r\}}$, where $s_1 < s_2 < \cdots < s_r$ are the nonzero elements of $S$. Here, it is to be understood that $x^{\times}_{d^\ast\{0,s\}}$ is the unit element if $d^{\ast}\{ 0,s\}$ has less than two nonzero elements. 
	Similarly, given $g_I = (g_t, g^{\times}_T) \in \G_I$, let $g_d^{\times}(S)$ denote the group element given by the monoidal product $g^{\times}_{d^\ast\{0,s_1\}} \otimes \cdots \otimes g^{\times}_{d^\ast\{ 0,s_r\}}$, with the understanding that $g^{\times}_{d^\ast\{0,s\}}$ is the unit element if $d^\ast\{ 0, s \}$ has less than two nonzero elements.  
	{\bf End of Notation.}
	
	\medskip
	\noindent
	The behaviour of $d^{\mathscr{A}}$ in terms of the description of $\mathscr{A}_I({\bf m})$ given in lemma \ref{TPC1} and the notation introduced above is as follows. 
	On objects, the functor $d^{\mathscr{A}} : \mathscr{A}_I({\bf m}) \longrightarrow \mathscr{A}_J({\bf n})$ is given by $d^{\mathscr{A}}(x_I) = x_J$, with $J := d(I)$ and with $x_J = (x_s, x^{\times}_S) \in {\bf X}^J \times \G_J^{+} = {\bf X}_J$ defined as: 
	\[ x_s = x_d(s ) \, x^{\times}_{d^\ast \{ 0,s \}}, \quad \mbox{and} \quad x^{\times}_S = x_d^{\times}(S)^{-1} \, \sigma_d(S) \, x^{\times}_{d^\ast S}, \]
	where $\sigma_d(S)$ denotes the permutation that sequentially orders all the elements of the ordered set $\{ d^{-1} \{ s_1 \}, \ldots, d^{-1} \{s_r\} \}$, where each of the sets $d^{-1} \{ s_j \}$ are expressed in order induced from $\bf m$. 
	
	\medskip
	\noindent
	On morphisms $(x_I, g_I) \in {\bf{X}}_I \times \G_I$, with  $x_I = (x_t, x^{\times}_T)$ and $g_I = (g_t, g^{\times}_T)$, we define $d^{\mathscr{A}}(x_I, g_I)$ as $(d^{\mathscr{A}}(x_I), d_x(g_I)) \in {\bf{X}}_J \times \G_J$. Here $d_x(g_I)$ is  the element $(g_s, g^{\times}_S) \in \G^J \times \G_J^{+} = \G_J$ defined as:
	\[g_s = (x^{\times}_{d^\ast \{ 0,s \}})^{-1} \, g_d( s ) \, x^{\times}_{d^\ast \{0,s \}} \, g^{\times}_{d^\ast \{0,s\}} \quad \mbox{and} \quad g_S^{\times} = (x_S^{\times})^{-1} \, g_d^{\times} (S)^{-1} \, x_S^{\times} \, g^\times_{d^\ast S}. \]
\end{lemma}
\begin{proof}
	Let us express a morphism $\{\alpha_S\}$ in $\mathscr{A}_I(\bf m)$ as an element $(x_I, g_I) \in {\bf{X}}_I \times \G_I$ in light of lemma \ref{TPC1}. Let $d : {\bf m} \longrightarrow {\bf n}$ be a morphism in $\mathscr{F}$ and let $S \subseteq {\bf n}$ be a pointed subset. In terms of the notation introduced above, and lemma \ref{TPC1}, we have a commutative diagram: 
	\[
	\xymatrix{a_{d^\ast S} \ar[d]^{\alpha_{d^\ast S}} \ar[r]^{F \quad \quad \quad \quad} &  x_d(s_1) \otimes \cdots \otimes x_d(s_r) \ar[d]^{g_d(s_1) \otimes \cdots \otimes g_d(s_r)} \ar[r]^{G} & a_{d^\ast\{0,s_1\}} \otimes \cdots \otimes a_{d^\ast\{0,s_r\}} \ar[d]^{\alpha_{d^\ast\{0,s_1\}} \otimes \cdots \otimes \alpha_{d^\ast\{0,s_r\}} } \\
		b_{d^\ast S} \ar[r]^{F'\quad \quad \quad \quad} & y_d(s_1) \otimes \cdots \otimes y_d(s_r) \ar[r]^{G'} & b_{d^\ast\{0,s_1\}} \otimes \cdots \otimes b_{d^\ast\{0,s_r\}}, }
	\]
	where the morphisms $F,G,F',G'$ may be expressed as:
	\[ F = (x^{\times}_{d^\ast S})^{-1} \, \sigma_d(S)^{-1}, \quad G = x_d^{\times}(S), \quad F' = (g^{\times}_{d^\ast S})^{-1} \, (x^{\times}_{d^\ast S})^{-1} \, \sigma_d(S)^{-1}, \quad G' = x_d^{\times}(S) \, g_d^{\times} (S). \]
	Unraveling the construction \ref{ConstM} of $d^{\mathscr{A}}$ on objects, shows that $d^{\mathscr{A}}(x_I)$ can be expressed in terms of the elements $a_{d^\ast\{0, s\}}$ and the composite morphism $FG$. Similarly, unraveling the construction of $d^{\mathscr{A}}$ on morphisms, shows that $d^{\mathscr{A}}(x_I, g_I)$ can be expressed in terms of the elements $\alpha_{d^\ast \{0,s\}}$, and the composite $FG(F'G')^{-1}$. We leave as an easy exercise for the reader to identify these composites with the formulas in the statement of lemma \ref{retracts}.
\end{proof}

\bigskip
\noindent
Let us explore some of the additional structure implicit in the categories $\mathscr{A}_I(\bf n)$. We will notice in theorem \ref{GRS} that the topological categories $\mathscr{A}_I(\bf n)$ admit retracts $\mathscr{A}_I^S(\bf n)$ where $S$ ranges over pointed subsets of ${\bf n}$. Furthermore, these subcategories will vary functorially in a suitable sense both under inclusions of subsets $S \subseteq R$, and along maps $d : {\bf m} \longrightarrow {\bf n}$ in $\mathscr{F}$. This will lead us an important notion of {\em rigid sections} of the $\Gamma$-space $\BA$ described in theorem \ref{perm}, that plays a crucial role in what follows. 

\medskip
\noindent
We begin with the following two important definitions: \ref{RigSet} and \ref{RigS}:

\bigskip
\begin{defn} \label{RigSet} (The category $\mathscr{F}(\bf n)$ of rigid subsets of $\bf n$ and the functor $d_I$)

	\medskip
	\noindent
	For $n \in \N$, define a Reedy category $\mathscr{F}(\bf n)$ called the category of rigid subsets of ${\bf n}$\footnote{see \cite{RV}, for background on Reedy categories}. Objects of $\mathscr{F}(\bf n)$ are given by pointed subsets $S \subseteq \bf n$, and morphisms are generated by inclusions and retractions, where given $S \subseteq R$, a retraction $R \twoheadrightarrow S$ is defined as the map sending the complement of $S \subseteq R$ to $0$. Denote the raising and lowering subcategories of $\mathscr{F}({\bf m})$, i.e the categories of inclusions and retractions by $\mathscr{F}({\bf m})_{+}$ and $\mathscr{F}({\bf m})_{-}$ respectively. One readily checks that any morphism in $\mathscr{F}(\bf n)$ can be uniquely factored as a composition of a morphism in $\mathscr{F}({\bf m})_{-}$ followed by a morphism in $\mathscr{F}({\bf m})_{+}$.  
	
	
	
	\medskip
	\noindent
	{\bf The functor $d_I$}:
	
	\noindent
	For a fixed $I \in \N^{\times m}$, let $I_0 \subseteq {\bf m}$ be the subset $I_0 = \{ s \in {\bf m} \, |  \, s \neq 0, i_s = 0 \}$. 
	Let $d : {\bf m} \longrightarrow {\bf n}$ be a morphism in $\mathscr{F}$, and let $U \in \mathscr{F}({\bf m})$ be any object. Then we define the object $d_I(U)$ of  $\mathscr{F}({\bf n})$ as 
	\[ d_I(U) := \{ 0 \} \cup \{ t \in d(U) \, | \, t \neq 0, d^{-1} \{t\} \subseteq U \cup I_0\}.\]
	Hence, nonzero elements in $d_I(U)$ are precisely those elements in $d(U)$, all of whose preimages under $d$ belong to $U \cup I_0$. Since $d_I(U) \subseteq d_I(V)$ for $U \subseteq V$, we see that $d_I$ induces a functor 
	\[ d_I : \mathscr{F}({\bf m})_{\pm} \longrightarrow \mathscr{F}({\bf n})_{\pm}, \quad \quad U \mapsto d_I(U).\]
\end{defn}

\medskip
\noindent
The following useful lemma is straightforward, and is left to the reader to verify: 

\medskip
\begin{lemma}\label{RigSet2}
	The functor $d_I :  \mathscr{F}({\bf m})_{\pm} \longrightarrow \mathscr{F}({\bf n})_{\pm}$ has the property that $d_I(U) = d(U)$ if $d$ is injective, and $d_I(d^{-1}(S)) = S$ if $d$ is surjective. Furthermore, $d_I$ composes strictly under morphisms of $\mathscr{F}$. In other words, given a pair of composable morphisms $c,d$ in $\mathscr{F}$, one has a strict equality $(c \circ d)_I = c_J \circ d_I$ of functors on $\mathscr{F}({\bf n})_{\pm}$, where $J = d(I)$ (see example \ref{exampleN}). 
\end{lemma}

\medskip
\begin{remark} \label{d:S}
	Given $d : {\bf m} \longrightarrow {\bf n}$ and a pointed subset $U \subseteq {\bf m}$, Friedlander defines a set $(d:U) \subseteq {\bf n}$ (see \cite{F}, Definition 2.1) that appears to resemble our definition of $d_I(U)$. However, these sets differ in several ways. For instance, unlike $d_I(U)$, the set $(d:U)$ is independent of the multiindex $I$. We will point out in the appendix why we believe that the definition $(d:U)$ in \cite{F} is incorrect. Unfortunately, this mistake pervades most constructions in \cite{F} and calls into question the Model category structure defined in \cite{F} (Definition 2.3). Also see remark \ref{genReedy}. 
\end{remark}

\bigskip
\begin{defn} \label{RigS} (The category of $\Gamma$-spaces with rigid sections)
	
	\noindent
	We define the category to $\Gamma$-spaces endowed with rigid sections as follows: 
	
	\smallskip
	\noindent
	{\bf Objects:}
	
	\noindent
	{\bf (1)} A $\Gamma$-space $\mathscr{E}$ in $\mathscr{F}[\sSet_*]/\mathscr{N}$ is a $\Gamma$-space with rigid sections if it is endowed with subspaces $\mathscr{E}_I^S ({\bf n}) \subseteq \mathscr{E}_I({\bf n})$ (called rigid sections) for each $I \in \N^{\times n}$ and every pointed subset $S \subseteq {\bf n}$. Furthermore, we require that $\mathscr{E}_I^{\bf n} ({\bf n}) = \mathscr{E}_I({\bf n})$. These sections are rigid in the following sense: Given $R \subseteq S$, one has inclusions and retracts that describe a functor on the category $\mathscr{F}(\bf n)$ of rigid subsets of ${\bf n}$ (see definition \ref{RigSet} above):
	\[ \mathscr{E}_I^R({\bf n}) \hookrightarrow \mathscr{E}_I^S({\bf n}), \quad \quad \mathscr{E}_I^S({\bf n}) \twoheadrightarrow \mathscr{E}_I^R(\bf n). \] 
	
	\medskip
	\noindent
	{\bf (2)} Let $d : {\bf m} \longrightarrow {\bf n}$ be a morphism in $\mathscr{F}$ and let $I \in \N^{\times m}$. Then $d^{\mathscr{E}} : \mathscr{E}_I({\bf m}) \longrightarrow \mathscr{E}_{d(I)}(\bf n)$ is required to lift to a map\footnote{Since the lifts in this definition are unique if they exist, we have reused the notation $d^{\mathscr{E}}$ for the lifts to avoid excessive notation.} $d^{\mathscr{E}} : \mathscr{E}_I^{U}({\bf m}) \longrightarrow \mathscr{E}_{d(I)}^{d_I(U)}(\bf n)$
	so that the following diagram commutes:
	\[
	\xymatrix{
		\mathscr{E}^{U}_I({\bf m}) \ar[d]^{d^{\mathscr{E}}}     \ar@{^{(}->}[r] &	\mathscr{E}_I({\bf m}) \ar[d]^{d^{\mathscr{E}}} 
		\\
		\mathscr{E}_{d(I)}^{d_I(U)}(\bf n)  \ar@{^{(}->}[r]  & \mathscr{E}_{d(I)}(\bf n).} 
	\]
	Consider functors defined on the raising and lowering subcategories $\mathscr{F}({\bf m})_{\pm}$ by $U \mapsto \mathscr{E}_I^U({\bf m})$ and $U \mapsto \mathscr{E}_{d(I)}^{d_I(U)}({\bf n})$. 
	Then we also require that the lift $d^{\mathscr{E}} : \mathscr{E}_I^{U}({\bf m}) \longrightarrow \mathscr{E}_{d(I)}^{d_I(U)}(\bf n)$ be a natural transformation of functors on both the categories $\mathscr{F}(\bf m)_{+}$ and $\mathscr{F}({\bf m})_{-}$. 
	
	\smallskip
	\noindent
	Given a $\Gamma$-space $\mathscr{E}$ with rigid sections, we can define a $\Gamma$-space $\mathscr{B}$ as $\mathscr{B}_I := \mathscr{E}_I^{\{ 0 \}}$. The terminal retraction in $\mathscr{F}(\bf n)$ then gives rise to a map of $\Gamma$-spaces $\pi : \mathscr{E} \longrightarrow \mathscr{B}$ in the category  $\mathscr{F}[\sSet_*]/\mathscr{N}$ that admits a section (the $\{ 0 \}$-section). 	
	
	\smallskip
	\noindent
	{\bf Morphisms:} \\
	\noindent
	A morphism $f : \mathscr{E} \longrightarrow \tilde{\mathscr{E}}$ between two $\Gamma$-spaces with rigid sections is a map in $\mathscr{F}[\sSet_*]/\mathscr{N}$ that commutes with the structure described above. It is said to be an equivalence if $f^S_I : \mathscr{E}^S_I \longrightarrow \tilde{\mathscr{E}}^S_I$ is a weak equivalence of simplicial sets for each triple $(n, S \subseteq {\bf{n}}, I \in \N^{\times n})$. Two $\Gamma$-spaces with rigid sections are said to be weakly equivalent if they are connected by a zigzag of equivalences. 
\end{defn}



\smallskip
\noindent
The importance of definition \ref{RigS} is given by the following theorem:

\medskip
\begin{thm}\label{GRS}
	Let $\mathscr{A}$ be a graded topological permutative category. Given $I \in \N^{\times n}$ and a pointed subset $S \subseteq { \bf n}$ with non-zero elements $s_1 < s_2 <  \ldots < s_r$, let $\mathscr{A}_I^S$ be the full subcategory of $\mathscr{A}_I$ with object space given by 
	\[ {\bf{X}}_I^S = {\bf{X}}^I_S \times \G_I^{+}, \quad \mbox{where} \quad {\bf{X}}^I_S = \prod_{k=1}^r {\bf{X}}^{i_{s_k}}, \quad \mbox{and as before} \quad \G_I^{+} = \prod_{T \subseteq \bf n} \G^{\times}_{\mu_I(T)}. \]
	${\bf{X}}_I^S$ is seen as a subspace of ${\bf{X}}_I$ via the inclusion along the basepoints in ${\bf{X}}^{i_s}$ for $s \notin S$. Then the $\Gamma$-space $\BA$ constructed in theorem \ref{perm} is an object in $\mathscr{F}[\sSet_*]/\mathscr{N}$ that admits rigid sections $\BA_I^S(\bf n)$. Notice that in terms of definition \ref{CS}, we may identify $\BA_I^S({\bf n})$ with $\B({\bf{X}}_I^S, \G_I)$. 
\end{thm}
\begin{proof}
	Notice that the object-space ${\bf{X}}_I$ of $\mathscr{A}_I(\bf n)$ defined in lemma \ref{TPC1} has a retract given by the pointed subspace ${\bf{X}}_I^S$. Since the monoids $\G_i$ preserve the basepoints of ${\bf{X}}^i$, we see that $\G_I$ acts on ${\bf{X}}_I^S$. This shows that the topological full  subcategory $\mathscr{A}_I^S(\bf n)$ of $\mathscr{A}_I(\bf n)$ is a retract. It is clear from the definition that given $R \subseteq S$, one has inclusions $\mathscr{A}_I^R({\bf n}) \hookrightarrow \mathscr{A}_I^S(\bf n)$ and retracts $\mathscr{A}_I^S({\bf n}) \twoheadrightarrow \mathscr{A}_I^R(\bf n)$ that define a functor on  $\mathscr{F}(\bf n)$.

	\medskip
	\noindent
	Now recall from lemma \ref{retracts} that if we are given a morphism $d : \bf m \longrightarrow \bf n$ in $\mathscr{F}$, then the functor $d^{\mathscr{A}} : \mathscr{A}_I({\bf m}) \longrightarrow \mathscr{A}_{d(I)}({\bf n})$ was defined on objects by $d^{\mathscr{A}}(x_t, x^{\times}_T) = (x_s, x^{\times}_S)$, where given a nonzero $s$, the element $x_s$ was defined as the right action on the element $x_d(s )$ by the invertible element $x^{\times}_{d^\ast \{ 0,s \}}$, where $x_d(s )$ was the sequential monoidal product of all elements $x_t$ over $t \in d^{-1} \{ s \}$ (or the unit $\ast$ if $d^{-1}\{s \}$ is empty). 
	Now let $j_1 < j_2  < \cdots <  j_n$ be the nonzero elements of $J := d(I)$. Since our monoidal product of non unit objects was assumed to factor through the smash product, we notice that $x_s \in {\bf{X}}^{j_s}$ will be the basepoint if any of the elements $x_t$ for $t \in d^{-1} \{ s \}$ are basepoints of ${\bf{X}}^{i_t}$, unless $i_t = 0$ (since ${\bf{X}}^0$ is the unit).  
	Hence the restriction of the map $d^{\mathscr{A}} : \mathscr{A}({\bf m}) \longrightarrow \mathscr{A}(\bf n)$ to $\mathscr{A}_I^U({\bf m})$ factors uniquely through $\mathscr{A}_{d(I)}^{d_I(U)}({\bf n})$, giving rise to $d^{\mathscr{A}} : \mathscr{A}_I^U({\bf m}) \longrightarrow \mathscr{A}_{d(I)}^{d_I(U)}({\bf n})$.
	Taking the singular nerve, we obtain a  $\Gamma$-space $\BA$ (see theorem \ref{perm}). The structure described above translates to the conditions required in theorem \ref{RigS}, with $\BA$  endowed with the rigid sections $\BA_I^S(\bf n)$. 
\end{proof}

\smallskip
\begin{remark}\label{triv}
	If all the objects ${\bf{X}}^i$ in a graded topological permutative category $\mathscr{A}$ are trivial, i.e. $\mbox{Ob}(\mathscr{A}) = (\N, +)$, then it is easy to see that all inclusions and retractions are isomorphisms for the $\Gamma$-space with rigid sections $\BA$ construction in theorem \ref{GRS} above. 
\end{remark}

\smallskip
\begin{example} \label{N} (The $\Gamma$-space $\mathscr{N}$) 
	
	\noindent
	Consider the graded topological permutative category with the discrete space $\N$ of objects and trivial monoids $\G_i$. It is easy to see that the underlying $\Gamma$-space is none other than the $\Gamma$-space $\mathscr{N}$ of example \ref{exampleN}. Since the above category is terminal among all graded topological permutative categories, we see that all the $\Gamma$-spaces obtained via theorem \ref{GRS} admit canonical maps to $\mathscr{N}$, i.e. are objects in $\mathscr{F}[\sSet_*]/\mathscr{N}$. Remark \ref{E-infinity} identifies the spectrum $\nabla(\mathscr{N})$ of definition \ref{nabla} with the Eilenberg-MacLane spectrum $\HZ$. 
\end{example}

\smallskip
\begin{remark} \label{Sigma}
	In all subsequent examples, the monoids $\G_i$ are either related to the general linear group $\GL_i$ or monoids of self homotopy equivalences of an $i$-fold smash product of a fixed pointed space. As such, each monoid $\G_i$ contains  a subgroup $\Sigma_i$ identified with the permutation group. The morphism $\sigma$ in \ref{TPC} is represented by obvious permutations under this inclusion. 
\end{remark}

\smallskip
\begin{example} \label{GLPCY} (The $\Gamma$-spaces $\BGL^{top}(\C)$ and $\BGL^{top}(\C, \mathbb{S})$)
	
	\noindent
	Let $\GL^{top}_i(\C)$ denote the general linear group seen as a topological group in its usual topology. We define $\GL^{top}_0(\C)$ to be the trivial group. One may define a graded topological permutative category $\GL^{top}(\C)$ with 
	\[ \mbox{Ob}(\GL^{top}(\C)) = (\N, +), \quad \mbox{Mor}(i,j) = \emptyset, \, \,  \mbox{if} \, \,  i \neq j, \, \,  \mbox{and} \quad \mbox{Mor}(i,i) = \GL^{top}_i(\C) \]
	with the objects endowed with the discrete topology. The monoidal structure on morphisms is given by the obvious block diagonal inclusion. Hence, we obtain a special $\Gamma$-space $\BGL^{top}(\C)$ in $\mathscr{F}[\sSet_*]/\mathscr{N}$. Using remark \ref{E-infinity}, we see that the spectrum $\nabla(\BGL^{top}(\C))$ is equivalent to connective complex K-theory $\ku$. 
	
	\noindent
	Now let $\mathbb{S}$ denote the one-point compactification of the topological space $\C$. For $i>0$, let $\mathbb{S}^i$ denote the $i$-fold smash product of $\mathbb{S}$, which is homeomorphic to the one-point compactification of $\C^i$. As such, $\GL^{top}_i(\C)$ acts on $\mathbb{S}^i$. Setting $\mathbb{S}^0$ to be a singleton, define a graded topological permutative category $\GL^{top}(\C, \mathbb{S})$ as: 
	\[ \mbox{Ob}(\GL^{top}(\C, \mathbb{S})) = \coprod_{i \in \N} \mathbb{S}^i, \quad  \mbox{Mor}(\mathbb{S}^i,\mathbb{S}^j) = \emptyset, \, \,  \mbox{if} \, \,  i \neq j, \quad  \mbox{Mor}(\mathbb{S}^i, \mathbb{S}^i) = \mathbb{S}^i \times \GL^{top}_i(\C). \]
	The monoidal structure is defined as the $\GL^{top}_i(\C) \times \GL^{top}_j(\C)$ equivariant smash product
	\[ \mathbb{S}^i \times \mathbb{S}^j \longrightarrow \mathbb{S}^{i+j} \]
	if $i,j > 0$ (or the canonical isomorphism otherwise). Hence we obtain a special $\Gamma$-space $\BGL^{top}(\C, \mathbb{S})$ with rigid sections, endowed with a projection map in $\mathscr{F}[\sSet_*]/\mathscr{N}$:
	\[ \pi_{\G}(\mathbb{S}): \BGL^{top}(\C, \mathbb{S}) \longrightarrow \BGL^{top}(\C). \]
	The above map $\pi_{\G}(\mathbb{S})$ can be identified with the map $\pi$ in definiton \ref{RigS}. 
	By theorem \ref{stablefib}, we see that the spectrum $\nabla(\BGL^{top}(\C, \mathbb{S}))$ is also equivalent to $\ku$. 
\end{example}

\medskip
\noindent
With an eye towards section \ref{FS}, it is helpful to consider more general topological fields $\F$ besides $\C$. For instance, we might consider $\C$ with the discrete topology, as well as the field $\overline{\F}_q$ for some prime $q$, also with the discrete topology. We now explore such examples. 


\smallskip
\begin{example} \label{NormalizerPCX} (The $\Gamma$-spaces $\BN(\F)$ and $\BN(\F, Y)$)) 
	
	\noindent
	Given a topological field $\F$, let $\GL_1(\F)$ denote the topological group of units in $\F$. For $i > 0$, define the topological group $\No_i(\F) := \Sigma_i \ltimes \GL_1(\F)^{\times i}$ seen as a subgroup of $\GL_i(\F)$ generated by the permutation matrices $\Sigma_i$ and the diagonal matrices $\GL_1(\F)^{\times i}$. We define $\No_0(\F)$ to be the trivial group. One may define a graded topological permutative category $\No(\F)$ with 
	\[ \mbox{Ob}(\No(\F)) = (\N, +), \quad \mbox{Mor}(i,j) = \emptyset, \, \,  \mbox{if} \, \,  i \neq j, \, \,  \mbox{and} \quad \mbox{Mor}(i,i) = \No_i(\F). \]
	with the objects endowed with the discrete topology. The monoidal structure on morphisms is given by the obvious block diagonal inclusion. Hence, we obtain a special $\Gamma$-space $\BN(\F)$ in $\mathscr{F}[\sSet_*]/\mathscr{N}$. 
	
	\noindent
	Now assume that $\GL_1(\F)$ admits a right action on a pointed simplicial topological space $Y$. For $i > 0$, $\No_i(\F)$ admits a right action on the $i$-fold smash product of $Y$ denoted by $Y^i$ in a canonical way. Setting $Y^0$ to be the constant point simplicial topological space, define a simplicial graded topological permutative category $\No(\F,Y)$ as: 
	\[ \mbox{Ob}(\No(\F, Y)) = \coprod_{i \in \N} Y^i, \quad  \mbox{Mor}(Y^i,Y^j) = \emptyset, \, \,  \mbox{if} \, \,  i \neq j, \quad  \mbox{Mor}(Y^i,Y^i) = Y^i \times \No_i(\F). \]
	with the monoidal structure on the objects as the smash product $Y^i \times Y^j \longrightarrow Y^{i+j}$ if $i,j > 0$ (or the canonical isomorphism otherwise). The monoidal structure on morphisms induced by the obvious block diagonal inclusion that is compatible with the smash product. Invoking remark \ref{trisimp} we therefore obtain a special $\Gamma$-space $\BN(\F, Y)$ with rigid sections, endowed with a projection map $\pi_{\No}(Y) : \BN(\F, Y) \longrightarrow \BN(\F)$ in $\mathscr{F}[\sSet_*]/\mathscr{N}$ which can be identified with the map $\pi$ in definiton \ref{RigS}. \end{example}

\smallskip
\begin{example} \label{CXX} (The $\Gamma$-spaces $\BGx_X$ and $\BGx_X(X)$)
	
	\noindent
	Let $X$ denote a pointed simplicial set. For $i>0$, let $X^i$ denote the $i$-fold smash product of $X$, and let $X^0$ denote the constant point simplicial set. Let $\G(|X^i|)$ denote the topological monoid of pointed self homotopy equivalences of the realization $|X^i|$ with the compactly generated compact-open topology. The graded topological permutative category $\G_X$ has
	\[ \mbox{Ob}(\G_X) = (\N, +), \quad \mbox{Mor}(i,j) = \emptyset, \, \,  \mbox{if} \, \,  i \neq j, \, \,  \mbox{and} \quad \mbox{Mor}(i,i) = \G(|X^i|). \]
	The monoidal structure on morphisms is induced via the canonical homomorphism:
	\[ \G(|X^i|) \times \G(|X^j|) \longrightarrow \G(|X^{i+j}|) \]
	induced by the smash product map: $|X^i| \times |X^j| \longrightarrow |X^{i+j}|$, 
	when $i,j > 0$ (or the canonical isomorphism otherwise). Let $\BGx_X$ denote the corresponding special $\Gamma$-space over $\mathscr{N}$. We are only aware of a nice description of $\nabla(\BGx_X)$ when $X$ is equivalent to a sphere (see proof of theorem \ref{PLAC}). 
	
	\noindent
	Now consider the canonical (right) action of $\G(|X^{i}|)$ on $|X^i|$. This allows us to define a graded topological permutative category $\G_X(X)$ with 
	\[ \mbox{Ob}(\G_X(X)) = \coprod_{i \in \N} |X^i|, \quad  \mbox{Mor}(|X^i|,|X^j|) = \emptyset, \, \,  \mbox{if} \, \,  i \neq j, \quad  \mbox{Mor}(|X^i|,|X^i|) = |X^i| \times \G(|X^i|) . \]
	with the monoidal structure induced by the smash product:
	\[ |X^i| \times |X^j| \longrightarrow |X^{i+j}|,  \]
	if $i,j > 0$ (or the canonical isomorphism otherwise).
	As before, we obtain a special $\Gamma$-space $\BGx_X(X)$ with rigid sections that comes endowed with a projection map in $\mathscr{F}[\sSet_*]/\mathscr{N}$ $\pi(X) : \BGx_X(X) \longrightarrow \BGx_X$ which can be identified with the map $\pi$ in definition \ref{RigS}.
\end{example}

\medskip
\noindent
We now introduce an important type of special $\Gamma$-space with rigid sections; namely an $X$-fibration for some fixed $X$. This is a key notion relevant to the proof of the $p$-local stable Adams conjecture where $X$ is taken to be a simplicial $p$-local 2-sphere. As was indicated in the introduction, the $J$-homomorphism will arise as a map that classifies such a fibration. 

\medskip
\begin{defn}\label{fib} (The category of sectioned $X$-fibrations)

	\noindent
	Let $X$ be a pointed simplicial set with $X^i$ the $i$-fold smash product of $X$ and $X^0$ the constant point simplicial set. Let $\mathscr{E}$ denote a special $\Gamma$-space with rigid sections and let $\pi : \mathscr{E} \longrightarrow \mathscr{B}$ denote the map described in definition \ref{RigS}. The category of sectioned $X$-fibrations has:
	
	\smallskip
	\noindent
	{\bf Objects:}

	\noindent
	Objects are $\Gamma$-spaces $\mathscr{E}$ with rigid sections called $X$-fibrations, so that $\pi : \mathscr{E} \longrightarrow \mathscr{B}$ satisfies:
	
	\noindent
	{\bf (1)} Given $I = (i_1, \ldots, i_n) \in \N^{\times n}$, the induced map: $\pi_I^S : \mathscr{E}^S_I({\bf n}) \longrightarrow \mathscr{B}_I(\bf n)$, is a fibration with fiber equivalent to space $X^{i_{s_1}} \times \ldots \times X^{i_{s_r}}$, where $s_1 < s_2 < \ldots < s_r$ denotes the nonzero elements of $S \subseteq \bf n$. In particular, the fibration $\pi_I^S$ admits a $\{ 0 \}$-section over $\mathscr{B}_I(\bf n)$. 
	
	\smallskip
	\noindent
	{\bf (2)} This condition demands that the fibers of the map $\pi_I$ defined above factor compatibly along the elements of the rigid subsets $S \subseteq \bf n$. 
	More precisely, given any pointed $S \subseteq \bf n$ and for $0 \neq s \in \bf n$, consider maps $d_s : {\bf n} \longrightarrow \bf{1}$ with $d_s^{-1}(1) = \{ s \}$. We require that the collection of maps $d_s$ induce the following pointed fiberwise equivalence over the map $\mathscr{B}_I({\bf n}) \longrightarrow \prod_{0 \neq s \in S} \mathscr{B}_{i_s}({\bf 1})$:
	\[ \prod_{0 \neq s \in S} d^{\mathscr{E}}_s : \mathscr{E}^{S}_I({\bf n}) \longrightarrow \prod_{0 \neq s \in S} \mathscr{E}_{i_s}({\bf 1}). \]
	{\bf (3)} Given the factorization in the previous conditon, this condition demands that morphisms in $\mathscr{F}$ induce the fiberwise smash product on the factors up to homotopy. More precisely, for a morphism $d : \bf m \longrightarrow \bf n$ in $\mathscr{F}$, and $0 \neq s \in \bf n$, let $T_s = d^{-1}\{ s \} \cup \{0 \} \subseteq \bf m$. Note that $d_I(T_s) = \{0,s\}$ for all $I \in \N^{\times m}$. Let $p\mathscr{E}_I^{T_s}({\bf m})$ denote the homotopy pushout of the diagram (denoted $\mbox{hPO}$) below, which is easily seen to be a diagram over the space $\mathscr{B}_I({\bf m})$: 
	\begin{equation}\label{hPO}
		\xymatrix{
			\hocolim_T \mathscr{E}_I^T({\bf m})\ar[d] \ar[r] \ar@{}[dr]|{\text{hPO}}&  \mathscr{E}^{T_s}_I({\bf m}) \ar[d]
			\\
			\mathscr{B}_I({\bf m}) \ar[r] & p\mathscr{E}_I^{T_s}({\bf m}),  }
	\end{equation}
	where the homotopy colimit is taken over the poset of pointed subsets $T \subset T_s$ so that $d_I(T) = \{ 0\}$ (this is equivalent to the existence of $t \in T_s/T$ such that $i_t \neq 0$). By \cite{MV} (Prop. 5.10.2), we see that the fiber of the map $p\mathscr{E}_I^{T_s}({\bf m}) \longrightarrow \mathscr{B}_I({\bf m})$ is equivalent to its homotopy fiber. Before precisely stating condition { \bf(3)}, we briefly digress in order to identify this fiber with $X^{j_s}$, where $J = d(I)$. 
	
	\smallskip
	\noindent
	The map defined in condition {\bf (2)} induces a map of homotopy colimits over the category defined above. We may easily identify the homotopy colimit for the codomain of this map with a certain fiberwise smash product $\widehat{\prod}$. This gives rise to a commutative diagram \ref{hPO2} below for which the vertical maps are fiberwise equivalence over the map $\mathscr{B}_I({\bf n}) \longrightarrow \prod_{0 \neq t \in T_s} \mathscr{B}_{i_t}({\bf 1})$ and the bottom horizontal map is the map that takes the fiberwise smash product of those factors for which $i_t \neq 0$: 
	\begin{equation}\label{hPO2}
		\xymatrix{
			\mathscr{E}^{T_s}_I({\bf m}) \ar[d]^{\prod_{0 \neq t \in T_s} d^{\mathscr{E}}_t} \ar[r] & p\mathscr{E}_I^{T_s}({\bf m}) \ar[d]^{\widehat{\prod}_{0 \neq t \in T_s} d^{\mathscr{E}}_t}  \\
			\prod_{0 \neq t \in T_s} \mathscr{E}_{i_t}({\bf 1}) \ar[r] & \widehat{\prod}_{0 \neq t \in T_s} \mathscr{E}_{i_t}({\bf 1}).}
	\end{equation}
	
	\smallskip
	\noindent
	By factoring $d$ as a surjection followed by an injection, part {\bf (2)} of definition \ref{RigS} shows that the restriction of $d^{\mathscr{E}} : \mathscr{E}_I^{T_s}({\bf m}) \longrightarrow \mathscr{E}_{d(I)}^{\{ 0, s \}}({\bf n})$ to $\mathscr{E}_I^T({\bf m})$ factors through $\mathscr{B}_I({\bf m})$ for all such $T$. Hence $d^{\mathscr{E}}$ admits the following canonical factorization over $d^{\mathscr{B}} : \mathscr{B}_I({\bf m}) \longrightarrow  \mathscr{B}_{d(I)}({\bf n})$:
	\[ \mathscr{E}^{T_s}_I({\bf m}) \longrightarrow p\mathscr{E}_I^{T_s}({\bf m}) \llra{p{d}^{\mathscr{E}}_s } \mathscr{E}_{d(I)}^{\{ 0, s \}}({\bf n}). \]
	We now state condition {\bf (3)} as the requirement that $p{d}^{\mathscr{E}}_s$ be a pointed fiberwise equivalence over $d^{\mathscr{B}} : \mathscr{B}_I({\bf m}) \longrightarrow  \mathscr{B}_{d(I)}({\bf n})$ for all $s \neq 0$. Note that this fiber is equivalent to $X^{j_s}$, where $J = d(I)$.

	\smallskip
	\noindent
	{\bf Morphisms:} \\
	\noindent
	A morphism $f : (\mathscr{E}, \pi, \mathscr{B}) \longrightarrow (\tilde{\mathscr{E}}, \tilde{\pi}, \tilde{\mathscr{B}})$ of sectioned $X$-fibrations is a morphism $f : \mathscr{E} \longrightarrow \tilde{\mathscr{E}}$ of $\Gamma$-spaces with rigid sections so that the following commutative diagram of spaces is a pointed fiberwise equivalence for each triple $(n, S \subseteq {\bf{n}}, I \in \N^{\times n})$:
	\[
	\xymatrix{
		\mathscr{E}^S_I \ar[d]^{\pi^S_I} \ar[r]^{f^S_I} &  \tilde{\mathscr{E}}^S_I \ar[d]^{\tilde{\pi}^S_I} 
		\\
		\mathscr{B}_I \ar[r]^{f^{\{0\}}_I} & \tilde{\mathscr{B}}_I. }
	\]
	We say that a morphism $f$ is an equivalence of sectioned $X$-fibrations, if $f$ is an equivalence of $\Gamma$-spaces with rigid sections. Two $X$-fibrations are said to be weakly equivalent if they are connected by a zigzag of equivalences. 
\end{defn}

\medskip
\begin{remark}\label{mu2}
	Even though part {\bf (3)} above demands that $d^{\mathscr{E}}$ factor through the homotopy pushout $p\mathscr{E}_I^{T_s}({\bf m})$, part {\bf (2)} of definition \ref{RigS} shows that one has a further factorization of $p{d}^{\mathscr{E}}_s$ of the form
	\[ p{d}^{\mathscr{E}}_s : p\mathscr{E}_I^{T_s}({\bf m}) \longrightarrow  P\mathscr{E}_I^{\, T_s}({\bf m}) \longrightarrow \mathscr{E}_{d(I)}^{\{ 0, s \}}({\bf n})
	\]
	where $P\mathscr{E}_I^{\, T_s}({\bf m})$ is the space over $\mathscr{B}_I({\bf m})$ defined by replacing all homotopy colimits and pushouts in diagram \ref{hPO} by honest colimits and pushouts. 
\end{remark}

\begin{remark} \label{factorization1}
	Friedlander classifies sectioned $X$-fibrations in \cite{F} (Theorem 6.1). However, the reader should notice that our definition of sectioned $X$-fibrations is stronger than the one used by Friedlander in \cite{F,F3} in that we use the full structure of the Reedy category $\mathscr{F}(\bf n)$ on our rigid sections, while Friedlander only works with the subcategory $\mathscr{F}(\bf n)_+$. We have chosen to preserve this extra structure in our definitions, since it appears naturally in all examples of interest, and because it simplifies the proof of Friedlander's theorem that is re-proven in the appendix. 
\end{remark}

\begin{remark}\label{universalqfib}
	In example \ref{CXX}, the map $\pi(X)$ can be seen as a universal quasi-fibration of $\Gamma$-spaces that classifies sectioned $X$-fibrations in much the same way that the universal quasi-fibration classifies fibrations \cite{M1}. We will make this precise in the appendix.  
\end{remark}


\section{From graded algebraic permutative categories to $\Gamma$-spaces with rigid sections.}\label{APC}

\medskip
\noindent
The construction of May \cite{M} described in the previous section can be applied to permutative categories internal to more general categories than just topological spaces. In this section we will describe examples of permutative categories that are analogs of the examples seen in the previous section, but which are internal to simplicial complex schemes. These examples will give rise to $\Gamma$-spaces with rigid sections when we apply certain functors (with values in simplicial sets) to our algebraic constructions. This will allow us to leverage the symmetries inherent in the algebraic context (see Section \ref{FS}). As in definition \ref{TPC}, our definition below refers to these permutative categories internal to complex schemes as graded algebraic permutative categories. 

\medskip
\begin{defn} \label{algperm} (Graded Algebraic Permutative Category)
	
	\noindent
	A graded algebraic permutative category $(\mathscr{A}, \otimes, \ast, \sigma)$ is the following data: 
	
	\medskip
	\noindent
	{\bf (1)} A collection of objects represented by pointed simplicial complex schemes $\mathcal{X}^i, i \in \N$ which are locally of finite type, and with $\ast := \mathcal{X}^0 = \mbox{Spec}(\C)$ seen as a constant simplicial scheme. 
	
	\smallskip
	\noindent
	In addition, we have a collection of morphisms given by complex group schemes $\mathcal{G}_i$ admitting right actions on the pointed simplicial complex schemes $\mathcal{X}^i$, with $\mathcal{G}_0$  the trivial complex group scheme. 
	
	\medskip
	\noindent
	{\bf (2)} A monoidal structure $\otimes$ represented by a collection of maps of pointed simplicial complex schemes on objects, and compatible group homomorphisms on morphisms:
	\[ \otimes : \mathcal{X}^i \times \mathcal{X}^j \longrightarrow \mathcal{X}^{i+j}, \quad \quad \otimes : \mathcal{G}_i \times \mathcal{G}_j \longrightarrow \mathcal{G}_{i+j}, \]
	for which $\ast$ is a strict unit and so that for $i,j > 0$ the sub-simplicial schemes $\ast_i \times \mathcal{X}^j$ and $\mathcal{X}^i \times \ast_j$ map under $\otimes$ to the basepoint of $\mathcal{X}^{i+j}$, where $\ast_i$, $\ast_j$ denote the respective basepoints of $\mathcal{X}^i$ and $\mathcal{X}^j$. 
	
	\medskip
	\noindent
	{\bf (3)} Furthermore, one has a natural isomorphism $\sigma$ represented by a lift:
	\[
	\xymatrix{
		& \mathcal{X}^{i+j} \times \mathcal{G}_{i+j} \ar[d]^{s,t} \\
		\mathcal{X}^i \times \mathcal{X}^j \ar[r]^{\otimes \times \otimes(\tau) \, \, \,\,  } \ar[ur]^{\sigma} & \mathcal{X}^{i+j} \times \mathcal{X}^{i+j}\,}
	\]
	where $s, t$ denote the source and target maps, and $\otimes(\tau)$ denotes monoidal structure applied to the twist map $\tau$ that switches the factors. Moreover, we require $\sigma$ to satisfy the algebraic analog of the identities given in definition \ref{TPC}. As before, a functor of algebraic permutative categories is defined as a functor that preserves the obvious structure. 
\end{defn}

\medskip
\noindent
May's machine that turns a topological permutative category into a $\Gamma$-topological permutative category can easily be seen to extend to a functorial construction that takes an algebraic permutative category $\mathscr{A}$ and constructs a collection $\mathscr{A}(\bf n)$ of categories internal to simplicial complex schemes. 

\medskip
\noindent
Adapting the constructions of lemmas \ref{TPC1} and \ref{retracts}, we see that the objects of $\mathscr{A}(\bf n)$ are given by a coproduct of pointed simplicial complex schemes $\mathcal{X}_I$, indexed by elements $I = (i_1, \ldots, i_n)$ in $\N^{\times n}$, and with $\mathcal{X}_{\underline{0}} = \mbox{Spec}(\C)$ for $\underline{0} = (0,\ldots,0) \in \N^{\times n}$. We may explicitly describe the pointed simplicial complex schemes $\mathcal{X}_I$ as
\[ \mathcal{X}_I =  \mathcal{X}^I \times \mathcal{G}_I^{+}, \quad \mbox{where} \quad \mathcal{X}^I = \prod_{k=1}^n \mathcal{X}^{i_k} \quad \mbox{and} \quad \mathcal{G}_I^{+} = \prod_{T \subseteq \bf n} \mathcal{G}_{\mu_I(T)}, \quad \mbox{with} \quad \mu_I(T) := \sum_{t \in T} i_t. \]
In the above, $T$ denotes those pointed subsets of $\bf n$ that contain at least two nonzero elements. Similarly, the morphisms are generated by complex group schemes $\mathcal{G}_I$ acting on $\mathcal{X}_I$, with $\mathcal{G}_{\underline{0}}$  the trivial group scheme. These group schemes $\mathcal{G}_I$ can be explicitly described as
\[ \mathcal{G}_I =  \mathcal{G}^I \times \mathcal{G}_I^{+}, \quad \mbox{where} \quad  \mathcal{G}^I = \prod_{k=1}^n \mathcal{G}_{i_k}, \]
admitting factorwise right actions on $\mathcal{X}_I$ in the obvious fashion. In analogy with theorem \ref{GRS}, for any pointed subset $S \subseteq \bf n$, let $s_1 < s_2 <  \ldots < s_r$ be the non-zero elements in $S$. Then notice that $\mathcal{X}_I$ has a retract given by the complex simplicial scheme $\mathcal{X}_I^S$:
\[ \mathcal{X}_I^S =  \mathcal{X}^I_S \times \mathcal{G}_I^{+}, \quad \mbox{where} \quad \mathcal{X}^I_S = \prod_{k=1}^r \mathcal{X}^{i_{s_k}}, \]
where $\mathcal{X}_I^S$ is seen as a subscheme of $\mathcal{X}_I$ via the inclusion along the basepoints in $\mathcal{X}^{i_s}$ for $s \notin S$. Since the complex group schemes $\mathcal{G}_i$ preserve the basepoints of $\mathcal{X}^i$, we see that $\mathcal{G}_I$ acts on $\mathcal{X}_I^S$. As in theorem \ref{GRS}, we have corresponding algebraic subcategories $\mathscr{A}_I^S(\bf n)$. 
Given $R \subseteq S$, one has inclusions $\mathscr{A}_I^R({\bf n}) \hookrightarrow \mathscr{A}_I^S(\bf n)$ and retracts $\mathscr{A}_I^S({\bf n}) \twoheadrightarrow \mathscr{A}_I^R(\bf n)$ give rise to a functor on the category $\mathscr{F}(\bf n)$ (see definition \ref{RigSet}). 
Given any pointed subset $U \subseteq {\bf m}$ and $d : {\bf m} \longrightarrow {\bf n}$ a morphism in $\mathscr{F}$, define the pointed subset $d_I(U) \subseteq {\bf n}$ as in \ref{RigSet}, Then point {\bf (2)} of definition \ref{algperm} shows that the restriction of the map $d^{\mathscr{A}} : \mathscr{A}({\bf m}) \longrightarrow \mathscr{A}(\bf n)$ to $\mathscr{A}_I^U({\bf m})$ factors uniquely through $\mathscr{A}_{d(I)}^{d_I(U)}({\bf n})$, giving rise to $d^{\mathscr{A}} : \mathscr{A}_I^U({\bf m}) \longrightarrow \mathscr{A}_{d(I)}^{d_I(U)}({\bf n})$. This leads us to the following definition.

\medskip
\begin{defn} \label{functor} (The $\Gamma$-space $\mbox{E}\mathscr{BA}$ and its rigid sections)
	
	\noindent
	Given a graded algebraic permutative category $\mathscr{A}$, let $\mathscr{A}_I$ denote the category internal to simplicial complex schemes generated by the action of $\mathcal{G}_I$ on $\mathcal{X}_I$ for $I \in \N^{\times n}$. Let $\Delta \circ \Ne \mathscr{A}_I$ be the pointed simplicial complex scheme obtained by taking the diagonal of the bisimplicial complex scheme given by the nerve $\Ne \mathscr{A}_I$. Given a functor $\mbox{E}$ from pointed simplicial complex schemes with values in $\sSet_*$ with the property that $\mbox{E}(\mbox{Spec}(\C)) = \ast$, we define a $\Gamma$-space over $\mathscr{N}$ denoted by $\mbox{E}\mathscr{BA}$, as:
	\[ {\bf n} \mapsto \coprod_{I \in \N^n} \mbox{E} ( \Delta \circ \Ne \mathscr{A}_I). \]
	Notice that $\mbox{E}\mathscr{BA}(\bf n) $ admit rigid sections $\mbox{E}\mathscr{BA}_I^S(\bf n)$ for all pointed subsets $S \subseteq \bf n$ and $I \in \N^{\times n}$, that define a functor on $\mathscr{F}(\bf n)$ as in definition \ref{RigS}. 
\end{defn}

\medskip
\noindent
There are two important (and closely related) functors $\mbox{E}$ we consider. The first functor we consider is the rigid Cech \'etale type of a simplicial complex scheme, or more precisely, its Bousfield-Kan $p$-completion that we describe below. The other functor we consider is the functor of complex points of the simplicial complex scheme. These functors allow us to construct various algebraic analogs of the topological permutative categories considered in section \ref{ToPC}.

\medskip
\begin{defn} \label{EHT} (The rigid Cech \'etale type of complex schemes: \cite{F2}, Section 4)
	
	\noindent
	Fix an algebraically closed field $\varOmega$ so that for any complex scheme $\mathcal{X}$ we consider, all residue fields of $\mathcal{X}$ can be embedded in $\varOmega$ (though no embedding is chosen). Let $\overline{\mathcal{X}}$ denote the set of geometric points of $\mathcal{X}$ defined as arbitrary maps of the form $\overline{x} : \mbox{Spec}(\varOmega) \rightarrow \mathcal{X}$ (see remark \ref{geompoints} below). 
	
	
	\smallskip
	\noindent
	Now, given a simplicial complex scheme $\mathcal{X} = \mathcal{X}_\bullet$ locally of finite type, let $\mbox{RC}(\mathcal{X})$ is a cofiltered category with objects given by rigid \'etale coverings $U \rightarrow \mathcal{X}$, where $U := \{ U_\bullet \rightarrow \mathcal{X}_\bullet \}$ denotes a simplicial family of \'etale coverings with $U_k := \coprod_{\overline{x} \in \overline{\mathcal{X}}_k} U_k(\overline{x})$, and with $U_k(\overline{x})$ a connected \'etale open of a geometric point $\overline{x} \in \overline{\mathcal{X}}_k$. A morphism in $\mbox{RC}(\mathcal{X})$ is a morphiusm $U_{\bullet} \rightarrow V_{\bullet}$ of simplicial schemes over $\mathcal{X}_\bullet$, that fixes the geometric points. We notice that there is at most one morphism between any two objects. This justifies the notion of rigidity in the definition.  
	
	\smallskip
	\noindent
	The rigid Cech \'etale type of $\mathcal{X}$, denoted by $\mathcal{X}_{ret}$ is defined as an inverse system:
	\[ \mathcal{X}_{ret} : \mbox{RC}(\mathcal{X}) \longrightarrow \sSet, \]
	sending $\{U_\bullet \rightarrow \mathcal{X}_\bullet \}$ to the components of the diagonal Cech nerve: $\{ k \mapsto \pi_0(U_{k}^{\times k}) \}$ with the $k$-fold product  taken over $\mathcal{X}_k$. 
	
	\smallskip
	\noindent
	Given $f : \mathcal{X} \longrightarrow \mathcal{Y}$ a map of simplicial complex schemes of locally finite type, then $f$ induces a functor $f^\ast : \mbox{RC}(\mathcal{Y}) \longrightarrow \mbox{RC}(\mathcal{X})$ endowed with a natural transformation $\mathcal{X}_{ret} \circ f^\ast \Rightarrow \mathcal{Y}_{ret}$, where $f^\ast(U_\bullet)_k(\overline{x})$ is defined as the connected neighborhood of $\overline{x}$ in $ \mathcal{X}_k \times_{\mathcal{Y}_k} U_k(f(\overline{x}))$. 
\end{defn}

\smallskip
\begin{remark}\label{geompoints}
	In definition \ref{Galois} we will consider (Galois) symmetries of $\mathcal{X}$ that are {\em NOT} morphisms over $\mbox{Spec}(\C)$. To ensure functoriality of the rigid Cech \'etale type under such morphisms, we view $\mathcal{X}$ as a scheme over $\mbox{Spec}(\Z)$ and not $\mbox{Spec}(\C)$. In particular, even though $\varOmega$ may be endowed with the structure of a $\C$-algebra, our geometric points have not been defined as maps over $\mbox{Spec}(\C)$. 
\end{remark}

\smallskip
\begin{remark}\label{hypercovers}
	The usual rigid \'etale type (\cite{F2}, Definition 4.4) is defined in terms of \'etale hypercovers of schemes. This is in contrast with the rigid Cech \'etale type that we have chosen to work with above, that uses the simpler notion of \'etale covers. For the simplicial schemes $\mathcal{X}_\bullet$ we encounter in this article, each scheme $\mathcal{X}_n$ is a finite coproduct of quasi-projective complex schemes. For such simplicial schemes, the rigid Cech \'etale type is equivalent to the usual \'etale type (\cite{F2}, Prop. 8.2).
\end{remark}

\smallskip
\begin{defn} (The $p$-completed Cech \'etale types: \cite{F2}, Chapter 8; \cite{F}, Section 9)\label{EHTb}
	
	\noindent
	Given a simplicial complex scheme $\mathcal{X}$ locally of finite type as above, the $p$-completed Cech \'etale homotopy type of $\mathcal{X}$ denoted $\mbox{\'et}_p(\mathcal{X})$, is the following Kan-complex:
	\[ \mbox{\'et}_p(\mathcal{X}) := \holim_{RC(\mathcal{X})} (\mathcal{X}_{ret})\hat{\, }_p, \]
	where $(\mathcal{X}_{ret})\hat{\, }_p$ denotes the (pointwise) Bousfield-Kan p-completion of the simplicial set $\mathcal{X}_{ret}$, and $\holim$ denotes the homotopy inverse limit  over $RC(\mathcal{X})$ (see \cite{BK} for the definition of homotopy inverse limit). 
	The pullback construction with respect to maps $f : \mathcal{X} \longrightarrow \mathcal{Y}$ given in definition \ref{EHT} describes $\mbox{\'et}_p(\mathcal{X})$ as a functor on the category of simplicial complex schemes locally of finite type. 
	
	\smallskip
	\noindent
	Let us define the related functor $\mbox{\'Et}_p$ which will be relevant for us (see remark \ref{Etale}):
	\[ \mbox{\'Et}_p(\mathcal{X}) := \Map_{\sSet}(\mbox{\'et}_p(\mbox{Spec}(\C)), \mbox{\'et}_p(\mathcal{X})). \]
	Example \ref{SpecC} shows that $\mbox{\'et}_p(\mbox{Spec}(\C))$ is contractible and therefore $\mbox{\'Et}_p(\mathcal{X})$ is equivalent to $\mbox{\'et}_p(\mathcal{X})$.
\end{defn}

\smallskip
\begin{defn} (Galois symmetries on simplicial complex schemes defined over $\Z$) \label{Galois}
	
	\noindent
	Let $\mbox{Gal}(\C)$ be defined as the group of (not necessarily continuous) field automorphisms of $\C$. Let $\mathcal{X} := \mathcal{X}_{\Z} \times_{\mbox{Spec}(\Z)} \mbox{Spec}(\C)$ be a simplicial complex scheme that is obtained from a simplicial locally noetherian scheme by base change from $\Z$ to $\C$. Then the $\mbox{Gal}(\C)$ action on $\mbox{Spec}(\C)$ extends to an action on $\mathcal{X} =  \mathcal{X}_{\Z} \times_{\mbox{Spec}(\Z)} \mbox{Spec}(\C)$ by acting on the second factor. This  $\mbox{Gal}(\C)$ induces a natural automorphism of $\mathcal{X}_{ret}$ by pulling back pointed \'etale covers, as in definition \ref{EHT}. 
	
	\medskip
	\noindent
	Notice that a geometric point of $\mathcal{X}$ is the same data as a pair:  an $\varOmega$-valued point of $\mathcal{X}_\Z$, and a field embedding $\C \hookrightarrow \varOmega$.  The induced action of $\mbox{Gal}(\C)$ on the set of geometric points has no effect on the $\varOmega$-valued points of $\mathcal{X}_\Z$, and is given by precomposition on the set of field embeddings $\C \hookrightarrow \varOmega$. In particular, the induced  $\mbox{Gal}(\C)$-action on the set of geometric points of $\mathcal{X}$ is free. 
	
	\medskip
	\noindent
	On taking $p$-completions and the inverse limit, we obtain a $\mbox{Gal}(\C)$ action on $\mbox{\'et}_p(\mathcal{X})$. It is straightforward to verify that the $\mbox{Gal}(\C)$ action on $\mbox{\'et}_p(\mathcal{X})$ commutes with morphisms $f : \mathcal{X} \longrightarrow \mathcal{Y}$ that are obtained by extending coefficients from $\Z$ to $\C$. The $\mbox{Gal}(\C)$ action on $\mbox{\'et}_p(\mathcal{X})$ induces an action on $\mbox{\'Et}_p(\mathcal{X})$ by conjugation on the space of maps. 
\end{defn}

\smallskip
\begin{example} (The rigid Cech \'etale type of $\mbox{Spec}(\C)$) \label{SpecC}
	
	\noindent
	Consider $\mbox{Spec}(\C)$ as a constant simplicial complex scheme. Since connected \'etale covers of $\mbox{Spec}(\C)$ are automorphisms of $\mbox{Spec}(\C)$, each geometric point of $\mbox{Spec}(\C)$ has a unique connected \'etale neighborhood up to (unique) isomorphism. Thus $\mbox{RC}(\mbox{Spec}(\C))$ is equivalent to the discrete category with a single object given by $U := \coprod \mbox{Spec}(\C) \longrightarrow \mbox{Spec}(\C)$, where the coproduct is indexed by the set $\overline{\mbox{Spec}}(\C)$ of geometric points of $\mbox{Spec}(\C)$, and each component of $U$ restricts to the identity map of $\mbox{Spec}(\C)$. Since the value of $\mbox{Spec}(\C)_{ret}$ on $U$ is a free $\mbox{Gal}(\C)$-space with $n$-simplices: $\overline{\mbox{Spec}}(\C)^{\times (n+1)}$, we call this space $\mbox{EGal}(\C)$ and may easily verify that it is contractible. It follows that $\mbox{\'et}_p(\mbox{Spec}(\C))$ is also contractible being equivalent to $\mbox{EGal}(\C)\hat{\,}_p$.
\end{example}

\begin{remark} \label{Etale}
	The homotopy type $\mbox{\'et}_p(\mathcal{X})$ is the one that is typically used in applications (see for instance \cite{F2}). However, we prefer working with $\mbox{\'Et}_p(\mathcal{X})$ for the following reason: The pointed simplicial complex schemes $\mathcal{X}$ we consider in this article will all be obtained by base change from a pointed scheme defined over $\Z$. Therefore, the Galois group $\mbox{Gal}(\C)$ naturally acts on $\mbox{\'et}_p(\mathcal{X})$ by \ref{Galois}. Unfortunately, this action will not preserve any choice of basepoints. Replacing $\mbox{\'et}_p(\mathcal{X})$ by $\mbox{\'Et}_p(\mathcal{X})$ (endowed with the conjugation $\mbox{Gal}(\C)$-action), resolves this issue by furnishing us with basepoints which are fixed under $\mbox{Gal}(\C)$. Another crucial fact about $\mbox{\'Et}_p(\mathcal{X})$ is that it supports a canonical map from the ($p$-completed) functor of points as we shall see in theorem \ref{ptstoet} below. 
\end{remark}

\begin{remark} \label{contractibility}
	Even though the value of the functor $\mbox{\'Et}_p$ on the terminal scheme $\mbox{Spec}(\C)$ is a contractible simplicial set, it is not equal to a constant point. Hence, as indicated in the note in definition \ref{Gamma}, one must normalize the construction in definition \ref{functor} to obtain a $\Gamma$-space. 
\end{remark}

\medskip
\begin{defn} (The functor of complex points) 
	
	\noindent
	Given a simplicial complex scheme $\mathcal{X} = \mathcal{X}_\bullet$, define the simplicial set of $\C$-points of $\mathcal{X}_\bullet$ (with the discrete topology) denoted by $\mathcal{X}(\C)_\bullet$, where $\mathcal{X}(\C)_k$ is the set of $\C$-points of $\mathcal{X}_k$, i.e. the set of sections of the structure map $\mathcal{X}_k \rightarrow \mbox{Spec}(\C)$. Notice that if $\mathcal{X}$ is obtained by base change of a simplicial scheme over $\Z$, then $\mathcal{X}(\C)$ admits an induced action of the Galois group $\mbox{Gal}(\C)$ that commutes with maps of simplicial complex schemes that are defined over $\Z$. 
\end{defn}

\medskip

\begin{thm} \label{ptstoet}
	Let $\mathcal{X}$ denote a simplicial complex scheme which is locally of finite type, then there is a natural transformation 
	\[ \eta : \mathcal{X}(\C)\hat{\, }_p \longrightarrow \mbox{\'Et}_p(\mathcal{X}), \]
	where $\mathcal{X}(\C)\hat{\, }_p$ denotes the Bousfield-Kan $p$-completion of the simplicial set $\mathcal{X}(\C)$. Moreover, if $\mathcal{X}$ is obtained by base change from a simplicial scheme over $\Z$, then for any element $\rho \in \mbox{Gal}(\C)$:
	\[ \rho \circ \eta = \eta \circ \rho^{-1}, \]
	where the $\rho$-action on $\mbox{\'Et}_p(\mathcal{X})$ is understood as conjugation on $\Map_{\sSet}(\mbox{\'et}_p(\mbox{Spec}(\C)), \mbox{\'et}_p(\mathcal{X}))$. 
	
	\smallskip
	\noindent
	In particular, given a graded algebraic permutative category $\mathscr{A}$, definition \ref{functor} and the above natural transformation furnish us with an $\mbox{Gal}(\C)$-equivariant rigid section preserving natural transformation of $\Gamma$-spaces
	\[ \mathscr{BA}(\C)\hat{\, }_p \longrightarrow \mbox{\'Et}_p(\mathscr{BA}). \]
\end{thm}
\begin{proof}
	Consider the discrete simplicial complex scheme $\mathcal{X}(\C)_\C := \coprod_{\mathcal{X}(\C)_\bullet} \mbox{Spec}(\C)$. Each geometric point corresponds to a map from $\mbox{Spec}(\C)$ to $\mathcal{X}$. Hence, we have a tautologically defined natural map of simplicial complex schemes $ \mathcal{X}(\C)_\C \longrightarrow \mathcal{X}$
	which induces a natural map $\mbox{\'et}_p(\mathcal{X}(\C)_\C) \longrightarrow \mbox{\'et}_p(\mathcal{X})$. Recall from example \ref{SpecC} that the functor $(\mbox{Spec}(\C))_{ret}$ took the constant value $\mbox{EGal}(\C)$. Since $\mathcal{X}(\C)_{\C}$ is a discrete simplicial complex scheme, the same argument given for $\mbox{Spec}(\C)$ extends to show that $\mbox{RC}(\mathcal{X}(\C)_{\C})$ is equivalent to a discrete category with a single object given by $U \times \mathcal{X}(\C) \longrightarrow \mathcal{X}(\C)_{\C}$, with $U$  as in the case of $\mbox{Spec}(\C)$. The value of $(\mathcal{X}(\C)_\C)_{ret}$ on this object is given by the product simplicial set $\mbox{EGal}(\C) \times \mathcal{X}(\C)$. Hence,  $\mbox{\'et}_p(\mathcal{X}(\C)_\C)$ is equivalent to the simplicial set $(\mbox{EGal}(\C) \times \mathcal{X}(\C))\hat{\, }_p$. Invoking the natural map between the product of $p$-completions and the completion of the product (\cite{BK}, Ch.1, 7.2), we obtain a natural composite map of simplicial sets:
	\[ \mbox{\'et}_p(\mbox{Spec}(\C)) \times \mathcal{X}(\C)\hat{\, }_p \simeq \mbox{EGal}(\C)\hat{\, }_p \times \mathcal{X}(\C)\hat{\, }_p
	\longrightarrow (\mbox{EGal}(\C) \times \mathcal{X}(\C))\hat{\, }_p \simeq \mbox{\'et}_p(\mathcal{X}(\C)_\C) \longrightarrow  \mbox{\'et}_p(\mathcal{X}). \]
	The map $\eta$ we seek is given by taking the adjoint of the above composite. 
	
	\smallskip
	\noindent
	To complete the proof, it remains to show the behavior of $\eta$ under $\rho$. Working locally, we reduce the question to the case of a complex algebra $R$ of the form $R = R_\Z \otimes_\Z \C$ with the induced symmetry $\rho$ acting along the base $\C$. Given a maximal ideal $\mathcal{I} \subset R$ representing a $\C$-point $x$ of $\mbox{Spec}(R)$, let $\rho^{-1}(x)$ denote the $\C$-point represented by $\rho^{-1}(\mathcal{I})$. The proof of the theorem now follows from unraveling the following commutative square:
	\[
	\xymatrix{
		R  \ar[d]^{\rho} \ar[r]^{\rho^{-1}(x)\quad } & R/{\rho^{-1} (\mathcal{I}}) \ar[d]^{\rho} & \ar[l]_{\quad \iso} \C \ar[d]^{\rho}\\
		R   \ar[r]^{x} & R/{\mathcal I} & \ar[l]_{\iso} \C}\]
	with the horizontal maps on the right induced by the canonical inclusion $\C \longrightarrow R$. 
\end{proof}

\medskip
\noindent
\begin{remark} \label{points and CT}
	If $\mathcal{X}$ is a connected simplicial scheme of finite type over $\C$, then Friedlander's comparison theorem (\cite{F2}, Corollary 8.5) shows that there is a zigzag of equivalences $\alpha$ and $\beta$:
	\begin{equation} \label{ab} \mbox{\'et}_p(\mathcal{X}) \llla{\alpha} \mbox{S\'et}_p(\mathcal{X}) \llra{\beta} (\mbox{Sing} \mathcal{X}(\C)^{top})\hat{\, }_p, \end{equation}
	where the functors $\mbox{S\'et}_p(\mathcal{X})$ and  $(\mbox{Sing} \mathcal{X}(\C)^{top})\hat{\, }_p$ are defined as follows. 
	
	\smallskip
	\noindent
	$\mbox{S\'et}_p(\mathcal{X})$ is the functor similar to $\mbox{\'et}_p(\mathcal{X})$ but with components of the diagonal Cech nerve in the definition of $\mathcal{X}_{ret}$ replaced instead by the triple diagonal of the singular simplices on the underlying topological space i.e. $\{ k \mapsto \mbox{Sing}_k({U_k(\C)^{top}})^{\times k} \}$. Here $U_k(\C)^{top}$ denotes the complex points of the \'etale cover $U_k$ endowed with the analytic topology. 
	Similarly, the functor $(\mbox{Sing} \mathcal{X}(\C)^{top})\hat{\, }_p$
	denotes the $p$-completion of the diagonal of the bisimplicial set of singular simplices on the simplicial topological space of complex points of $\mathcal{X}$ with the analytic topology. The maps $\alpha$ and $\beta$ are the maps induced by taking components of singular simplices, and mapping $U_k^{\times k}$ to $\mathcal{X}_k$ respectively. 
	
	\medskip
	\noindent
	Define $\tilde{\alpha}$ and $\tilde{\beta}$ as the maps induced on the mapping spaces out of  $\mbox{\'et}_p(\mbox{Spec}(\C))$ into the zigzag of equation \ref{ab} above. Then, from the definition of $\mbox{S\'et}_p(\mathcal{X})$ and the proof of theorem \ref{ptstoet} we see that $\eta$ lifts canonically to the functor $\mbox{S\'Et}_p(\mathcal{X}):=\Map_{\sSet}(\mbox{\'et}_p(\mbox{Spec}(\C)), \mbox{S\'et}_p(\mathcal{X}))$: 
	\[
	\xymatrix{
		&  \mbox{S\'Et}_p(\mathcal{X})\ar[d]^{\tilde{\alpha}} \ar[r]^{\tilde{\beta} \quad \quad \quad \quad \quad \quad\quad} & \Map_{\sSet}(\mbox{\'et}_p(\mbox{Spec}(\C)), (\mbox{Sing}\mathcal{X}(\C)^{top})\hat{\, }_p) \\
		\mathcal{X}(\C)\hat{\, }_p \ar[r]^{\eta} \ar[ur]^{\tilde{\eta}} &  \mbox{\'Et}_p(\mathcal{X}). & }
	\]
	From the definitions of these maps, we see that the map $\tilde{\beta} \circ \tilde{\eta}$ can be described as the composite
	\[ \mathcal{X}(\C)\hat{\, }_p \longrightarrow (\mbox{Sing} \mathcal{X}(\C)^{top})\hat{\, }_p \hookrightarrow \Map_{\sSet}(\mbox{\'et}_p(\mbox{Spec}(\C)), (\mbox{Sing}\mathcal{X}(\C)^{top})\hat{\, }_p), \]
	where the first map above is induced by the canonical inclusion of $\mathcal{X}(\C)$ (with the discrete topology) into $\mbox{Sing} \mathcal{X}(\C)^{top}$, and the second map is induced by the projection of $\mbox{\'et}_p(\mbox{Spec}(\C))$ to a point. 
\end{remark}

\medskip
\noindent
For the remainder of this section, we will give several important examples of algebraic permutative categories.
In order to construct these examples, we will introduce two types of simplicial schemes $\mathcal{D}^i$ and $\mathcal{S}^i$. The schemes $\mathcal{S}^i$ are closely related to the ones used in \cite{F} and \cite{BhK} (see remark \ref{error} below), while the schemes $\mathcal{D}^i$ are the ones recently introduced by Friedlander \cite{F3}. 

\bigskip
\begin{defn} \label{SGlscheme} (The simplicial schemes $\mathcal{D}^i$ and $\mathcal{S}^i$) 
	
	\smallskip
	\noindent
	{\bf The scheme $\mathcal{D}$:}
	
	\noindent
	Define the pointed simplicial complex schemes $\mathcal{D}^i$ by setting $\mathcal{D}^0 = \mbox{Spec}(\C)$, and defining  $\mathcal{D}^i$ for $i > 0$ to have $s$-simplices 
	\[ \mathcal{D}^i_s = \coprod_{a : [s] \rightarrow [1]} \mathcal{D}^i_a, \quad \mbox{where} \quad 
	\mathcal{D}^i_a = \begin{cases}
		\mbox{Spec}(\C) \quad \mbox{if } \quad a^{-1}(1) = \emptyset \\
		\mathbb{A}^i_\C  \quad \quad \mbox{if} \, \quad a^{-1}(0) = \emptyset \\
		\mathbb{A}^i_\C -\{0\} \quad \mbox{otherwise}, 
	\end{cases}
	\]
	where $a : [m] \rightarrow [1]$ denotes a morphism in the simplex category of ordered sets $[k] = \{0, 1, \ldots, k \}$. The face and degeneracy maps being induced by the canonical maps between the complex schemes. The basepoint is chosen to be the scheme $\mathcal{D}^i_a = \mbox{Spec}(\C)$, where $a : [0] \longrightarrow [1]$, with $a(0) = 0$.

	\smallskip
	\noindent
	{\bf The scheme $\mathcal{S}^i$:}
	
	\noindent
	We define the simplicial complex scheme $\mathcal{S}$ to have $s$-simplices: 
	\[ \mathcal{S}_s = \coprod_{a : [s] \rightarrow [1]} \mathcal{S}_a, \quad \mbox{where} \quad 
	\mathcal{S}_a = \begin{cases}
		\mbox{Spec}(\C) \quad \mbox{if } \quad a^{-1}(1) = \emptyset \\
		\mbox{Spec}(\C)  \quad \mbox{if} \, \, \quad a^{-1}(0) = \emptyset \\
		\mathbb{A}^1_\C -\{0\} \, \quad \mbox{otherwise}. 
	\end{cases}
	\]
	The basepoint is chosen to be the scheme $\mathcal{S}_a = \mbox{Spec}(\C)$, where $a : [0] \longrightarrow [1]$, with $a(0) = 0$. Setting $\mathcal{S}^0 = \mbox{Spec}(\C)$, define $\mathcal{S}^i$ as the $i$-fold smash product of $\mathcal{S}$ with itself when $i >0$. 
	
	\smallskip
	\noindent
	{\bf The product structure for $\mathcal{D}$:}
	
	\noindent
	Next, we construct a map of pointed simplicial complex schemes
	\[ \tilde{\tau}_{i,j} : \mathcal{D}^i \times \mathcal{D}^j \longrightarrow \mathcal{D}^{i+j}, \]
	which is defined to be the obvious isomorphism if $i$ or $j$ is $0$. Otherwise, we identify the left hand side as $ (\mathcal{D}^i \times \mathcal{D}^j)_s = \coprod_{(a, b) : [s] \rightarrow [1]} (\mathcal{D}^i_a \times \mathcal{D}^j_b)$, 
	and define $\tilde{\tau}_{i,j}$ on each component as the canonical map induced via the identification $\mathbb{A}^i_\C \times \mathbb{A}^j_\C = \mathbb{A}^{i+j}_\C$:
	\[ \tilde{\tau}_{i,j} : (\mathcal{D}^i_a \times \mathcal{D}^j_b) \longrightarrow \mathcal{D}_{a b}^{i+j}, \quad \mbox{where} \quad a b : [s] \longrightarrow [1], \quad \mbox{is defined as} \quad a b(r) = a(r) b(r). \]
	
	\noindent
	Notice that if $\mathcal{D}^i_a$ or $\mathcal{D}^j_b$ is $\mbox{Spec}(\C)$, then so is $\mathcal{D}^{i+j}_{a b}$. Similarly, if both $\mathcal{D}^i_a$ and $\mathcal{D}^j_b$ are the affine plane, then so is $\mathcal{D}^{i+j}_{a b}$. Hence $\tilde{\tau}_{i,j}$ is a well-defined morphism of pointed simplicial complex schemes that factors through the smash product when $i,j > 0$. 
	
	\medskip
	\noindent
	It is straightforward to verify that the algebraic complex general linear group scheme $\GLs_i$ acts on $\mathcal{D}^i$ so that the map $\tilde{\tau}_{i,j}$ is equivariant with respect to the diagonal inclusion $\GLs_i \times \GLs_j \rightarrow \GLs_{i+j}$. 
	
	\smallskip
	\noindent
	{\bf The product structure for $\mathcal{S}$:}
	
	\noindent
	As above, we have a map of pointed simplicial complex schemes 
	\[ \tau_{i,j} : \mathcal{S}^i \times \mathcal{S}^j \longrightarrow \mathcal{S}^{i+j} \]
	which is defined to be the obvious isomorphism if $i$ or $j$ is 0, and the smash product map otherwise. Let $\mathcal{N}_i$ denote the complex group scheme $\Sigma_i \ltimes \GLs_1^{\times i}$ seen as the canonical subgroup scheme of $\GLs_i$. Then $\mathcal{N}_i$ acts canonically on $\mathcal{S}^i$, and  $\tau_{i,j}$ is equivariant with respect to the diagonal inclusion $\mathcal{N}_i \times \mathcal{N}_j \longrightarrow \mathcal{N}_{i+j}$.
\end{defn}

\smallskip
\begin{remark}\label{error}
	In \cite{F} (page 139) and \cite{BK} (section 4.1), certain schemes $\mbox{S}^{2i}$ were defined. However, the monoidal structure on $\mbox{S}^{2i}$ was not well defined. Our schemes $\mathcal{S}^i$ in definition \ref{SGlscheme} above are equivalent to $(\mbox{S}^2)^{\wedge i}$ so that the monoidal structure $\tau_{i,j}$ above is well-defined by construction. 
\end{remark}

\begin{remark} \label{models}
	Notice that if $i>0$, then as a simplicial complex scheme, $\mathcal{D}^i$ is the mapping cone of the canonical inclusion $\mathbb{A}^i_\C-\{0\} \subset \mathbb{A}^i_\C$, where $\mathbb{A}^i_\C$ denotes the $i$-dimensional complex affine scheme. Similarly, $\mathcal{S}^i$ is defined as the $i$-fold smash product of $\mathcal{S}$, with $\mathcal{S}$  the mapping cone of the collapse map $\mathbb{A}^1_\C-\{0\} \rightarrow \mbox{Spec}(\C)$. Hence both these schemes serve as models for the $2i$-sphere. 
\end{remark}

\begin{remark} \label{base}
	Notice that all objects and maps defined in definition \ref{SGlscheme} are defined over $\Z$. In particular, one has a $\Z$-form of the above definition by simply replacing the base $\C$ by $\Z$ throughout. 
\end{remark}

\begin{remark} \label{twoS}
	Notice that there is a canonical simplicial projection $z : \mathcal{D}^1 \longrightarrow \mathcal{S}$ that collapses the affine line. This map yields a zigzag of $\mathcal{N}_i$-equivariant maps
	\[ \mathcal{S}^i \llla{\quad z^{i}} (\mathcal{D}^1)^{\wedge i} \llra{\tilde{\tau}^{i}} \mathcal{D}^i, \]
	where $(\mathcal{D}^1)^{\wedge 0}$ is understood as $\mbox{Spec}(\C)$, and so that the following commutative diagram is $\mathcal{N}_i \times \mathcal{N}_j$ equivariant
	\[
	\xymatrix{
		\mathcal{S}^i \times \mathcal{S}^j \ar[d]^{\tau_{i,j}} & (\mathcal{D}^1)^{\wedge i} \times  (\mathcal{D}^1)^{\wedge j} \ar[l]_{z^i \times z^j \quad } \ar[d]^{\wedge} \ar[r]^{\quad \quad \tilde{\tau}^i \times \tilde{\tau}^j} & 
		\mathcal{D}^i \times \mathcal{D}^j \ar[d]^{\tilde{\tau}_{i,j}} \\
		\mathcal{S}^{i+j} & (\mathcal{D}^1)^{\wedge i+j} \ar[l]_{z^{i+j}} \ar[r]^{\tilde{\tau}^{i+j}} & \mathcal{D}^{i+j}.}
	\]
\end{remark}

\medskip
\noindent
Having defined the families of simplicial schemes $\mathcal{D}^i, (\mathcal{D}^1)^i$ and $\mathcal{S}^i$, we may define graded algebraic permutative categories with these as objects. 

\medskip
\begin{defn} \label{Spcat} (The $\Gamma$-spaces $\mbox{\'Et}_p \BGLs(\mathcal{D})$, $\mbox{\'Et}_p \mathscr{B}\mathcal{N}(\mathcal{D}^1)$ and $\mbox{\'Et}_p \mathscr{B}\mathcal{N}(\mathcal{S})$)
	
	\noindent
	Recall that $\GLs_i$ denotes the $i \times i$ general linear group seen as a complex group scheme. We define a topological permutative category $\GLs(\mathcal{D})$ with objects and nontrivial morphisms given by 
	\[ \mbox{Ob}(\GLs(\mathcal{D})) = \{ \mathcal{X}^i :=  \mathcal{D}^i \}_{i \in \N} , \quad \quad \mbox{Mor}(\mathcal{D}^i,\mathcal{D}^i)  = \mathcal{D}^i \times \GLs_i, \]
	and with the monoidal structure $\tilde{\tau}_{i,j}$. Similarly, recall the complex group scheme $\mathcal{N}_i := \Sigma_i \ltimes \GLs_1^{\times i}$ seen as the canonical subgroup scheme of $\GLs_i$. Define graded algebraic permutative categories $\mathcal{N}(\mathcal{D}^1)$ and $\mathcal{N}(\mathcal{S})$ with objects and nontrivial morphisms given by 
	\[ \mbox{Ob}(\mathcal{N}(\mathcal{D}^1)) = \{ \mathcal{X}^i := (\mathcal{D}^1)^{\wedge i} \}_{i \in \N} , \quad  \quad  \mbox{Mor}((\mathcal{D}^1)^{\wedge i},(\mathcal{D}^1)^{\wedge i})  = (\mathcal{D}^1)^{\wedge i} \times \mathcal{N}_i, \]
	\[ \mbox{Ob}(\mathcal{N}(\mathcal{S})) = \{ \mathcal{X}^i := \mathcal{S}^i\}_{i \in \N} , \quad  \quad  \quad \mbox{Mor}( \mathcal{S}^i, \mathcal{S}^i)  =  \mathcal{S}^i \times \mathcal{N}_i,\quad \quad \quad \]
	and with the monoidal structure induced by $\tilde{\tau}$ and $\tau$ respectively. 
\end{defn}

\medskip
\noindent
Applying the functor $\mbox{\'Et}_p$ as described in definition \ref{functor} to definition \ref{Spcat} gives rise to $\Gamma$-spaces with rigid sections $\mbox{\'Et}_p \BGLs(\mathcal{D})$, $\mbox{\'Et}_p \mathscr{B}\mathcal{N}(\mathcal{D}^1)$ and $\mbox{\'Et}_p \mathscr{B}\mathcal{N}(\mathcal{S})$ that satisfy:

\medskip
\begin{thm}\label{special}
	The $\Gamma$-spaces with rigid sections $\mbox{\'Et}_p \BGLs(\mathcal{D})$, $\mbox{\'Et}_p \mathscr{B}\mathcal{N}(\mathcal{D}^1)$ and $\mbox{\'Et}_p \mathscr{B}\mathcal{N}(\mathcal{S})$ fit on a zigzag of $\Gamma$-spaces with rigid sections:
	\[
	\mbox{\'Et}_p \mathscr{B}\mathcal{N}(\mathcal{S}) \llla{z}  \mbox{\'Et}_p \mathscr{B}\mathcal{N}(\mathcal{D}^1)  \llra{\tilde{\tau}} \mbox{\'Et}_p \BGLs(\mathcal{D}). \]
	Furthermore, the above zigzag is $\mbox{Gal}(\C)$-equivariant and the map $z$ above is an equivalence of $\Gamma$-spaces with rigid sections. 
\end{thm}
\begin{proof}
	The existence of the zigzag follows from remark \ref{twoS}. Now, by \ref{base}, the maps $\tilde{\tau}$ and $z$ are defined over $\Z$. Therefore, using definition \ref{Galois} we see that the above zigzag is $\mbox{Gal}(\C)$-equivariant. Finally, notice also that since the map $z : \mathcal{D}^1 \longrightarrow \mathcal{S}$ is a homotopy equivalence on complex points with the analytic topology. Hence the comparison theorem (see remark \ref{points and CT}), shows that the induced map $z$ on the $\Gamma$-spaces with rigid sections above is an equivalence.
\end{proof}

\section{The fracture square and the $p$-local stable Adams conjecture.} \label{FS}

\bigskip
\noindent
In this section we extend the zigzag \ref{special} by another zigzag \ref{zig-zagpts}, to obtain a longer zigzag (equation \ref{zigzag}) of $\Gamma$-spaces which we straighten functorially using lemma \ref{straight}. This allows us to construct a fracture square that gives rise to a sectioned $X$  fibration of $\Gamma$-spaces for $X$ a $p$-local $2$-sphere $\widetilde{\pi}_{(p)}(\mathbb{S}) : \widetilde{\mathscr{B}}{\GL}_{(p)}^{top}(\C, \mathbb{S})  \longrightarrow \widetilde{\mathscr{B}}{\GL}_{(p)}^{top}(\C)$. This fibration will support an action of the Adams operation $\psi^q$ for any prime $q \neq p$. We then feed this fibration into our version of Friedlander's classification result \ref{infinityfib} to prove the stable $p$-local Adams conjecture.  

\medskip
\noindent
Let us fix a prime $p$ and let $q \neq p$ be any other prime. Let $\W$ be the ring of Witt vectors on the separable closure $\overline{\F}_q$ of the finite field $\F_q$ of $q$-elements. In particular, the Frobenius $\Fr^q$ lifts to an automorphism of $\W$.

\medskip
\noindent
Next, let us fix an embedding or rings $\iota : \W \subset \C$ once and for all, and extend the Frobenius automorphism of $\W$ to a Galois symmetry; that is, a field automorphism of $\C$ that we denote by $\psi$. 

\medskip
\begin{remark} \label{FL}
	Observe that the units $\overline{\F}_q^{\times}$ belong to $\W^{\times}$ under the Teichm\"uller lift and can be identified as a subgroup of the roots of unity in $\C$ under the embedding $\iota : \W \subset \C$. By construction, the Galois symmetry $\psi$ described above, and the $q$-power self-map $(\,)^{\wedge q}$, that is defined as $w \mapsto w^q$ for $w \in \C^{\times}$, both agree with the action of the Frobenius automorphism $\Fr^q$ on the image of $\overline{\F}_q^{\times}$ in $\C^{\times}$. This observation has important consequences: 
	
	\medskip
	\noindent
	The three maps above (the $q$-power map $(\,)^{\wedge q}$, the Frobenius $\Fr^q$ and the map $\psi$) induce self-maps of topological spaces we consider in this article that are defined functorially on the underlying groups of units on which these maps are defined. For instance, recall from remark \ref{base} that the simplicial scheme $\mathcal{S}$ defined in \ref{SGlscheme} admits a form over $\Z$.
	
	\medskip
	\noindent
	On taking points over the respective fields, we obtain the following commutative diagram of simplicial topological spaces induced by the inclusion $\overline{\F}_q^\times \subset \C^{\times}$ for both maps, with $\C$  given the discrete topology on the right and the analytic topology on the left.  $\overline{\F}_q$ is given the discrete toplogy
	\begin{equation}\label{Teichunits} 
		\xymatrix{
			\mathcal{S}^{top}(\C)^i \ar[d]_{(\,)^{\wedge q}} &  \mathcal{S}(\overline{\F}_q)^i \ar[l]_{\quad  \iota^{top}} \ar[r]^{\iota} \ar[d]^{\Fr^q} & \mathcal{S}(\C)^i \ar[d]^{\psi}  \\
			\mathcal{S}^{top}(\C)^i &  \mathcal{S}(\overline{\F}_q)^i  \ar[l]_{\quad \iota^{top}} \ar[r]^{\iota} &  \mathcal{S}(\C)^i.
		}
	\end{equation}
\end{remark}

\medskip
\begin{defn}\label{Adams-2} (The Adams operations $\psi^q$)
	
	\noindent
	Recall the topological groups $\No_i(\F) := \Sigma_i \ltimes \GL_1(\F)^{\times i}$ defined in \ref{NormalizerPCX}, where $\F$ is a topological field. Let $\No_i^{top}(\C)$ be the topological group when $\C$ is given the analytic topology. The Teichmuller lift  $\overline{\F}_q^{\times} \subset \W^{\times} \subset \C^{\times}$ induces a zigzag of topological groups, where the group on the left given the analytic topology and the other two given the discrete topology: 
	\[\No_i^{top}(\C) \llla{\quad \iota^{top}} \No_i(\overline{\F}_q) \llra{\iota} \No_i(\C). \]
	As in remark \ref{FL}, the $q$-power self-map of $\C^{\times}$ induces an endomorphism of $\No_i^{top}(\C)$, and the map $\psi$ induces an endomorphism of $\No_i(\C)$ compatibly with the endomorphism of $\No_i(\overline{\F}_q)$ induced by $\Fr^q$. 
	\[ \xymatrix{
		\No_i^{top}(\C)  \ar[d]_{(\,)^{\wedge q}} &  \No_i(\overline{\F}_q) \ar[l]_{\quad  \iota^{top}} \ar[r]^{\iota} \ar[d]^{\Fr^q} &  \No_i(\C) \ar[d]^{\psi}  \\
		\No_i^{top}(\C)  & \No_i(\overline{\F}_q) \ar[l]_{\quad \iota^{top}} \ar[r]^{\iota} &   \No_i(\C).
	}\]
	For the sake of simplicity, maps of topological spaces (or topological groups) induced by any of these three self-maps ($(\,)^{\wedge q}$, $\Fr^q$ and $\psi$) will all be denoted by $\psi^q$, with the domain of any particular map  understood from the underlying ring: $\C^{top}$, $\overline{\F}_q$ or $\C$ respectively. We give these maps the uniform name of ``Adams operations". Note that we do not require $\psi^q$ to be invertible.  
\end{defn}

\medskip
\begin{lemma} \label{Adams}
	Defining $\mathcal{S}^i$ as the $i$-fold smash product of $\mathcal{S}$, then the action of the groups $\No_i$ on $\mathcal{S}^i$ is compatible with the Adams operations. In other words, the following zigzag cube of simplicial topological spaces commutes, with all vertical maps  $\psi^q$ and incoming maps  the group action: 
	\[
	\xymatrix{
		& \mathcal{S}^{top}(\C)^i  \times \No_i^{top}(\C)
		\ar'[d]_{\psi^q}[dd]  \ar[dl] & & \ar[dl] 
		\mathcal{S}(\overline{\F}_q)^i \times \No_i(\overline{\F}_q)
		\ar[ll]_{\iota^{top}} \ar[rr]^{\iota} \ar'[d]_{\psi^q}[dd] && \mathcal{S}(\C)^i \times  \No_i(\C)
		\ar[dl] \ar[dd]_{\psi^q} \\
		\mathcal{S}^{top}(\C)^i \ar[dd]_{\psi^q} & &  \mathcal{S}(\overline{\F}_q)^i \ar[ll]_{\quad \quad \quad \iota^{top}} \ar[rr]^{\quad \quad  \iota} \ar[dd] && \mathcal{S}(\C)^i \ar[dd] &  \\
		&  \mathcal{S}^{top}(\C)^i  \times \No_i^{top}(\C) \ar[dl] & & \ar'[l][ll]_{\quad \quad \quad \quad \quad \iota^{top}}  \mathcal{S}(\overline{\F}_q)^i \times \No_i(\overline{\F}_q)\ar[dl] \ar'[r]^{\quad \quad \quad \quad \iota}[rr] && \mathcal{S}(\C)^i \times  \No_i(\C) \ar[dl]  \\ 
		\mathcal{S}^{top}(\C)^i & &  \mathcal{S}(\overline{\F}_q)^i  \ar[ll]_{\iota^{top}} \ar[rr]^{\iota} &&  \mathcal{S}(\C)^i & 
	}
	\]
\end{lemma}
\begin{proof}
	In order to verify the commutativity of the above cube, let us first make some simplifications. Recall the definition of the simplicial scheme $\mathcal{S}$ given in \ref{SGlscheme} has an integral form we denote by $\mathcal{S}_\Z$ (see remark \ref{base}). By definition, the simplices of $\mathcal{S}_\Z$ are a coproduct of two types of schemes: the terminal scheme, and the multiplicative scheme $\GL_1$. This allows us to reduce the question to the case where $\mathcal{S}^i$ is replaced by the scheme $\GL_1^{\times i}$, and with $\No_i$ acting in the canonical fashion on $\GL_1^{\times i}$. The action of the subgroup of $\No_i$ given by the symmetric group, on $\GL_1^{\times i}$ is given by permuting the factors of $\GL_1^{\times i}$. Since the maps $\iota, \iota^{top}$ and the Adams operations $\psi^q$ in the above diagram preserve the product decomposition, it is sufficient to restrict to the case when $i=1$. In this case, the commutativity of the above cube is equivalent to the following two statements. 
	
	\smallskip
	\noindent
	Firstly, the action of $\psi^q$ on $\GL_1(\C), \GL_1(\overline{\F}_q)$ and $\GL_1^{top}(\C)$ commutes with products (these are inward vertical faces). Secondly, the inclusion maps $\iota^{top} : \GL_1(\overline{\F}_q) \longrightarrow \GL_1^{top}(\C)$ and $\iota : \GL_1(\overline{\F}_q) \longrightarrow \GL_1(\C)$ are group homomorphisms (these are the horizontal faces). The first statement follows from the definition of $\psi^q$ (see remark \ref{FL}) and the second one from the fact that the Teichm\"uller lift is multiplicative. 
\end{proof}

\medskip
\begin{thm} \label{zig-zagpts}
	There exists a zigzag of 
	simplicial topological permutative categories giving rise to a zigzag diagram of special $\Gamma$-spaces with rigid sections endowed with Adams operations $\psi^q$
	\[ \BN^{top}(\C, \mathcal{S}^{top}(\C)) \longrightarrow \BN^{top}(\C, \mathcal{S}^{top}(\C))\hat{\, }_p \llla{\, \, \iota^{top}} \BN(\overline{\F}_q, \mathcal{S}(\overline{\F}_q))\hat{\, }_p \llra{\iota}  \BN(\C, \mathcal{S}(\C))\hat{\, }_p. \]
	Furthermore, the map $\iota^{top}$ is an equivalence and all maps are $\psi^q$-equivariant.
\end{thm} 
\begin{proof}
	As displayed in lemma \ref{Adams}, the zigzag $\mathcal{S}^{top}(\C)^i \llla{\, \iota^{top}} \mathcal{S}(\overline{\F}_q)^i  \llra{\iota} \mathcal{S}(\C)^i$ given by the $i$-fold smash product of the zigzag of equation \ref{Teichunits} is compatible with the actions of the respective groups $\No_i^{top}(\C), \No_i(\overline{\F}_q)$ and $\No_i(\C)$ in a $\psi^q$-equivariant fashion. It therefore induces a zigzag of the corresponding simplicial topological permutative categories as constructed in example \ref{NormalizerPCX}. Hence these maps preserve the rigid sections. 
	Since the Adams operations induce natural endo-transformations of these permutative categories, it is clear that all the corresponding $\Gamma$-spaces admit Adams operations. Finally, using the zigzag cube constructed in lemma \ref{Adams}, one concludes that the zigzag of $\Gamma$-spaces are equivariant with respect to the Adams operations. So it only remains to prove that $\iota^{top}$ is an equivalence. 
	
	\medskip
	\noindent
	Since the $\Gamma$-spaces we are considering are special, the question can be reduced down to showing that the following maps of topological spaces given by the respective homotopy orbits induce isomorphisms in cohomology with $\F_p$-coefficients (see definition \ref{CS}):
	\[ \iota^{top} : |\mathcal{S}(\overline{\F}_q)^i| \times_{\No_i(\overline{\F}_q)} \ENo_i(\overline{\F}_q)  \longrightarrow  |\mathcal{S}^{top}(\C)^i| \times_{\No_i^{top}(\C)} \ENo_i^{top}(\C), \]
	where $|\mathcal{S}(\overline{\F}_q)^i|$ and $|\mathcal{S}^{top}(\C)^i|$ denote the geometric realizations of the respective simplicial topological spaces. 
	Notice that both groups $\No_i(\overline{\F}_q)$ and $\No_i^{top}(\C)$ map to $\Sigma_i$ with kernels  $\GL_1(\overline{\F}_q)^{\times i}$ and $\GL_1^{top}(\C) ^{\times i}$ respectively. Hence both the above homotopy orbits fiber over $\BS_i$ reducing the question down to showing that the map on fibers is a cohomology isomorphism with coefficients in $\F_p$. This map on the fiber is given by the inclusion
	\[ \iota^{top} : |\mathcal{S}(\overline{\F}_q)^i| \times_{\GL_1(\overline{\F}_q)^{\times i}}\EGL_1(\overline{\F}_q)^{\times i} \longrightarrow |\mathcal{S}^{top}(\C)^i| \times_{\GL_1^{top}(\C)^{\times i}}\EGL_1^{top}(\C)^{\times i}. \]
	
	\noindent
	From the definition \ref{SGlscheme}, recall that $\mathcal{S}^i$ was defined as the $i$-fold smash product of $\mathcal{S}$. Moreover, as a $\GL_1$-space, one can identify $\mathcal{S}$ with the unreduced suspension of $\GL_1$, endowed with one of the suspension vertices as the basepoint. It follows that $|\mathcal{S}(\overline{\F}_q)|$ has the structure of a $\GL_1(\overline{\F}_q)$-CW complex with two zero-cells that are both fixed under $\GL_1(\overline{\F}_q)$, and one principal one-cell. The same structure holds for $|\mathcal{S}^{top}(\C)|$ as a $\GL^{top}_1(\C)$-CW complex. Taking the product CW structure, we obtain equivariant CW decompositions of the spaces $|\mathcal{S}(\overline{\F}_q)^i|$ and $|\mathcal{S}^{top}(\C)^i|$ under the groups $\GL_1(\overline{\F}_q)^{\times i}$ and $\GL_1^{top}(\C)^{\times i}$ respectively (with a total of $(2^i+1)$-equivariant cells under the respective groups). By construction, the collection of (equivariant) cells of both spaces are in bijection, and the corresponding isotropy group of a cell in $|\mathcal{S}(\overline{\F}_q)^i|$ includes into the isotropy group of the corresponding cell in $|\mathcal{S}^{top}(\C)^i|$ as an inclusion of the form $\GL_1(\overline{\F}_q)^{\times j} \longrightarrow \GL_1^{top}(\C)^{\times j}$ for some $j \leq i$. Taking the Borel construction converts this to the inclusion on classifying spaces which is known to be an isomorphism in cohomology with $\F_p$-coefficients (\cite{FM}, Proposition 2.3). Since we have a finite number of equivariant cells, what we seek to prove follows from an easy induction on the skeleta.
\end{proof}

\begin{remark} \label{warning}
	Notice that the above theorem is not true if we replace $\mathcal{S}$ by $\mathcal{D}$. This is because the scheme $\mathcal{D}$ is defined using the affine line $\mathbb{A}$ and not the punctured line $\mathbb{A}-\{0\}$. Hence the map $\iota$, from $\overline{\F}_q$-points to $\C$-points, is not well defined for $\mathcal{D}$. This necessitates the use of a zigzag \ref{zigzag} below. 
\end{remark}

\medskip
\noindent
Taking stock of what we have achieve so far, notice that theorem \ref{special} and theorem \ref{zig-zagpts} give rise to two zigzags of special $\Gamma$-spaces with rigid sections. The zigzag of \ref{special} is constructed using the \'etale homotopy type functor, while the zigzag of \ref{zig-zagpts} is constructed via the functor of points. Both zigzags are equivariant with respect to Adams operations, with the action of the Adams operation on the zigzag \ref{special} defined as the symmetry $\psi^{-1}$ as in remark \ref{FL} where we recall that $\psi$ was the field automorphism of $\C$ that lifts the Frobenius on $\mathbb{W}$. We may glue the zigzags of theorem \ref{zig-zagpts} and theorem \ref{special} using theorem \ref{ptstoet} which furnishes us with the map $\BN(\C, \mathcal{S}(\C))\hat{\, }_p  \rightarrow \mbox{\'Et}_p \mathscr{B}\mathcal{N}(\mathcal{S})$ in the diagram \ref{zigzag} below. Note that all wrong way maps in this diagram are equivalences.

\begin{equation} \label{zigzag}
	\BN^{top}(\C, \mathcal{S}^{top}(\C)) \rightarrow  \cdots \rightarrow  \BN(\C, \mathcal{S}(\C))\hat{\, }_p  \rightarrow \mbox{\'Et}_p \mathscr{B}\mathcal{N}(\mathcal{S}) \leftarrow \mbox{\'Et}_p \mathscr{B}\mathcal{N}(\mathcal{D}^1)   \rightarrow  \mbox{\'Et}_p \BGLs(\mathcal{D}).
\end{equation}
Our next task is to straighten such a zigzag functorially. To that end we prove a general straightening lemma. 

\smallskip
\begin{defn} \label{straight0} (The category of zigzags)
	
	\noindent
	Define a category with objects given by zigzags of finite length consisting of morphisms of $\Gamma$-spaces with rigid sections, and such that $g_i$ is an equivalence for all $i$:
	\[ \mathscr{X} \llra{f_0} \mathscr{Y}_0 \llla{g_0} \mathscr{X}_1 \llra{f_1} \mathscr{Y}_1 \cdots \mathscr{Y}_{k-1} \llla{g_{k-1}} \mathscr{X}_k \llra{f} \mathscr{Y}. \]
	Morphisms between any two such zigzags are defined as commutative ladders of maps in the category of $\Gamma$-spaces with rigid sections:
	\[ 
	\xymatrix{
		\mathscr{X} \ar[r]^{f_0} \ar[d] & \mathscr{Y}_0 \ar[d]& \ar[l]_{g_0} \ar[d] \ar[r]^{f_1\quad \quad} \mathscr{X}_1 & \mathscr{Y}_1 \cdots \mathscr{Y}_{k-1} \ar@{-->}[d]& \ar[l]_{\quad \quad g_{k-1}} \mathscr{X}_k \ar[r]^{f} \ar[d] & \mathscr{Y} \ar[d]\\
		\mathscr{X}' \ar[r]^{f'_0} & \mathscr{Y}'_0 & \ar[l]_{g'_0}  \ar[r]^{f'_1\quad \quad} \mathscr{X}'_1 & \mathscr{Y}'_1 \cdots \mathscr{Y}'_{k-1} & \ar[l]_{\quad \quad g'_{k-1}} \mathscr{X}'_k \ar[r]^{f'} &\mathscr{Y}'.}
	\]
	A morphism is said to be an equivalence if all the vertical maps above are equivalences of $\Gamma$-spaces with rigid sections. 
\end{defn}

\smallskip
\begin{lemma} \label{straight}
	There exists a straightening functor from the category of zigzags above to the category whose objects are morphisms of $\Gamma$-spaces with rigid sections and whose morphisms are commuting diagrams. Furthermore, given a zigzag: 
	\[ \mathscr{X} \llra{f_0} \mathscr{Y}_0 \llla{g_0} \mathscr{X}_1 \llra{f_1} \mathscr{Y}_1 \cdots \mathscr{Y}_{k-1} \llla{g_{k-1}} \mathscr{X}_k \llra{f} \mathscr{Y}, \]
	the value of the straightening functor is given by a morphism of the form $\widehat{f} : \widehat{\mathscr{X}} \longrightarrow \mathscr{Y}$, where  $\widehat{\mathscr{X}}$ is a $\Gamma$-space with rigid sections weakly equivalent to $\mathscr{X}$. Under this equivalence, $\widehat{f}$ is homotopic to the composite $f \circ g_{k-1}^{-1} \circ f_{k-1} \circ g_{k-2}^{-1} \cdots \circ f_0$ when evaluated at each triple $(n, S \subseteq {\bf n}, I \in \mathbb{N}^{\times n})$. 
\end{lemma}
\begin{proof} 
	In order to prove the above lemma, we will describe a straightening construction on a zigzag of simplicial sets which is natural with respect to morphisms of zigzags. More precisely, assume that the following is a zigzag of simplicial sets
	\[ X \llra{f_0} Y_0 \llla{g_0} X_1 \llra{f_1} Y_1 \cdots Y_{k-1} \llla{g_{k-1}} X_k \llra{f} Y, \]
	such that $g_i$ is a weak equivalence for all indices $i$. Then we will show that there exists a morphism $\widehat{f} : \widehat{X} \longrightarrow Y$ such that $\widehat{X}$ is weakly equivalent to $X$ and so that under this equivalence, $\widehat{f}$ is homotopic to the composite $f \circ g_{k-1}^{-1} \circ f_{k-1} \circ g_{k-2}^{-1} \cdots \circ f_0$. Moreover, this pair $(\widehat{f}, \widehat{X})$ will be functorial with respect to maps of zigzags of the form above. Applying this construction pointwise to the value of the $\Gamma$-spaces with rigid sections at each triple $(n, S \subseteq {\bf n}, I \in \mathbb{N}^{\times n})$, and invoking functoriality of the construction allows us to deduce the lemma for a zigzag of $\Gamma$-spaces with rigid sections.
	
	\medskip
	\noindent
	In order to make the construction on simplicial sets, it is clear that one only needs to consider a zigzag with $3$-morphisms
	\[ X \llra{f_0} Y_0 \llla{g_0} X_1 \llra{f_1} Y_1, \]
	since then one may repeat the process by taking $\widehat{f}$ as $f_0$. 
	To construct $\widehat{X}$, we fix a functorial factorization of maps of simplicial sets into acyclic cofibrations, followed by fibrations, and functorially factor $f_0$ as $f_0 = \tilde{f_0} \circ i$, where $i : X \longrightarrow \tilde{X}$ is an acyclic cofibration, and $\tilde{f_0} : \tilde{X} \longrightarrow Y_0$ is a fibration. Now define $\widehat{X}, j$ via the following pullback (denoted $\mbox{PB}$), and define $\widehat{f}$ as $f_1 \circ j$.
	\[
	\xymatrix{
		\widehat{X} \ar[d] \ar[r]^{j} \ar@{}[dr]|{\text{PB}}& X_1 \ar[d]^{g_0} 
		\\
		\tilde{X} \ar[r]^{\tilde{f_0}} & Y_0. }
	\]
	Since the classical model structure on simplicial sets is known to be right-proper, the map $\widehat{X} \longrightarrow \tilde{X}$ is a weak equivalence. Naturality of the construction is easily verified. 
\end{proof}

\smallskip
\begin{remark} \label{normalization}
	Depending on the model for functorial factorization on simplicial sets used in the above proof, the value of the $\Gamma$-space $\widehat{X}$ thus obtained on the singleton may just be equivalent (but not equal) to a constant point. We may need to normalize $\widehat{X}$ as in Definition \ref{Gamma} to obtain a $\Gamma$-space with rigid sections whose value on the singleton is the constant point simplicial set. 
\end{remark}

\medskip
\begin{defn} \label{straight2} (Equivalence of maps of $\Gamma$-spaces with rigid sections) 
	
	\noindent
	Define an equivlence relation on the collection of maps of $\Gamma$-spaces with rigid sections for which the generators of the equivalence relation are as follows. Consider a commutative diagram of $\Gamma$-spaces with rigid sections as below, in which all horizontal maps are equivalences
	\[ \xymatrix{
		\tilde{\mathscr{X}} \ar[d] \ar[r] &    \mathscr{Z}  \ar[d] & \ar[l] \mathscr{X} \ar[d] \\
		\tilde{\mathscr{Y}} \ar[r] & \mathscr{W}  &   \ar[l] \mathscr{Y}.
	}\]
	Then this diagram defines an equivalence between the maps $\mathscr{X} \longrightarrow \mathscr{Y}$ and $\tilde{\mathscr{X}} \longrightarrow \tilde{\mathscr{Y}}$.
\end{defn}

\smallskip
\begin{claim}\label{straight3}
	Given a zigzag diagram of $\Gamma$-spaces with rigid sections as in definition \ref{straight0}
	\[ \mathscr{X} \llra{f_0} \mathscr{Y}_0 \llla{g_0} \mathscr{X}_1 \llra{f_1} \mathscr{Y}_1 \cdots \mathscr{Y}_{k-1} \llla{g_{k-1}} \mathscr{X}_k \llra{f} \mathscr{Y}, \]
	let $\widehat{f} : \widehat{\mathscr{X}} \longrightarrow \mathscr{Y}$ denote the straightening of the zigzag from lemma \ref{straight}. Then the equivalence class of the map $\widehat{f}$ as defined in definition \ref{straight2} does not change under the following four situations: 
	
	\smallskip
	\noindent
	{\em (1)} We replace the zigzag by an equivalent zigzag as defined in definition \ref{straight0}. 
	
	\noindent
	{\em (2)} We replace any $f_i$ or any $g_j$ by a homotopic map. 
	
	\noindent
	{\em (3)} If any morphism $g_i$ is the identity, we may drop it and and replace $f_i$ with the composite $f_{i+1} \circ f_i$. Similarly, for $i > 0$, if $f_i$ is the identity, we drop it and replace $g_i$ by $g_{i-1} \circ g_i$. 
	
	\noindent
	{\em (4)} If any morphism $g_{i-1}$ equals the following morphism $f_i$, then we drop the pair of morphisms $g_{i-1}$ and $f_i$ and replace $\mathscr{Y}_i$ with $\mathscr{Y}_{i-1}$. Similarly, if $f_i$ equals equals the following morphism $g_i$, then we may drop the pair and replace $\mathscr{X}_i$ with $\mathscr{X}_{i+1}$. Here $\mathscr{X}_0$ is understood as $\mathscr{X}$.
\end{claim}
\begin{proof}
	The proof of {\em (1)} follows from the naturality of the construction of $\widehat{f}$. To prove {\em (2)}, assume that one has a homotopy $F_i : \mathscr{X}_i \wedge (\Delta[1] \coprod \ast) \longrightarrow \mathscr{Y}_i$ between maps $f_i$ and $\tilde{f}_i$, where $\mathscr{X}_0$ is defined as $\mathscr{X}$. 
	Define a new zigzag by replacing the space $\mathscr{X}_i$ by $\mathscr{X}_i \wedge (\Delta[1] \coprod \ast)$ and keeping all other spaces the same. We redefine the map $f_i$ as the homotopy $F_i$ and the map $g_{i-1} : \mathscr{X}_i \wedge (\Delta[1] \coprod \ast) \longrightarrow \mathscr{Y}_{i-1}$ as the projection to $\mathscr{X}_i$ followed by the map $g_{i-1}$. This new zigzag interpolates (via equivalences) the two zigzags with maps $f_i$ and $\tilde{f}_i$ respectively as the endpoints of $\Delta[1]$. Now one uses property {\em (1)} above to establish the required equivalence. The same proof holds when $f_i$ is replaced by $g_j$. 
	
	\smallskip
	\noindent
	To prove {\em (3)}, we only consider the case when $g_i$ is the identity. The case of $f_i$ is similar. it is enough to consider the case when $g_0$ is the identity. Recall that in the proof of lemma \ref{straight}, we reduced the question to zigzag diagrams of simplicial sets and began the process of straightening by functorially replacing $f_0$ by an acyclic cofibration $i : X \longrightarrow \tilde{X}$ followed by a fibration $\tilde{f}_0$. Then we consider the pullback
	\[
	\xymatrix{
		\widehat{X} \ar[d] \ar[r]^{j} \ar@{}[dr]|{\text{PB}}& X_1 \ar[d]^{g_0} 
		\\
		\tilde{X} \ar[r]^{\tilde{f_0}} & Y_0. }
	\]
	By redefining $f_0$ as $f_1 \circ j$, we continue the straightening process. Notice that if $g_0$ is the identity, then $j = \tilde{f}_0$. Hence one has a map $i : X \longrightarrow \widehat{X}$ that allows us to construct an equivalence between the original zigzag where $g_0$ has been dropped and $f_0$ replaced by $f_1 \circ f_0$, and the new zigzag with $f_0 = f_1 \circ j$. Since all constructions are functorial, this extends from simplicial sets to $\Gamma$ spaces with rigid sections. 
	
	\smallskip
	\noindent
	Finally, for the first part of {\em (4)}, assume that $g_{i-1} = f_i$. Consider the commutative diagram
	\[ \xymatrix{
		\mathscr{Y}_{i-1} \ar[d]^{=} &   \ar[l]_{g_{i-1}} \mathscr{X}_i \ar[r]^{f_i} \ar[d]^{g_{i-1}} &  \mathscr{Y}_i \ar[d]^{=} \\
		\mathscr{Y}_{i-1}  & \mathscr{Y}_{i-1}  \ar[l]_{=} \ar[r]^{=} &   \mathscr{Y}_i.
	}\]
	We may extend both rows as in the original zigzag to get an equivalence of zigzags. The proof now follows using parts {\em (1)} and {\em (3)}. The proof of the second par of {\em (4)} is similar. 
\end{proof}

\medskip
\noindent
Applying lemma \ref{straight} to straighten diagram \ref{zigzag}, we obtain a map $\kappa$ of $\Gamma$-spaces with rigid sections endowed with Adams operations
\begin{equation} \label{kappa}
	\kappa : \widehat{\mathscr{B}}{\No}^{top}(\C, \mathcal{S}^{top}(\C))  \longrightarrow  \mbox{\'Et}_p \BGLs(\mathcal{D}).
\end{equation}
Ignoring Adams operations, we will now identify the map $\kappa$, up to equivalence. 

\medskip
\begin{lemma}\label{normalization2}
	The map $\kappa$ given in \ref{kappa} above is equivalent (as in definition \ref{straight2}) to the canonical map of $\Gamma$-spaces with rigid sections
	\[ \mathscr{B}{\No}^{top}(\C, \mathbb{S}) \longrightarrow \BGL^{top}(\C, \mathbb{S})\hat{\, }_p \]
	where the co-domain is the $p$-completion of the $\Gamma$-space in example \ref{GLPCY}, and the domain is the $\Gamma$-space of example \ref{NormalizerPCX} where $Y$ is taken to be the compactification $\mathbb{S}$ of the complex numbers $\C$ given the usual (analytic) topology.
\end{lemma}
\begin{proof}
	In order to identify $\kappa$, we will simplify the zigzag diagram \ref{zigzag} via a sequences of equivalences. Since the zigzag \ref{zigzag} was obtained by glueing the two zigzags from theorem \ref{special} and theorem \ref{zig-zagpts}, we will begin by focusing on the former. Let us recall from remark \ref{points and CT} that one has a diagram of zigzags
	
	\begin{equation}  \label{Spcat3} \xymatrix{
			\mathscr{B}{\No}^{top}(\C, \mathcal{S}^{top}(\C))\hat{\, }_p  &   \ar[l]_{z} \mathscr{B}{\No}^{top}(\C, \mathcal{D}^{top}(\C)^1)\hat{\, }_p  \ar[r]^{\tilde{\tau}} &  \BGL^{top}(\C, \mathcal{D}^{top}(\C))\hat{\, }_p\\
			\mbox{S\'Et}_p \mathscr{B}\mathcal{N}(\mathcal{S}) \ar[u]_{\beta} \ar[d]^{\tilde{\alpha}}&   \ar[l]_{z} \mbox{S\'Et}_p \mathscr{B}\mathcal{N}(\mathcal{D}^1)  \ar[r]^{\tilde{\tau}} \ar[d]^{\tilde{\alpha}} \ar[u]_{\beta} &  \mbox{S\'Et}_p \BGLs(\mathcal{D}) \ar[d]^{\tilde{\alpha}} \ar[u]_{\beta}\\
			\mbox{\'Et}_p \mathscr{B}\mathcal{N}(\mathcal{S}) &   \ar[l]_{z} \mbox{\'Et}_p \mathscr{B}\mathcal{N}(\mathcal{D}^1)  \ar[r]^{\tilde{\tau}} &  \mbox{\'Et}_p \BGLs(\mathcal{D}).} 
	\end{equation} 
	In defining the upward pointing maps $\beta$ we have composed the maps $\tilde{\beta}$ of remark \ref{points and CT}, which take values in the mapping space
	$\Map_{\sSet}(\mbox{\'et}_p(\mbox{Spec}(\C)), (\mbox{Sing}\mathcal{X}(\C)^{top})\hat{\, }_p)$, 
	with the map that evaluates at the canonical geometric point in $\mbox{\'et}_p (\mbox{Spec}(\C))$. 
	
	\medskip
	\noindent
	Moreover, from remark \ref{points and CT}, we also know that the map $\eta : \BN(\C, \mathcal{S}(\C))\hat{\, }_p  \rightarrow \mbox{\'Et}_p \mathscr{B}\mathcal{N}(\mathcal{S})$ (which is used to glue the two zigzags \ref{special} and \ref{zig-zagpts}), lifts to the interpolating singular etale type $\mbox{S\'Et}_p \mathscr{B}\mathcal{N}(\mathcal{S})$. Using claim \ref{straight2} (part {\em (1)}), we see that the equivalence class of $\kappa$ does not change if one uses the top row of diagram \ref{Spcat3} above instead of the bottom row (which was the one described in theorem \ref{special}), when glueing to the zigzag of theorem \ref{zig-zagpts}. 
	
	\medskip
	\noindent
	Once we perform the glueing of the zigzag of theorem \ref{zig-zagpts} and the top row of \ref{Spcat3}, we obtain a longer zigzag in which the map $\iota^{top}$ occurs sequentially. Hence, by claim \ref{straight3} (part {\em (4)}) we may drop $\BN(\overline{\F}_q, \mathcal{S}(\overline{\F}_q))\hat{\, }_p$. This shows that the map $\kappa$ obtained by straightening the zigzag diagram \ref{zigzag}, is equivalent to the map obtained by straightening the following zigzag
	
	\begin{equation} \label{zigzag2}
		\BN^{top}(\C, \mathcal{S}^{top}(\C)) \rightarrow \BN^{top}(\C, \mathcal{S}^{top}(\C))\hat{\, }_p \llla{z} \mathscr{B}{\No}^{top}(\C, \mathcal{D}^{top}(\C)^1)\hat{\, }_p   \llra{\tilde{\tau}} \BGL^{top}(\C, \mathcal{D}^{top}(\C))\hat{\, }_p . 
	\end{equation}
	Our final task is to relate the simplicial topological spaces $\mathcal{D}^{top}(\C)^i$ and $\mathcal{S}^{top}(\C)$ with the compactification of $\C$. For this, recall from definition \ref{SGlscheme} that $\mathcal{D}^{top}(\C)^i$ is the simplicial mapping cone of the canonical inclusion of topological spaces $(\C^i-\{0\}) \subset \C^i$. Similarly, $\mathcal{S}^{top}(\C)$ was defined as the simplicial mapping cone of the collapse map $(\C-\{0\}) \rightarrow \ast$. We now define simplicial subspaces of these spaces. 
	
	\medskip
	\noindent
	Define the $\U(1)$-equivariant simplicial compact form $\mathcal{S}^{top}(\C)_c \subset \mathcal{S}^{top}(\C)$ by replacing the space $\C-\{0\}$ by the unit circle $\mbox{S}^1$. Similarly, we define the $\U(1)$-equivariant simplicial compact form $\mathcal{D}^{top}(\C)^1_c \subset \mathcal{D}^{top}(\C)^1$ by replacing the space $\C-\{0\}$ by the unit circle $S^1$ and the affine line $\C$ by the unit disc $\mbox{D}$. Finally, for arbitrary $i > 0$, define the $\U(i)$-equivariant subspace $\mathcal{D}^{top}(\C)^i_+ \subset \mathcal{D}^{top}(\C)^i$ to be the simplicial subspace with $\C^i$ unchanged, and the topological space $\C^i - \{ 0 \}$ replaced by $\overline{\mbox{D}}^i := \{ v \in \C^i - \{ 0\}, \, \,  \mbox{such that} \, \, |v| \geq 0\}$. 
	
	\medskip
	\noindent
	Furthermore, the smash product on $\mathcal{D}^{top}(\C)^i$ restricts to a $\U(i) \times \U(j)$-equivariant smash product 
	\[ \mathcal{D}^{top}(\C)^i_+  \times \mathcal{D}^{top}(\C)^i_+ \longrightarrow  \mathcal{D}^{top}(\C)^{i+j}_+ \]
	With the above replacements, the geometric realization of the simplicial topological space $\mathcal{D}^{top}(\C)_c^1$ is homeomorphic to a $2i$-sphere where one hemesphere (the ``algebraic hemesphere") is canonically identified with the disc $\mbox{D}$ and the other hemisphere (the ``topological hemisphere") with the cone on $\mbox{S}^1$. Similarly, the geometric realization of the simplicial topological space $\mathcal{S}^{top}(\C)_c$ is homeorphic to the $2$-sphere given by the unreduced suspension on $\mbox{S}^1$. As in the case of $\mathcal{D}^{top}(\C)_c^1$, we may decompose $\mathcal{S}^{top}(\C)_c$ into two hemispheres, one of which contains the basepoint. 
	
	\medskip
	\noindent
	Now the subspace inclusion into the zigzag \ref{twoS} induces an equivalence of  $\Sigma_i \ltimes (\U(1)^{\times i})$-equivariant zigzags which is compatible with the smash product:
	\begin{equation} \label{zigzag3}	 \xymatrix{ |\mathcal{S}^{top}(\C)_c^{\wedge i}| \ar[d] & \ar[l]_{z^i} {|\mathcal{D}^{top}(\C)_c^1}^{\wedge i}| \ar[r]^{\tilde{\tau}^i} \ar[d] & |\mathcal{D}^{top}(\C)_+^i|  \ar[d]\\
			|\mathcal{S}^{top}(\C)^{\wedge i}| & \ar[l]_{z^i} {|\mathcal{D}^{top}(\C)^1}^{\wedge i}| \ar[r]^{\tilde{\tau}^i} & |\mathcal{D}^{top}(\C)^i|. }
	\end{equation}
	Recall that $z^i$ is given by pinching away the hemisphere $\mbox{D}$. However, if we replace $z^i$ by the homeomorphism $z_{alg}^i$ that identifies $|\mathcal{D}^{top}(\C)_c^1|$ with $|\mathcal{S}^{top}(\C)_c|$ by canonically identifying their respective hemispheres, then one may obtain  $\Sigma_i \ltimes (\U(1)^{\times i})$-equivariant homotopy between $z^i$ and $z_{alg}^i$ that is compatible with the smash product. This homotopy simply deforms the the hemisphere $\mbox{D} \subset |\mathcal{D}^{top}(\C)_c^1|$ to a point. 
	
	\smallskip
	\noindent
	One has a $\Sigma_i \ltimes (\U(1)^{\times i})$-equivariant commutative diagram below with $\mathbb{S}^{2i}$ denoting the compactification of the topological space $\C^i$, which we have identified with $|\mathcal{S}^{top}(\C)_c^{\wedge i}|$
	\begin{equation} \label{zigzag4} \xymatrix{
			|\mathcal{S}^{top}(\C)_c^{\wedge i}| \ar[d]^{=}&   \ar[l]_{z_{alg}^i} {|\mathcal{D}^{top}(\C)_c^1}^{\wedge i}| \ar[r]^{\tilde{\tau}^i} \ar[d]^{z_{alg}^i} &  |\mathcal{D}^{top}(\C)_+^i|\ar[d]^{h^i} \\
			\mathbb{S}^{2i}  & \mathbb{S}^{2i} \ar[l]_{=} \ar[r]^{=} &   \mathbb{S}^{2i}.
	} \end{equation} 
	The map $h^i$ is the $\U(i)$-equivariant map taking values in $\mathbb{S}^{2i}$ is defined as follows
	\[ h^i(r,w) = \frac{1}{r} \, w,  \quad \mbox{where} \quad 0 \leq r \leq 1, \, \, \mbox{and with} \, \, w \in \overline{\mbox{D}}^i, \, \mbox{if} \,\, r < 1 \, \, \mbox{and} \, \, w \in \C^i \, \, \mbox{if} \, \, r=1. \]
	In the above formula, $r$ represents the cone coordinate with $r=0$ denoting the basepoint. Applying the singular bar construction with respect to the group $\Sigma_i \ltimes (\U(1)^{\times i})$ for the two topological spaces on the left of the rows of \ref{zigzag3} and with respect to $\U(i)$ for the topological spaces on the right,  and then $p$-completing, the top row yields a zigzag equivalent to zigzag \ref{zigzag2}. 
	Invoking \ref{straight3} (parts {\em(1) and \em (2)}), we see that the map obtained by straightening \ref{zigzag2} is equivalent to the one obtained if we replace $z$ by $z_{alg}$. Then invoking \ref{straight3} (part {\em (1)}), we may work instead with the $p$-completed bar construction on the bottom row of \ref{zigzag4} without changing the equivalence class of the straightening. The $p$-completed bar construction on the bottom row of \ref{zigzag4} is equivalent to:
	\begin{equation}\label{zigzag5}
		\mathscr{B}{\No}^{top}(\C, \mathbb{S}) \lra \BN^{top}(\C, \mathbb{S})\hat{\, }_p     \llla{=} \BN^{top}(\C, \mathbb{S})\hat{\, }_p  \lra \BGL^{top}(\C, \mathbb{S})\hat{\, }_p.  
	\end{equation}
	Hence we conclude that the map $\kappa$ is equivalent to the straightening of \ref{zigzag5}, which is equivalent to the canonical map $\mathscr{B}{\No}^{top}(\C, \mathbb{S}) \longrightarrow \BGL^{top}(\C, \mathbb{S})\hat{\, }_p$ on invoking claim \ref{straight3} (part {\em (3)}). This is what we were seeking to prove. 
\end{proof}

\medskip
\noindent
The above discussion justifies the construction of the following fracture square: 

\medskip
\begin{defn} \label{AFS} (The arithmetic fracture square) 
	
	\noindent
	Let $l : \mbox{X} \rightarrow \mbox{X}_{\Q}$ denote the Bousfield-Kan localization with respect to the rationals \cite{BK}, and let $\kappa_{\Q}$ denote the rationalization of the map $\kappa$ constructed in \ref{kappa} above. Define a $\Gamma$-space $\BGL_{(p)}^{top}(\C, \mathbb{S})$ with rigid sections, endowed with an Adams operation $\psi^q$, so that the following is a homotopy pullback (denoted $\mbox{hPB}$) for each triple $(n, S \subseteq {\bf{n}}, I \in \N^{\times n})$:
	\[
	\xymatrix{
		\BGL_{(p)}^{top}(\C, \mathbb{S})^S_I  \ar[d] \ar[r] \ar@{}[dr]|{\text{hPB}} & \mbox{\'Et}_p \BGLs(\mathcal{D})^S_I \ar[d]^{l} \\
		(\widehat{\mathscr{B}}{\No}^{top}(\C, \mathcal{S}^{top}(\C))_{\Q})^S_I \ar[r]^{\quad \kappa_{\Q}} & (\mbox{\'Et}_p \BGLs(\mathcal{D})_{\Q})^S_I.}
	\]
	As in the proof of lemma \ref{straight}, the homotopy pullback is constructed by first functorially replacing  $\kappa_{\Q}$ evaluated on a triple $(n, S \subseteq {\bf{n}}, I \in \N^{\times n})$ by a fibration, and then defining $\BGL_{(p)}^{top}(\C, \mathbb{S})^S_I$ as the pullback. Note that we may need to normalize the value on the singleton (see remark \ref{normalization}). 
\end{defn}

\medskip
\noindent
The next theorem uses the fracture square above to construct a sectioned fibration endowed with an Adams operation $\psi^q$. 

\medskip
\begin{thm} \label{main}
	The $\Gamma$-space with rigid sections $\BGL_{(p)}^{top}(\C, \mathbb{S})$ constructed in definition \ref{AFS} is equivalent to a $\Gamma$-space $\widetilde{\mathscr{B}}{\GL}_{(p)}^{top}(\C, \mathbb{S})$ with rigid sections, and endowed with an Adams operation $\psi^q$. Moreover, $\widetilde{\mathscr{B}}{\GL}_{(p)}^{top}(\C, \mathbb{S})$ satisfies the following two properties.
	
	\smallskip
	\noindent
	Firstly, $\widetilde{\mathscr{B}}{\GL}_{(p)}^{top}(\C, \mathbb{S})$ is equivalent to the $p$-localization of the $\Gamma$-space with rigid sections constructed in example \ref{GLPCY}. Secondly, the map of $\Gamma$-spaces described in definition \ref{RigS}:
	\[ \widetilde{\pi}_{(p)}(\mathbb{S}) : \widetilde{\mathscr{B}}{\GL}_{(p)}^{top}(\C, \mathbb{S})  \longrightarrow \widetilde{\mathscr{B}}{\GL}_{(p)}^{top}(\C) := \widetilde{\mathscr{B}}{\GL}_{(p)}^{top}(\C, \mathbb{S})^{\{ 0\}}\]
	is a sectioned $\mathbb{S}_{(p)}$-fibration, where $\mathbb{S}_{(p)}$ is any fixed model for the $p$-local 2-sphere. 
	
	\smallskip
	\noindent
	In particular, we conclude that the $p$-localization of the map $\pi_{\G}(\mathbb{S})$ from example \ref{GLPCY} admits a $\psi^q$-equivariant model as a sectioned $\mathbb{S}_{(p)}$-fibration. 
\end{thm}

\begin{proof}
	By lemma \ref{normalization2} we know that the map $\kappa$ in \ref{kappa} is equivalent to the canonical map $\mathscr{B}{\No}^{top}(\C, \mathbb{S}) \longrightarrow \BGL^{top}(\C, \mathbb{S})\hat{\, }_p$. The map $\BN^{top}(\C, \mathbb{S})^{\{0\}}_I \longrightarrow \BGL^{top}(\C,\mathbb{S})^{\{0\}}_I$ induces a rational homology isomorphism for any $I \in \N^{\times n}$ since it is equivalent to the inclusion of the normalizer of the maximal torus in (products of) general linear groups (which is well known to be a rational equivalence). As a consequence, it follows that the map $\BN^{top}(\C, \mathbb{S})^S_I \longrightarrow \BGL^{top}(\C, \mathbb{S})^S_I$ is also a rational homology isomorphism for any $S \subseteq {\bf n}$. These maps therefore turn into equivalences under Bousfield-Kan localization with respect to $\Q$. We conclude that for any triple $(n, S \subseteq {\bf n}, I \in \N^{\times n})$, the diagram \ref{AFS} is equivalent to the arithmetic fracture square (\cite{BK}, Ch 6, Section 8) whose homotopy pullback is equivalent to the $p$-localization of $\BGL^{top}(\C, \mathbb{S})$. Since all the equivalences above are natural with respect to morphisms in $\mathscr{F}(\bf n)$ and $\mathscr{F}$, the equivalence between $\mathscr{B}{\GL}_{(p)}^{top}(\C, \mathbb{S})$ and the $p$-localization of the $\Gamma$-space with rigid sections $\BGL^{top}(\C, \mathbb{S})$, defined in example \ref{GLPCY}, is an equivalence of $\Gamma$-spaces with rigid sections. 
	
	\smallskip
	\noindent
	Next we apply the functor \ref{Reedy} to $\mathscr{B}{\GL}_{(p)}^{top}(\C, \mathbb{S})$ to obtain an equivalent $\Gamma$-space $\widetilde{\mathscr{B}}{\GL}_{(p)}^{top}(\C, \mathbb{S})$ with rigid sections endowed with an Adams operation $\psi^q$. Consider the map $\widetilde{\pi}_{(p)}(\mathbb{S})$ as constructed in definition \ref{RigS}:
	\[ \widetilde{\pi}_{(p)}(\mathbb{S}) : \widetilde{\mathscr{B}}{\GL}_{(p)}^{top}(\C, \mathbb{S}) \longrightarrow \widetilde{\mathscr{B}}{\GL}_{(p)}^{top}(\C) :=  \widetilde{\mathscr{B}}{\GL}_{(p)}^{top}(\C, \mathbb{S})^{\{ 0\}}. \]
	
	\smallskip
	\noindent
	It remains to show that $\widetilde{\pi}_{(p)}(\mathbb{S})$ is a sectioned $\mathbb{S}_{(p)}$-fibration. To establish this, let us begin by noticing that lemma \ref{Reedy} shows that $\widetilde{\pi}_{(p)}(\mathbb{S})$ is a fibration of simplicial sets. It follows from the equivalence of $\widetilde{\mathscr{B}}{\GL}_{(p)}^{top}(\C, \mathbb{S})$ with the $p$-localization of $\BGL^{top}(\C, \mathbb{S})$ established above, that the fibers of $\widetilde{\pi}_{(p)}(\mathbb{S})$ have the right homotopy type. The only property in definition \ref{fib} that requires some explaination, is condition ${\bf (3)}$. This readily follows from the equivalence of the rigid sections of $\widetilde{\pi}_{(p)}(\mathbb{S})$ with those of the $p$-localization of $\pi_{\G}(\mathbb{S})$. 
\end{proof}

\medskip
\begin{remark} \label{Adamsjust}
	The original $\Gamma$-space with rigid sections described in example \ref{GLPCY} has by now been heavily obscured by several equivalences. The reader may therefore wish to know why we are justified in calling the map $\psi^q$ ``the Adams operation". Notice that from our construction of $\psi^q$, we can easily see that $\psi^q$ is of degree $q^i$ on any fiber $\mathbb{S}^i$ of $\widetilde{\pi}_{(p)}(\mathbb{S})$. Notice also that $\psi^q$ has the same effect on the rational cohomology of $\widetilde{\mathscr{B}}{\GL}_{(p)}^{top}(\C)_I$ as that of an unstable Adams operations by the same name. It is well-known that this uniquely determines the underlying map $\psi^q$ on $\widetilde{\mathscr{B}}{\GL}_{(p)}^{top}(\C)_I$ up to homotopy (as a map of spaces) and therefore justifies our calling $\psi^q$ the Adams operation. 
\end{remark}

\medskip
\noindent
We are now ready to apply our framework to the proof of the stable $p$-local Adams conjecture. 

\bigskip
\noindent
{\bf  The stable $p$-local Adams conjecture}

\medskip
\noindent
In theorem \ref{infinityfib}, given a sectioned $X$-fibration $\pi : \mathscr{E} \longrightarrow \mathscr{B}$, we functorially construct a $\Gamma$-space $\BGx_X(\mbox{P}_\pi)$ that is equivalent to $\mathscr{B}$, and comes endowed with a map 
\[ \mathfrak{J}_\pi^{\{0\}} : \BGx_X(\mbox{P}_\pi) \longrightarrow \BGx_X
\]
that ``classifies" the sectioned $X$-fibration $\pi$. 
Now, let us take for $X$ some fixed model $\mathbb{S}_{(p)}$ for the $p$-local two sphere and let $\pi$ be the fibration constructed in theorem \ref{main} above. For the sake of simplicity, let us denote the map $\mathfrak{J}_\pi^{\{0\}}$ by $\mathfrak{J}$, and rename the domain of $\mathfrak{J}$ by $\BGL_{(p)}^{top}(\C)$. We therefore obtain a well-defined map of $\Gamma$-spaces:
\[ \mathfrak{J} : \BGL_{(p)}^{top}(\C) \longrightarrow \BGx_{\mathbb{S}_{(p)}} \]
Furthermore, by naturality, this map is fixed under the Adams operation $\psi^q$, i.e. we have $\mathfrak{J} = \mathfrak{J} \circ \psi^q$. Let $\Cq$ denote the homotopy co-equalizer (see below) of $\psi^q$ and the identity map, on $\BGx_{\mathbb{S}_{(p)}}$. This homotopy co-equalizer is also the homotopy orbits of the action of $\psi^q$ seen as the action of the monoid $\N$. One therefore has an extension: 
\[ \mathfrak{J}_q : \Cq \longrightarrow \BGx_{\mathbb{S}_{(p)}}. \]
Of course, such an extension $\mathfrak{J}_q$ may not be unique. However, we may make a canonical choice for $\mathfrak{J}_q$ as follows. Consider the following model for $\Cq$ given by the mapping cylinder of $\psi^q$, $\Delta[1] \times_{\psi^q} \BGL^{top}_{(p)}(\C)$ which is defined as the pushout of $\Gamma$-spaces:
\[
\xymatrix{
	{(\Delta[0]_0 \coprod \Delta[0]_1 \coprod \ast) \wedge \BGL_{(p)}^{top}(\C)} \ar[d] \ar[r]^{\quad \quad \quad \quad id \coprod \psi^q} \ar@{}[dr]|{\text{PO}}& \BGL_{(p)}^{top}(\C) \ar[d] \\
	(\Delta[1] \coprod \ast) \wedge \BGL^{top}_{(p)}(\C)  \ar[r] & \Delta[1] \times_{\psi^q} \BGL^{top}_{(p)}(\C),}
\]
where the left vertical map is the canonical inclusion. Consider the composite map:
\[ (\Delta[1] \coprod \ast) \wedge \BGL^{top}_{(p)}(\C) \longrightarrow \BGL^{top}_{(p)}(\C) \longrightarrow \BGx_{\mathbb{S}_{(p)}}, \]
where the first map is the canonical projection, and the second map is $\mathfrak{J}$. Using naturality of $\mathfrak{J}$ as shown in theorem \ref{infinityfib}, the above composite factors uniquely through the pushout $\Delta[1] \times_{\psi^q} \BGL^{top}_{(p)}(\C)$.
The extension $\mathfrak{J}_q$ we seek is this canonical extension:
\begin{equation} \label{Jq} \mathfrak{J}_q : \Cq = \Delta[1] \times_{\psi^q}  \BGL^{top}_{(p)}(\C) \longrightarrow \BGx_{\mathbb{S}_{(p)}}. \end{equation}
Now let us recall the functor $\nabla$ (\ref{nabla}). By \cite{BF} (Lemmas 4.1 and 4.6), we see that $\nabla$ preserves colimits and that one has a canonical map which is an equivalence: 
\[  \Cqp \longrightarrow \nabla (\Cq), \]
where $\Cqp$ denotes the homotopy coequalizer of the identity map and $\nabla(\psi^q)$. 
Applying $\nabla$ to $\mathfrak{J}_q$ we obtain a map of spectra, which we shall denote by $J_q$: 
\begin{equation} \label{Jq2} J_q : \Cqp \longrightarrow \nabla (\Cq) \longrightarrow \nabla(\BGx_{\mathbb{S}_{(p)}}), \end{equation}
Now let $\Eqp$ denote the homotopy equalizer of $\nabla \psi^q$ and the identity map, on $\nabla (\BGx_{\mathbb{S}_{(p)}})$. It is straightforward to see that $\Cqp$ is weakly equivalent to $\Sigma \Eqp$. In particular, we see that:
\[ \pi_1(\Cqp) = \pi_0(\Eqp) = \Z, \quad \quad \pi_1(\nabla (\BGx_{\mathbb{S}_{(p)}})) = \Z_{(p)}^{\times}. \]

\smallskip
\noindent
\begin{thm} \label{imp}
	The map $J_q : \Cqp \longrightarrow \nabla (\BGx_{\mathbb{S}_{(p)}})$ sends a generator of $\pi_1(\Cqp)$ to the element $q \in \Z_{(p)}^{\times} = \pi_1(\nabla (\BGx_{\mathbb{S}_{(p)}}))$. 
\end{thm}
\begin{proof}
	Notice that $J_q$ identifies $\pi_0(\Cqp)$, with $\pi_0(\nabla (\BGL_{(p)}^{top}(\C))) = \Z$. Given a spectrum $\mathscr{S}$, let $\Omega^{\infty} \mathscr{S}$ denote the simplicial set $(\mbox{Q} \mathscr{S})_0$, where $\mbox{Q}$ was the functor that replaces a spectrum by an $\Omega$-spectrum as defined in \ref{Q}. Hence the components of $\Omega^{\infty} \Cqp$ and $\Omega^{\infty} \nabla (\BGL_{(p)}^{top}(\C))$ can be identified with $\Z$. We will use the notation $\Omega^{\infty}_1$ to denote the component indexed by $1 \in \Z$. 
	
	\medskip
	\noindent
	In analogy with the mapping cylinder constructed above, let us define the simplicial set $\Delta[1] \times_{\psi^q} \BGLL^{top}_1(\C)_{(p)}$ as the pushout diagram of simplicial sets: 
	\[
	\xymatrix{
		(\Delta[0]_0 \coprod \Delta[0]_1 ) \times \BGLL^{top}_1(\C)_{(p)} \ar[d] \ar[r]^{\quad \quad \quad id \coprod \psi^q} \ar@{}[dr]|{\text{PO}}& \BGLL^{top}_1(\C)_{(p)}  \ar[d] \\
		\Delta[1]  \times \BGLL^{top}_1(\C)_{(p)}  \ar[r] & \Delta[1] \times_{\psi^q} \BGLL^{top}_1(\C)_{(p)}.}
	\]
	By equations \ref{Jq}, \ref{Jq2} and the definition of the functor $\nabla$, one notices that there is a pair of canonical maps, which can be seen to be isomorphisms on the fundamental group: 
	\[ \Delta[1] \times_{\psi^q} \BGLL^{top}_1(\C)_{(p)} \longrightarrow \Omega_1^{\infty} \Cqp, \quad \quad \quad \BG(\mathbb{S}_{(p)})  \longrightarrow \Omega_1^{\infty} \nabla (\BGx_{\mathbb{S}_{(p)}}). \]
	Moreover, equations \ref{Jq} and \ref{Jq2} also show that one has a commutative diagram:
	\[
	\xymatrix{
		\Delta[1] \times_{\psi^q} \BGLL^{top}_1(\C)_{(p)} \ar[d]^{J_q} \ar[r] & \Omega_1^{\infty} \Cqp \ar[d]^{J_q} \\
		\BG(\mathbb{S}_{(p)}) \ar[r] & \Omega_1^{\infty} \nabla (\BGx_{\mathbb{S}_{(p)}}) .}
	\]
	Now the map $\psi^q$ induces an automorphism of degree $q$ on the fiber over the base point over $\BGLL^{top}_1(\C, \mathbb{S})_{(p)} \longrightarrow \BGLL^{top}_1(\C)_{(p)}$ using Theorem \ref{main}. In particular, the left vertical map in the above diagram sends the generator of $\pi_1(\Delta[1] \times_{\psi^q} \BGLL^{top}_1(\C)_{(p)}) = \Z$ to the element $q \in \Z_{(p)}^\times =  \pi_1(\BG(\mathbb{S}_{(p)}))$. The result follows by chasing the diagram. 
\end{proof}

\medskip
\noindent
Finally, we get to the $p$-local stable Adams conjecture. Recalling the notation from the introduction, we have: 

\bigskip
\begin{thm} \label{PLAC}
	Fix $p,q$ to be any primes such that $p \neq q$. Recall the classical $p$-local $J$-homomorphism: 
	\[ J : \Z \times \BUn_{(p)} \longrightarrow \Pic^{ev} {\bf S}_{(p)}, \]
	where $\Pic^{ev} {\bf S}_{(p)}$ is the even Picard space defined as $\Omega^{\infty} \pic^{ev} {\bf S}_{(p)}$. Then $J$ lifts to a stable map 
	\[ \underline{J} : \underline{\ku}_{(p)} \longrightarrow \pic^{ev}  {\bf S}_{(p)} \]
	such that $\underline{J}$ is invariant under precomposition with the Adams operation $\psi^q$. More precisely, $\underline{J}$ admits a  canonical factorization through a map $\underline{J}_q$: 
	\[
	\xymatrix{
		& \ar[dl]  \underline{\ku}_{(p)}   \ar[d]^{\underline{J}}  \\
		\underline{\ku}_{h\psi^q} \ar[r]^{\underline{J}_q} & \pic^{ev} {\bf S}_{(p)}}
	\]
	where $\underline{\ku}_{h\psi^q}$ denotes the homotopy co-equalizer (or homotopy orbits) of the $\psi^q$-action on $\underline{\ku}_{(p)}$. Furthermore, the map $\underline{J}_q$ sends a generator of $\pi_1(\underline{\ku}_{h\psi^q}) = \Z$ to the $q \in \Z_{(p)}^{\times} = \pi_1(\pic^{ev} {\bf S}_{(p)})$. 
\end{thm}
\begin{proof}
	Recall the classifying map of $\Gamma$-spaces:
	\[ \mathfrak{J} : \BGL_{(p)}^{top}(\C) \longrightarrow \BGx_{\mathbb{S}_{(p)}}. \]
	Evaluating $\mathfrak{J}$ on the pointed set ${\bf 1}$, taking the geometric realization and invoking theorem \ref{Segal}, we obtain an $E_\infty$-map which is equivalent to 
	\[ J({\bf 1}) : \coprod_{k \geq 0} (\BGb^{top}_k(\C))_{(p)} \longrightarrow \coprod_{k\geq 0} \BG(|\mathbb{S}_{(p)}^{\wedge k}|),
	\]
	where $\G(|\mathbb{S}_{(p)}^{\wedge k}|)$ is the grouplike topological monoid of self-homotopy equivalences of the $i$-fold smash product of $|\mathbb{S}_{(p)}|$, where $\mathbb{S}_{(p)}$ is a simplicial model $\mathbb{S}_{(p)}$ for the $p$-local 2-sphere. 
	From our construction, it is clear that $J(\bf 1)$ is equivalent to the unstable $p$-local $J$-homomorphism that classifies the localization of the underlying (even dimensional) spherical fibration of a complex vector bundle obtained by taking the one-point compatification of the fibers. 
	
	\medskip
	\noindent
	Our main task is to identify the $E_{\infty}$-structure on the domain and codomain of $J({\bf 1})$ with the standard $E_\infty$-structure that group completes to $\Omega^{\infty} \underline{\ku}_{(p)}$ and $\Pic^{ev}{\bf S}_{(p)}$ respectively. Having done so, we may apply 
	apply the geometric realization functor to theorem \ref{imp}, and invoke theorem \ref{Segal} to conclude the proof. 
	
	\medskip
	\noindent
	Notice that theorem \ref{main} shows that the values of the $\Gamma$-space $\BGL_{(p)}^{top}(\C)$ are naturally equivalent to the $p$-localization of the respective values of the $\Gamma$-space $\BGL^{top}(\C)$ of example \ref{GLPCY}. We have already identified the spectrum $\nabla(\BGL^{top}(\C))$ with $\ku$. It follows from this, and theorem \ref{Segal}, that $\nabla(\BGL_{(p)}^{top}(\C))$ is equivalent to $\underline{\ku}_{(p)}$.
	
	\medskip
	\noindent
	Before addressing the $E_\infty$-structure on the codomain of $J({\bf 1)}$ we need a brief digression. Notice that given a topological permutative category $\mathscr{A}$, there is an underlying infinity category $\mbox{N}(\mathscr{A})$ given by the homotopy coherent nerve of $\mathscr{A}$. The infinity category $\mbox{N}(\mathscr{A})$ admits a monoidal structure encoded by the operadic nerve $\mbox{N}^{\otimes}(\mathscr{A})$ (\cite{Lu}, 2.1.1.25). Now let us consider the case when $\mbox{N}(\mathscr{A})$ is an infinity groupoid,  i.e. when all morphism are invertible in the higher categorical sense as is the case in all our examples of topological permutative categories. Then $\mbox{N}(\mathscr{A})$ can be identified with the classifying space of $\mathscr{A}$. Moreover, the monoidal structure on $\mbox{N}(\mathscr{A})$ gives rise to an $E_\infty$-structure, which is precisely the $E_\infty$-structure described by Segal on the classifying space of the symmetric monoidal category $\mathscr{A}$ (\cite{Se}, Section 2). 
	
	\smallskip
	\noindent
	Now consider the graded topological permutative category $\G_X$ described in example \ref{CXX}, where $\mbox{X}$ is taken to be the simplicial $p$-local 2-sphere $\mathbb{S}_{(p)}$. The objects of $\G_{\mathbb{S}_{(p)}}$ are the natural numbers with the usual monoidal structure. The morphisms are automorphisms of the objects $n \in \N$, which are given by the space of self homotopy equivalences of $|{\mathbb S}_{(p)}^{\wedge n}|$ if $n > 0$, and trivial otherwise. 
	One has a map of symmetric monoidal infinity categories $\varphi : \mbox{N}(\G_{\mathbb{S}_{(p)}}) \longrightarrow \mathscr{P}^{ev}({\bf S}_{(p)})$, where $\mathscr{P}^{ev}({\bf S}_{(p)})$ denotes the even Picard category of the $p$-local sphere spectrum ${\bf S}_{(p)}$, defined as the symmetric monoidal infinity category of invertible ${\bf S}_{(p)}$-modules equivalent to even suspensions of ${\bf S}_{(p)}$, and equivalences thereof. The functor $\varphi$ sends $k \in \N$ to the module $|\mathbb{S}^{\wedge k}_{(p)}| \wedge {\bf S}_{(p)}$ if $k > 0$, and to the unit ${\bf S}_{(p)}$ if $k=0$. Since  $\mbox{N}(\G_{\mathbb{S}_{(p)}})$ and $\mathscr{P}^{ev}({\bf S}_{(p)})$ are infinity groupoids, $\varphi$ induces a map of $E_\infty$-spaces:
	\[ \varphi : \coprod_{k\geq 0} \BG(|\mathbb{S}_{(p)}^{\wedge k}|) \longrightarrow \Pic^{ev} {\bf S}_{(p)} = \Omega^{\infty} \pic^{ev} {\bf S}_{(p)}. \]
	On stabilization the above map (\cite{Se}, Section 4), we see that it is an equivalence on topological group completion and hence induces an equivalence of the underlying spectra. 
\end{proof}

\appendix

\section{Friedlander's classification result for sectioned $X$-fibrations}  \label{classifyX}

\noindent
{\bf{Comparison with work of E. Friedlander:}}

\medskip
\noindent
Friedlander proves a classification theorem for sectioned $X$-fibrations (see \cite{F}, Theorem 6.1; \cite{F3}, Theorem, 4.10). For our purposes, we will only require the existence of a functorial classifying map, for which we offer a simpler proof below. The reason this proof can be shortened is that we have made an assumption that rigid sections, by definition, vary functorially in the Reedy category of rigid subsets $\mathscr{F}(\bf n)$ (see definition \ref{RigSet}). In \cite{F} and \cite{F3} Friedlander defines a notion of sections that vary functorially only along the inclusion maps in $\mathscr{F}(\bf n)$. We observe that the full Reedy structure is available to us in all examples of interest since they all arise from graded topological or algebraic permutative categories. Hence, employing this structure does not result in any loss of generality; on the contrary, it simplifies Friedlander's argument. 

\medskip
\noindent
As was already pointed out in remark \ref{d:S}, there is one other crucial difference between our constructions and the one described by Friedlander in \cite{F}. Given any pointed subset $U \subseteq {\bf m}$, a morphism $d : {\bf m} \longrightarrow {\bf n}$ in $\mathscr{F}$ and $I \in \N^{\times n}$, Friedlander defines a set $(d:U) \subseteq {\bf n} $ (see \cite{F}, Definition 2.1) that is analogous to our definition of $d_I(U)$ (see definition \ref{RigSet}). Furthermore, Friedlander claims that for a sectioned $X$ fibration $\mathscr{E} \longrightarrow \mathscr{B}$, the induced map $d^{\mathscr{E}} : \mathscr{E}_I({\bf m}) \longrightarrow \mathscr{E}_{d(I)}({\bf n})$ maps $\mathscr{E}_I^U({\bf m})$ to $\mathscr{E}^{(d:U)}_{d(I)}({\bf n})$. This claim is false as one may see from the following example: Let $d : {\bf 2} \longrightarrow {\bf 1}$ be the map that sends both nonzero elements of ${\bf 2}$ to the nonzero element $1 \in {\bf 1}$. Let $I = (1,0) \in \N^{\times 2}$ and $U = \{0,1\} \subset {\bf 2}$. One checks that $d_I(U) = {\bf 1}$, while $(d:U) = \{ 0 \}$. Since $\mathscr{E}({\bf 1})$ admits the structure of a monoid induced by the map $d^{\mathscr{E}}$ with $d$ as given above, the discussion following example \ref{nabla} shows that the image of $\mathscr{E}_I^U({\bf 2})$ under $d^{\mathscr{E}}$ cannot factor through $\mathscr{E}^{\{ 0 \}}_{d(I)}({\bf 1})$. 

\medskip
\begin{lemma} \label{Reedy}
	There is a functor $\mathscr{E} \mapsto \rE$ on the category of $\Gamma$-spaces with rigid sections, endowed with a natural transsormation $\eta_r : \mathscr{E} \longrightarrow \rE$ so that the following hold:
	
	\smallskip
	\noindent
	Define $\rB$ as $\rE^{\{ 0 \}}$ and $\mathscr{B}$ as $\mathscr{E}^{\{0\}}$ as in definition \ref{RigS}. 
	Then the $\Gamma$-space $\rB$ is fibrant, and the map $\eta_r^{\{0\}}: \mathscr{B} \longrightarrow \rB$ is an acyclic cofibration in the model category structure on $\Gamma$-spaces.
	
	\smallskip
	\noindent
	Furthermore, for a fixed $I \in \N^{\times n}$, the functor $S \mapsto \rE_I^S$ defined on the category $\mathscr{F}(\bf n)$, is Reedy fibrant and the natural transformation $\eta_r : \mathscr{E}_I^S \longrightarrow \rE_I^S$ is a Reedy acyclic cofibration when seen as a natural transformation of functors on $\mathscr{F}({\bf n})$. 
	
	\smallskip
	\noindent
	As a standard consequence of the above, we see that for a fixed pointed subset $S \subseteq {\bf n}$ and $I \in \N^{\times n}$, the map described in definition \ref{RigS}, $\pi_I^S : \rE_I^S \longrightarrow \rB_I$ is a fibration of simplicial sets and $\eta_r : \mathscr{E}_I^S \longrightarrow \rE_I^S$ is an acyclic cofibration of simplicial sets. 
\end{lemma}

\begin{proof}
	The category of finite pointed sets $\mathscr{F}$ may be seen as a generalization of a  Reedy category \cite{BM}. As such, functorial constructions on diagrams on $\mathscr{F}$, i.e. $\Gamma$-spaces, are typically made by extending matching objects by latching objects using an inductive argument on the size of the finite set ${\bf n}$ (see for instance \cite{BF}, Lemma 3.6). Similarly, for a fixed value of $n$, the category $\mathscr{F}({\bf n})$ is a Reedy category. As before, functorial constructions of diagrams on $\mathscr{F}({\bf n})$ are made by extending matching objects by latching objects using an inductive argument on the size of the subset $S \subseteq {\bf n}$ (see \cite{RV}, Lemma 3.10). Since $\Gamma$-spaces with rigid sections are diagrams on both $\mathscr{F}$ and $\mathscr{F}({\bf n})$, (with some compatibility between them), we construct our functor $\rE$ by extending a suitable notion of matching object, by a similar notion of latching object (see equation \ref{factorization}) using simultaneous induction on $n$ and on the size of the rigid subsets $S \subseteq {\bf n}$. Our construction can be compared with \cite{F} (Proposition 2.5). 
	
	\medskip
	\noindent
	One begins the induction by defining $\rE^{\{ 0 \}}$ as the $\Gamma$-space $\rB$ that fits into the functorial factorization of the projection map $\mathscr{B} \longrightarrow \ast$ in the category of $\Gamma$-spaces
	\[ \mathscr{B} \llra{Q} \rB \llra{R} \ast, \]
	where $Q$ is an acyclic cofibration, and $R$ is a fibration in the model category structure on $\Gamma$-spaces. In particular, $\rB$ is fibrant and $\eta_r^{\{0\}} := Q$ is a an acyclic cofibration as $\Gamma$-spaces, by construction. By definition, the fibrancy of $\rB$ amounts to the following $\Sigma_n$-equivariant map being a fibration of simplicial sets:
	\begin{equation} \label{rB}
		\rB({\bf n}) \longrightarrow \invlim_{d \in n > t} \, \, \rB(\bf t), 
	\end{equation}
	where $d \in n > t$ denotes the category of all morphisms $d : {\bf n} \longrightarrow {\bf t}$ for $n > t$. The acyclic cofibrancy of $Q$ amounts to the following map being a acyclic cofibration in $\Sigma_n$-simplicial sets (see \cite{BF}, Proposition 3.3)
	\begin{equation} \label{rB2}
		\dirlim_{\beta \in m < n} \rB({\bf m}) \cup_{ \dirlim_{\beta \in m < n}  \mathscr{B}({\bf m})} \mathscr{B}({\bf n}) \longrightarrow \rB({\bf n}), 
	\end{equation}
	where $\beta \in m < n$ denotes the category of morphisms $\beta : {\bf m} \longrightarrow {\bf n}$ in $\mathscr{F}$ for $m < n$. Since $Q$ is a weak equivalence, the components of $\rB$ fibering over $I \in \N^{\times n}$ are well defined by those of $\mathscr{B}$. Hence $\rB$ is a $\Gamma$-space over $\mathscr{N}$.
	
	\medskip
	\noindent
	Now assume by induction that the functor 
	$\rB$ and the natural transformation $\eta_r$ have been functorially extended to a $\Gamma$-space $\rE$ with rigid sections on the full subcategory of $\mathscr{F}$ consisting of objects ${\bf k}$ for $k < n$. Further assume that $\rE^{\tilde{S}}_I$ has been constructed for all subsets $\tilde{S} \subseteq {\bf n}$ of cardinality less than some integer $i$ so that the properties described in definition \ref{RigS} hold. Let us recast these properties in the form of a functorial factorization \ref{factorization} constructed as follows.
	
	\medskip
	\noindent
	Given $I \in \N^{\times n}$, consider spaces $\partial^1_I, \partial^2_I, \Delta_I$ which are assumed to exist by induction: 
	\[ \partial^1_I := \dirlim_{\substack{\beta \in m < n \\ \beta(J) = I \\ T \subsetneq \tilde{S}}}   \mathscr{E}_J^{\beta^{-1}(\tilde{S})} \cup_{\mathscr{E}_J^{\beta^{-1}(T)}}    \rE_J^{\beta^{-1}(T)} \quad \quad \partial^2_I :=  \dirlim_{T \subsetneq \tilde{S}} \, \mathscr{E}_I^T  \quad \quad \Delta_I := \invlim_{\substack{\tilde{S} \twoheadrightarrow T\neq \tilde{S}, \\ d \in n > t}} \rE_{d(I)}^{d_I(T)}, \]
	where we have kept the objects of $\mathscr{F}$ implicit, and we will continue to do so for the rest of the proof. Notice that an arbitrary morphism $d$ in $\mathscr{F}$  factors as a surjection followed by an injection $d = \beta_0 \circ d_0$. This factorization is unique, up to an action of a symmetric group. It now follows by invoking lemma \ref{RigSet2} that our induction assumption above is equivalent to the existence of a functorial factorization as below for all objects $\tilde{S} \subseteq {\bf n}$ of cardinality less than $i$: 
	\begin{equation} \label{factorization}
		\dirlim_{\substack{\beta \in m < n \\ \beta(J) = I }} \rE_J^{\beta^{-1}(\tilde{S})} \cup_{\partial^1_I} (\mathscr{E}_I^{\tilde{S}} \cup_{\partial^2_I} \dirlim_{T \subsetneq \tilde{S}} \,  \rE_I^T) \llra{Q} \rE_I^{\tilde{S}} \llra{R} \invlim_{d \in n > t}  \rE_{d(I)}^{d_I(\tilde{S})} \times_{\Delta_I} \invlim_{\tilde{S} \twoheadrightarrow T\neq \tilde{S}} \rE_I^{T}
	\end{equation}
	in which $Q$ is an acyclic cofibration, and $R$ is a fibration of simplicial sets. We also demand that the symmetric group $\Sigma_n$, which we identify with the group of automorphisms of ${\bf n}$ in $\mathscr{F}$, acts so that $\sigma \in \Sigma_n$ maps $\rE_I^{\tilde{S}}$ to $\rE_{\sigma(I)}^{\sigma(\tilde{S})}$ compatibly with the action of $\sigma$ on the categories that index the two ends of \ref{factorization} given by sending $\beta \in m < n$ to $\sigma \circ \beta$, and sending $d \in n > t$ to $d \circ \sigma^{-1}$. 
	
	\smallskip
	\noindent
	Next, we construct factorization \ref{factorization} for subsets $S$ of cardinality $i$. 
	Given $S$ of cardiality $i$, consider the canonical map below, which exists by our induction assumption:
	\begin{equation} \label{pq2} \dirlim_{\substack{\beta \in m < n \\ \beta(J) = I }} \rE_J^{\beta^{-1}(S)} \cup_{\partial^1_I} (\mathscr{E}_I^S \cup_{\partial^2_I} \dirlim_{T \subsetneq S} \,  \rE_I^T) \llra{C_I(S)} \invlim_{d \in n > t}  \rE_{d(I)}^{d_I(S)} \times_{\Delta_I} \invlim_{S \twoheadrightarrow T\neq S} \rE_I^{T}.\end{equation}
	We functorially factor $C_I(S)$ (in the category of simplicial sets) into an acyclic cofibration $Q$ followed by a fibration $R$, and define $\rE_I^S$ as the extension 
	\begin{equation} \label{pq}
		\dirlim_{\substack{\beta \in m < n \\ \beta(J) = I }} \rE_J^{\beta^{-1}(S)} \cup_{\partial^1_I} (\mathscr{E}_I^S \cup_{\partial^2_I} \dirlim_{T \subsetneq S} \,  \rE_I^T)
		\llra{Q} \rE_I^S \llra{R} \invlim_{d \in n > t} \, \rE_{d(I)}^{d_I(S)} \times_{\Delta_I} \invlim_{S \twoheadrightarrow T\neq S} \rE_I^{T}. \end{equation}
	Now recall the action of $\sigma \in \Sigma_n$ on the domain and codomain of \ref{pq2} that acts on the indexing categories by sending $\beta \in m < n$ to $\sigma \circ \beta$ and $d \in n > t$ to $d \circ \sigma^{-1}$ respectively. This action maps $C_I(S)$ to $C_{\sigma(I)}(\sigma(S))$. Now using functoriality of the factorization \ref{pq}, we see that $\sigma \in \Sigma_n$ defines a canonical map from $\rE_I^S$ to $\rE_{\sigma(I)}^{\sigma(S)}$ giving rise to an action of $\Sigma_n$. Thus we have verified the induction step for the construction of $\rE_I^S$ based on the cardinality of $S$. 
	
	\medskip
	\noindent
	After completing the induction for all subsets $S \subseteq {\bf n}$ for a fixed value of $n$, we  may now proceed by induction on $n$. Recall that $\rE^{\{0\}} :=  \rB$ was previously defined. For $S = \{ 0 \}$, factorization \ref{factorization} is equivalent to equations \ref{rB} and \ref{rB2}. Hence we may increase $n$ and begin the induction process on the subsets $S \subseteq {\bf n}$ as before. By induction, we may thus construct the $\Gamma$-space with rigid sections $\rE$ as well as the natural transformation $\eta_r$ which is defined as the obvious component of the map $Q$. 
	
	\medskip
	\noindent
	Now consider the question of Reedy fibrancy which requires that the matching map $\mbox{M}_S$ below given by the second factor of $R$ in equation \ref{factorization} be a fibration of simplicial sets:
	\[ \mbox{M}_S :  \rE_I^S \longrightarrow \invlim_{S \twoheadrightarrow T\neq S} \rE_I^{T}. \]
	Since the map $R$ is a fibration, we see that proving $\mbox{M}_S$ is a fibration will follow by showing that following canonical map of simplicial sets is a fibration: 
	\[ \pi_{\Delta} : \invlim_{d \in n > t} \,\rE_{d(I)}^{d_I(S)} \longrightarrow \invlim_{\substack{S \twoheadrightarrow T\neq S \\ \, d \in n > t}} \rE_{d(I)}^{d_I(T)} = \Delta_I. \] 
	We prove that using an inductive lifting argument based on the cardinality of the set of nonzero elements in the image of the morphisms $d \in n > t$ as given below. 
	
	\medskip
	\noindent
	Consider the question of lifting any any horn inclusion $\Lambda^k[s]  \longrightarrow \Delta[s]$, along $\pi_{\Delta}$. We restrict the lifting diagram as follows, where $d_0 : {\bf n} \longrightarrow {\bf t}_0$ is chosen to be a fixed surjection: 
	\[
	\xymatrix{ \Lambda^k[s]  \ar[d] \ar[r] &
		\invlim_{d \in n > t} \,\rE_{d(I)}^{d_I(S)}   \ar[r] \ar[d]^{\pi_{\Delta}} & \rE_{d_0(I)}^{(d_0)_I(S)} \ar[r] \ar[d]^{\pi_1} & \invlim_{\alpha \in t_0 > t} \,\rE_{\alpha d_0(I)}^{\alpha_{d_0(I)}(d_0)_I(S)} \ar[d]^{\pi_2} \\
		\Delta[s] \ar[r] & \invlim_{\substack{S \twoheadrightarrow T \neq S \\ \, d \in n > t}} \rE_{d(I)}^{d_I(T)} \ar[r] & \invlim_{S \twoheadrightarrow T \neq S}  \,\rE_{d_0(I)}^{(d_0)_I(T)} \ar[r] & \invlim_{\substack{S \twoheadrightarrow T \neq S \\ \, \alpha \in t_0  > t}} \,\rE_{\alpha d_0(I)}^{\alpha_{d_0(I)}(d_0)_I(T)}. } 
	\]
	By induction, one assumes that $\Lambda^k[s]  \longrightarrow \Delta[s]$ has been endowed with a lift along $\pi_{\Delta}$ on factors indexed by maps $d : {\bf n} \longrightarrow {\bf p}$, with $n > p$ and whose image in ${\bf p} - \{0\}$ has cardinality less than $t_0$. In particular, we are given a lift of $\Lambda^k[s]  \longrightarrow \Delta[s]$ along $\pi_2$. 
	
	\medskip
	\noindent
	We now observe that the image of the functor $(d_0)_{I}$ on the poset $S \twoheadrightarrow T \neq S$ is the poset of retractions under $(d_0)_{I}(S)$ except possibly the identity retraction to $(d_0)_{I}(S)$ itself. Hence, one may replace $S$ by $(d_0)_{I}(S)$ in the last two inverse limits on the bottom row. Then using the fibrancy of the map $R$ on the right most square, we may extend the given lift along $\pi_1$. 
	
	\medskip
	\noindent
	Since any map $d : {\bf n} \longrightarrow {\bf p}$, with $n > p$ and whose image in ${\bf p} - \{0\}$ has cardinality $t_0$, factors through a surjection of the form $d_0$ (which is unique up to a free action of the permutation group $\Sigma_{t_0}$), our lift of $\Lambda^k[s]  \longrightarrow \Delta[s]$ along $\pi_1$ can be extended compatibly to all maps $d : {\bf n} \longrightarrow {\bf p}$, with $n > p$ and whose image in ${\bf p} - \{0\}$ has cardinality at most $t_0$. Using induction on $t_0$, we see that the trivial cofibration $\Lambda^k[s]  \longrightarrow \Delta[s]$ can be lifted along $\pi_{\Delta}$, showing that $\pi_{\Delta}$ is a fibration. 
	
	\medskip
	\noindent
	Next we proceed to show that the natural transformation $\eta_r : \mathscr{E}_I^S \longrightarrow \rE_I^S$ is a Reedy acyclic cofibration when seen as a natural transformation of functors on $\mathscr{F}({\bf n})$. The proof below is similar to
	\cite{F} (Proposition 2.5 (b)). The details are as follows.
	
	\medskip
	\noindent
	We work by induction on $n$ and assume that the above natural transformation is a Reedy acyclic cofibration for all integers less than $n$. 
	Since the map $Q$ of equation \ref{factorization} is an acyclic cofibration of simplicial sets by construction, we reduce our question to showing that the following canonical map is also an acyclic cofibration of simplicial sets:
	\[  \lambda :  \partial^1_I := \dirlim_{\substack{\beta \in m < n \\ \beta(J) = I \\ T \subsetneq S}}   \mathscr{E}_J^{\beta^{-1}(S)} \cup_{\mathscr{E}_J^{\beta^{-1}(T)}}    \rE_J^{\beta^{-1}(T)} \longrightarrow  \dirlim_{\substack{\beta \in m < n \\ \beta(J) = I }}  \rE_J^{\beta^{-1}(S)}
	.\] 
	Consider the question of lifting $\lambda$ along any fibration $\pi : E \rightarrow B$. We may restrict the lifting diagram as follows where $\beta_0 : {\bf m}_0 \longrightarrow {\bf n}$ is chosen to be a fixed injection and $J_0$ chosen to be a pointed subset (necessariy unique) in ${{\bf m}_0}$ so that $\beta_0(J_0) = I$: 
	\[
	\xymatrix{
		\mathscr{E}_J^{\alpha^{-1}(\beta_0^{-1}(S))} \cup_{\mathscr{E}_J^{\alpha^{-1}(T)}} \rE_J^{\alpha^{-1}(T)}  \ar[r]^{\quad \quad \quad i_2} \ar[d]^{j_2} & 
		\mathscr{E}_{J_0}^{\beta_0^{-1}(S)} \cup_{\mathscr{E}_{J_0}^{T}}\rE_{J_0}^{T}  \ar[r] \ar[d]^{j_1} &     \mathscr{E}_J^{\beta^{-1}(S)} \cup_{\mathscr{E}_J^{\beta^{-1}(T)}}    \rE_J^{\beta^{-1}(T)}    \ar[d]^{\lambda} \ar[r] & E \ar[d]^{\pi} \\
		\rE_J^{\alpha^{-1}(\beta_0^{-1}(S))}   \ar[r]  & \rE_{J_0}^{\beta_0^{-1}(S)} \ar[r] &  \rE_J^{\beta^{-1}(S)} \ar[r] & B,} 
	\]
	where we have used the convention of taking the direct limit over all variables. In the first column, these variables are $\{ \alpha, T, J\}$ with $\{ \alpha \in m < m_0, \alpha(J) = J_0, T \subsetneq \beta_0^{-1}(S)\}$. In the second column, the variable is $\{ T\}$, with $\{T \subsetneq \beta_0^{-1}(S)\}$ and it the third column they are $\{ \beta, T, J\}$, with $\{\beta \in m < n, \beta(J) = I, T \subsetneq S\}$. 
	
	\medskip
	\noindent
	As before, assume that a lift of $\lambda$ along $\pi$ has been constructed on the restriction of $\lambda$ to the image of $\beta : {\bf m} \longrightarrow {\bf n}$ with $m < n$ and whose image in ${\bf n}-\{0\}$ has cardinality less than $m_0$. Hence we are given a lift the map $j_2$ along $\pi$, giving rise to a canonical lift of the pushout $\mbox{P}(i_2,j_2)$ of the maps $i_2$ and $j_2$. 
	By the induction assumption on $n$, we see that the map from $\mbox{P}(i_2,j_2)$ to $\rE_{J_0}^{\beta_0^{-1}(S)}$ is an acyclic cofibration of simplicial sets. Hence we may extend our lift to $\rE_{J_0}^{\beta_0^{-1}(S)}$. 
	Now any morphism $\beta : {\bf m} \longrightarrow {\bf n}$ with $m < n$ and whose image in ${\bf n}-\{0\}$ has cardinality $m_0$, can be factored through an injection of the form $\beta_0$ which is unique up to a free action of the group $\Sigma_{m_0}$. This allows us to proceed by induction on the size of the image of $\beta$ to show that $\lambda$ may be lifted along $\pi$. Since $\pi$ was an arbitrary fibration, we conclude that $\lambda$ is an acyclic cofibration of simplicial sets. 
\end{proof}

\medskip
\begin{remark} \label{genReedy}
	Notice that factorization \ref{pq} shows that the $\Sigma_n$-action on simplices in the complement of the inclusion $\mathscr{E} \subseteq \rE$, need not be free. This observation is in contrast to the stronger condition of (equivariant) cofibrancy for the model category structure as claimed in \cite{F} (Proposition 2.5(b)). 
	
	\noindent
	The reader will notice that in contrast to \cite{F}, in this article we have refrained from describing a model structure on the category of $\Gamma$-spaces with rigid sections. However, lemma \ref{Reedy} suggests that an obvious choice for such a structure should be a hybrid of the Model structure on $\Gamma$-spaces, and the Model structure on functors defined on $\mathscr{F}({\bf n})$. 
\end{remark}


\bigskip
\noindent
Before we begin with Friedlander's classification theorem, let us set some notation. For a pointed simplicial set $X$ and $i > 0$, we define $X^i$ to be the $i$-fold smash product. Let $X^0$ be the constant point simplicial set. In what follows, we assume familiarity with the construction given in lemmas \ref{TPC1} and \ref{retracts}, example \ref{CXX}, and definition \ref{fib}. 

\medskip
\noindent
Applying lemma \ref{Reedy} to example \ref{CXX}, one obtains 

\medskip
\begin{defn} \label{univXfib}
	The universal sectioned $X$-fibration is defined as the map (see definition \ref{RigS})
	\[ \pi(X) : \rBGx_X(X) \longrightarrow \rBGx_X,\]
	where $\rBGx_X(X)$ is the value of the functor \ref{Reedy} applied to the $\Gamma$-space with rigid sections $\BGx_X(X)$ which was defined in example \ref{CXX}. 
\end{defn}

\medskip
\begin{thm} (\cite{F}) \label{infinityfib} There exists functors on the category of sectioned $X$-fibrations given by $\pi \mapsto \tilde{\pi}$ and $\pi \mapsto \tilde{\pi}(X)$, and natural equivalences $\eta_r \tilde{\iota}$ and $\Pi_\pi$, so that one has a zig-zag of natural transformations: 
	\[
	\xymatrix{
		\mathscr{E} \ar[d]^{\pi} \ar[r]^{\eta_r \tilde{\iota}} & \rtE \ar[d]^{\tilde{\pi}} & \rBGx_X(X, \mbox{P}_\pi)   \ar[l]_{\Pi_\pi \quad \quad } \ar[d]^{\tilde{\pi}(X)} \ar[r]^{\quad \mathfrak{J}_\pi} & \rBGx_X(X) \ar[d]^{\pi(X)}   \\
		\mathscr{B} \ar[r]^{\eta_r \tilde{\iota}^{\{0\}}} & \rtB &  \rBGx_X(\mbox{P}_\pi)  \ar[l]_{\Pi_\pi^{\{0\}}  \quad} \ar[r]^{\quad \,  \mathfrak{J}_\pi^{\{0\}}} & \rBGx_X.}	\]
	In the above diagram $\pi(X)$ denotes the constant functor on the category of sectioned $X$-fibrations supporting a (classifying) natural transformation $\mathfrak{J}_\pi$ from $\tilde{\pi}(X)$ to $\pi(X)$.
\end{thm}

\begin{proof}
	As in the proof of this theorem given in \cite{F}, the heart of the argument is to first functorially replace $\pi$ by an equivalent sectioned $X$-fibration which admits an underlying ``principal bundle". This principal bundle will then allow us to construct a zigzag of natural morphisms of sectioned $X$-fibrations to the universal sectioned $X$-fibration $\pi(X)$. 
	
	\smallskip
	\noindent
	Recall that the examples \ref{CXX} and \ref{CXX} are constructed as singular nerves as defined in \ref{CS}. More precisely, in example \ref{CXX}, we first construct a simplicial topological space as the bar construction $\Ne (|X|_I, \G_I)$ followed by taking the diagonal of the bisimplical set of levelwise singular simplices. In order to relate this construction to the proof below, we will provisionally work in the category of (possibly simplicial) topological spaces. We remind the reader that, by convention, all topological spaces and topological constructions like products and mapping spaces are in the category of compactly generated weak hausdorff spaces. Taking the diagonal of the bisimplicial set obtained by applying the singular simplices levelwise will allow us to eventually return to simplicial sets. 
	
	\smallskip
	\noindent
	Given $\pi$ as above, let $\rE \longrightarrow \rB$ be its replacement under the functor constructed in lemma \ref{Reedy}. Let $\mbox{E}_I^\bullet$ be defined as the geometric realization $|\rE_I^\bullet|$. Then $\mbox{E}_I^\bullet$ is a $\Gamma$-space with rigid sections taking values in topological spaces (i.e. satisfying the axioms of Definition \ref{RigS} with values in topological spaces). We define $\mbox{B}_I$ to be $\mbox{E}_I^{\{ 0 \}}$, so that the map $|\pi|_I^S : \mbox{E}_I^S \longrightarrow \mbox{B}_I$ is identified with the map induced by the terminal projection $S \twoheadrightarrow \{ 0 \}$ in $\mathscr{F}(\bf n)$. Since geometric realization preserves finite inverse limits and turns  Kan-fibrations into Serre fibrations, lemma \ref{Reedy} shows that the functor $S \mapsto \mbox{E}_I^S$ is a Reedy fibrant functor on $\mathscr{F}(\bf n)$ with values in topological spaces. 
	
	\smallskip
	\noindent
	Fix $I = (i_1, \ldots, i_n) \in \N^{\times n}$, and let $|X|^I_\bullet$ be the functor on $\mathscr{F}(\bf n)$, $S \mapsto |X^{i_{s_1}}| \times \cdots \times |X^{i_{s_r}}|$ where $s_1 < s_2 < \ldots < s_r$ are the nonzero elements of $S$, and inclusions and retractions of $\mathscr{F}(\bf n)$ are defined to induce the obvious inclusions and projections of the factors. It is easy to see that $|X|^I_{\bullet}$ is Reedy cofibrant. Next, we consider the space of natural transformations from the functor $|X|^I_{\bullet}$ to $\mbox{E}_I^{\bullet}$. Since $\{ 0 \}$ is both initial and terminal in $\mathscr{F}(\bf n)$, and because $|X|^I_{\{ 0 \}} = \ast$, and $\mbox{E}_I^{\{ 0 \}} = \mbox{B}_I$, the following diagram commutes functorially in $S$ for any natural transformation $p^{\bullet}_I$ from the functor $|X|^I_{\bullet}$ to $\mbox{E}_I^{\bullet}$:
	\[
	\xymatrix{ |X|^I_S
		\ar[d] \ar[r]^{p_I^S} & \mbox{E}_I^S  \ar[d]^{|\pi|_I^S} \\
		\ast  \ar[r]^{p_I^{\{ 0 \}}} & \mbox{B}_I.}
	\]
	This shows that $p_I^\bullet$ is a pointed map that takes values in the fibers of $|\pi|_I^{\bullet}$ over the image of $p_I^{\{ 0 \}}$. Extending the classical definition of the pricipalization of a fibration \cite{M1}, we define $\mbox{P}_{\pi,I}$ to be the space of all natural transformations $p^{\bullet}_I$ from the functor $|X|^I_{\bullet}$ to $\mbox{E}_I^{\bullet}$, with the property that $p_I^S$ restricts to a pointed equivalence with the corresponding fiber of $|\pi|_I^S$. We give $\mbox{P}_{\pi,I}$ the compactly generated compact open topology. 
	
	\smallskip
	\noindent
	By definition, $\mbox{P}_{\pi,I}$ admits a left action by the monoid $\G^I := \prod_{k=1}^n \G(|X^{i_k}|)$. Let us study this action in some more detail. Let $\mbox{j} : \mathscr{F}_1(\bf n) \rightarrow \mathscr{F}(\bf n)$ denote the inclusion of the full subcategory $\mathscr{F}_1(\bf n)$ of pointed subsets $S \subseteq {\bf n}$ that contains at most one nonzero element. Let $\mbox{K} : \mbox{E}_I^{\bullet} \longrightarrow \mbox{Rj}_\ast \mbox{j}^\ast \mbox{E}_I^{\bullet}$ denote the canonical map from $\mbox{E}_I^\bullet$ to the homotopy right Kan-extension $\mbox{Rj}_\ast$ along $\mbox{j}$, of the restriction $\mbox{j}^\ast$ of the functor $\mbox{E}_I^\bullet$ to $\mathscr{F}_1(\bf n)$. Now, by condition {\bf (2)} of \ref{fib}, we see that $\mbox{K}$ is a weak equivalence. It is easy to verify that the right-Kan extension $\mbox{Rj}_\ast \mbox{j}^\ast \mbox{E}_I^{\bullet}$ is also Reedy-fibrant. Using the adjunction for homotopy right Kan extensions (\cite{R}, 10.3.7), it follows that $\mbox{P}_{\pi,I}$ is detected on $\mathscr{F}_1(\bf n)$, and is consequently a (left) principal bundle for $\mbox{E}_I$ over $\mbox{B}_I$ with structure monoid equivalent to  $\G^I$. 
	
	\smallskip
	\noindent
	Given a morphism $d : \bf m \longrightarrow \bf n$ and $I \in \N^{\times m}$, define $J := d(I)$ as in example \ref{exampleN}. Let $\sigma_d$ be the permutation that reorders all the elements of the ordered set $\{ d^{-1} \{ 1 \}, \ldots, d^{-1} \{n\} \}$ to be in increasing order, where each of the sets $d^{-1} \{ t \}$ are expressed in order induced from $\bf m$. Then $\sigma_d^{-1}$ followed by the monoidal product within each factor, extends to a natural surjective map $d(X) : |X|^I \longrightarrow |X|^J$, where we shall use the notation $|X|^I$ for $|X|^I_{\bf m}$, and $|X|^J$ for $|X|^J_{\bf n}$.  Under $d(X)$, the subspace $|X|^I_{d^{-1}(S)} \subseteq |X|^I$ maps surjectively onto $|X|^J_S \subseteq |X|^J$. Now using remark \ref{mu2}, we see that $d^{\E} : \mbox{E}_I \longrightarrow \mbox{E}_J$ when composed with any natural transformation in $p_I^{\bullet} \in \mbox{P}_{\pi,I}$ factors canonically through $d(X)$. Next, using condition {\bf (2)} of definition \ref{RigS} one obtains a well-defined natural transformation $\mbox{P}(d)(p_I^\bullet) : |X|^J_{S} \longrightarrow \mbox{E}_{J}^{S}$ of functors defined on $\mathscr{F}({\bf n})$ with objects denoted by subsets $S \subseteq {\bf n}$, and that has the property
	\[ \mbox{P}(d)(p_I^\bullet)  \circ d(X)\, (\lambda) = d^{\E} \circ p_I^\bullet \, (\lambda), \quad \mbox{for} \quad \lambda \in |X|^I_{d^{-1}(S)}. \]
	Finally, we show that $\mbox{P}(d)(p_I^\bullet)$ is in fact a natural transformation in $\mbox{P}_{\pi,J}$. Let $0 \neq s \in {\bf n}$. The natural transformation $p_I^{\bullet}$ induces an equivalence between the fiber over $\mathscr{B}_I({\bf m})$ of the homotopy colimit diagram as given in condition {\bf (3)} of \ref{fib}, with the same diagram but with $\mathscr{E}_I^T$ replaced by $|X|^I_T$. The latter homotopy colimit is easily seen to be equivalent to $|X|^{d(I)}_{\{0,s\}}$ while condition {\bf (3)} of \ref{fib} shows that the former homotopy colimit is equivalent to $\mbox{E}_J^{\{0,s\}}$. Invoking condition  {\bf (2)} of \ref{fib}, it follows easily that $\mbox{P}(d)(p_I^\bullet)$ must also be an equivalence with the fiber of $\mbox{E}_J^{\bullet}$ over $\mathscr{B}_J({\bf n})$.
	
	\smallskip
	\noindent
	From the above observation, we obtain a well-defined map $ \mbox{P}(d) : \mbox{P}_{\pi,I} \longrightarrow \mbox{P}_{\pi,J}$. 
	Moreover, $\mbox{P}(d)$ is seen to be compatible under the homomorphism from $\G^I$ to $\G^J$ induced using the naturality of $d(X)$ (see example \ref{CXX}). As pointed out in \cite{F} (see Proposition 4.2 therin), the maps $\mbox{P}(d)$ do not compose functorially in $d$. This problem can be remedied by replacing the bar construction defining the universal fibration $\pi(X)$ of \ref{CXX} by a two-sided bar construction (see definition \ref{CS}). This replacement will give rise to a map $\tilde{\pi}(X) : \BGx_X(X, \mbox{P}_\pi) \longrightarrow \BGx_X(\mbox{P}_\pi)$ supporting a zigzag of equivalences of sectioned $X$-fibrations to $\pi$. We present details below. 
	
	\smallskip
	\noindent
	We now apply lemmas \ref{TPC1} and \ref{retracts} to example \ref{CXX}. For $I = (i_1, \ldots, i_n)$, we consider: 
	\[ |X|_I = |X|^I \times \G_I^{+}, \quad \mbox{where} \quad \G_I^{+} = \prod_{T \subseteq \bf n} \G^{\times}_{\mu_I(T)}, \quad \mbox{with} \quad \mu_I(T) := \sum_{t \in T} i_t. \]
	In the above, $T$ denotes those pointed subsets of $\bf n$ that contain at least two nonzero elements, and $\G^{\times}_{\mu_I(T)}$ denotes the maximal subgroup of invertible elements within the topological monoid $\G(|X^{\mu_I(T)}|)$. Now $|X|_I$ admits a right action by the topological monoid: 
	\[ \G_I = \G^I \times \G_I^{+}. \]
	More generally, let $S \subseteq {\bf n}$ be a pointed subset with nonzero elements $s_1 < s_2 < \cdots < s_r$. Then the subspace of $|X|_I$ defined as $|X|_I^S := |X|^I_S \times \G_I^{+}$ is preserved under the $\G_I$-action. 
	As observed in theorem \ref{GRS}, the singular nerve construction applied to the bar construction $\Ne(|X|_I^S, \G_I)$ agrees with $\BGx_X(X)^S_I$, where $\BGx_X(X)$ is the $\Gamma$-space of example \ref{CXX}. 
	
	\smallskip
	\noindent
	Let us now extend the canonical left $\G^I$ action on $\mbox{P}_{\pi,I}$, trivially to an action of $\G_I$, and let $\Ne(|X|_I^S, \G_I, \mbox{P}_{\pi.I})$ denote the simplicial topological space given by the two sided bar construction (see definition \ref{CS}). 
	Notice that the simplices of the one sided bar construction $\Ne((|X|_I^S, \G_I)$ can be identified as a product factor (as topological spaces) in those of $\Ne(|X|_I^S, \G_I, \mbox{P}_{\pi,I})$. More precisely, as topological spaces, the $q$-simplices factor as:
	\begin{equation} \label{2sbar} \Ne(|X|_I^S, \G_I, \mbox{P}_{\pi,I})_q = |X|_I^S \times \G_I^{\times q} \times \mbox{P}_{\pi,I} =  \Ne(|X|_I^S, \G_I)_q \times \mbox{P}_{\pi,I}. \end{equation}
	Next, we study the functoriality of the two sided bar construction under morphisms in $\mathscr{F}$ by taking advantage of the above factorization. We first introduce some notation: 
	
	\smallskip
	\noindent
	{\bf Start of Notation}: Let $x_I = (x_t, x^{\times}_T)$ be a point in $|X|_I = |X|^I \times \G_I^{+}$, where $x_t$ denotes the components in $|X|^I = \prod_{k=1}^m |X^{i_k}|$ indexed by some nonzero $t \in {\bf m}$, and $x^{\times}_T$ denotes the components in $\G_I^{+} = \prod_{T \subseteq \bf m} \G^{\times}_{\mu_I(T)}$. Let $d : {\bf m} \longrightarrow {\bf n}$ be a morphism in $\mathscr{F}$ and let $J := d(I)$. Given a pointed subset $S \subseteq \bf n$, let us define $\overline{x}_d^{\times}(S)$ as the automorphism of $|X|^J_S$ expressed as factorwise product of group elements:  $x^{\times}_{d^\ast\{0,s_1\}} \times \cdots \times x^{\times}_{d^\ast\{ 0,s_r \}}$, where $s_1 < s_2 < \cdots < s_r$ are the nonzero elements of $S$ and $x^{\times}_{d^\ast\{0,s\}}$ is understood to be the trivial automorphism if $d^\ast\{0,s\} := d^{-1}\{s\} \cup \{0\}$ has fewer than two nonzero elements.
	{\bf End of Notation}.
	
	\smallskip
	\noindent
	Now notice that theorem \ref{GRS} and lemmas \ref{TPC1} and  \ref{retracts} allow us to explicitly describe the map of one sided bar constructions $\Ne(d) : \Ne(|X|_I, \G_I) \longrightarrow \Ne(|X|_J, \G_J)$. 
	We extend $\Ne(d)$ to the two-sided bar construction: $d^{\Pm} : \Ne(|X|_I, \G_I, \mbox{P}_{\pi,I}) \longrightarrow \Ne(|X|_J, \G_J, \mbox{P}_{\pi,J})$, by defining it on $q$-simplices that are factored as in equation \ref{2sbar}: 
	\[ d^{\Pm}(x_I, g_1, \ldots, g_q, p^\bullet_I) := (\Ne(d)(x_I, g_1, \ldots, g_q), p^\bullet_J), \]
	where the natural transformation $p^\bullet_J \in \mbox{P}_{\pi,J}$ is defined by $p^S_J(\lambda) := \mbox{P}(d) (p^\bullet_I) (\lambda \, \overline{y}_d^{\times}(S)^{-1})$ for any  $\lambda \in |X|^J_S$ and with $\overline{y}_d^{\times}(S)$ defined in our notation above as the invertible element when $y_I \in |X|_I$ is taken to be the element obtained by the sequential right action of the elements $g_r$ on $x_I$, i.e.  $y_I := x_I g_1 \cdots g_q$. Using the formulas in lemmas \ref{TPC1} and \ref{retracts}, one checks that $d^{\Pm}$ is a well-defined map of simplicial topological spaces that composes functorially in $d$. 
	
	\smallskip
	\noindent
	Now let $\BGx_X(X, \mbox{P}_\pi)^S_I $ denote the singular nerve of $\Ne(|X|_I^S, \G_I, \mbox{P}_{\pi,I})$ as in definition \ref{CS}. One readily sees that these spaces describe a $\Gamma$-space $\BGx_X(X, \mbox{P}_\pi)$ with rigid sections, so that the morphisms $d : {\bf m} \longrightarrow {\bf n}$ in $\mathscr{F}$ are defined to act as the maps $d^{\Pm}$ above. Similarly, the spaces $\Ne(\G_I^{+}, \G_I, \mbox{P}_{\pi.I})$ describe a $\Gamma$-space $\BGx_X(\mbox{P}_\pi)_I$. Moreover, the projection map 
	\[ \tilde{\pi}(X) : \BGx_X(X, \mbox{P}_\pi) \longrightarrow \BGx_X(\mbox{P}_\pi) \]
	is a map in $\mathscr{F}[\sSet_*]/\mathscr{N}$ which can be identified with the map induced by the terminal retraction as described in definition \ref{RigS}. 
	
	\smallskip
	\noindent
	We now proceed to describe a natural zigzag (compare \cite{F}, Propositon 5.6):
	\[
	\xymatrix{
		\mathscr{E} \ar[d]^{\pi} \ar[r]^{\tilde{\iota}} & \tilde{\mathscr{E}} \ar[d]^{\tilde{\pi}} & \BGx_X(X, \mbox{P}_\pi)   \ar[l]_{\Pi_\pi \quad \quad } \ar[d]^{\tilde{\pi}(X)} \ar[r]^{\quad \mathfrak{J}_\pi} & \BGx_X(X) \ar[d]^{\pi(X)}   \\
		\mathscr{B} \ar[r]^{\tilde{\iota}^{\{0\}}} & \tilde{\mathscr{B}} &  \BGx_X(\mbox{P}_\pi)  \ar[l]_{\Pi_\pi^{\{0\}}  \quad} \ar[r]^{\quad \,  \mathfrak{J}_\pi^{\{0\}}} & \BGx_X.}
	\]
	The top row will be morphisms of $\Gamma$-spaces with rigid sections, and the vertical maps projections as defined in definition \ref{RigS}. The bottom horizontal maps are the retraction of the top horizontal maps to the subset $\{ 0 \} \subseteq {\bf n}$ (so the the diagram commutes by definition).  
	
	\smallskip
	\noindent
	We begin by defining $\tilde{\pi} : \tilde{\mathscr{E}} \longrightarrow \tilde{\mathscr{B}}$ to be the map $\mbox{Sing} \, \mbox{E} \longrightarrow \mbox{Sing} \, \mbox{B}$ obtained by taking the singular simplices on the map $\mbox{E} \longrightarrow \mbox{B}$, which we recall was defined as the geometric realization $|\rE| \longrightarrow |\rB|$. Notice that the natural equivalence $\eta_r$ induces natural equivalence of sectioned $X$-fibrations $\tilde{\iota}: \mathscr{E} \longrightarrow \tilde{\mathscr{E}}$.
	The map $\mathfrak{J}_\pi$ in the above zigzag is induced by the map of simplicial topological spaces that collapes all $\mbox{P}_{\pi,I}$ to a point. 
	The map $\Pi_\pi$ is induced by a map of simplicial topological spaces with values in the constant simplicial topological space $\mbox{E}_I$, defined on $q$-simplices factored as in equation \ref{2sbar}, as the evaluation:
	\[ \Pi_\pi(x_I, g_1, \ldots, g_q, p^\bullet_I) := p^{\bf n}_I ((x_I g_1 \cdots g_q)^I), \quad \mbox{where} \quad I \in \N^{\times n}, \]
	and $(x_I g_1 \cdots g_q)^I$ indicates the component along the topological space $|X|^I = \prod_{k=1}^n |X^{i_k}|$, of the element $x_I g_1 \cdots g_q \in |X|_I = |X|^I \times \G_I^{+}$, which is defined as the sequential right action of the elements $g_r$ on $x_I$. Since $\mbox{P}_{\pi,I}$ was a principal bundle over $\mbox{B}_I$, it follows that $\Pi_\pi$ is a natural equivalence of $\Gamma$-spaces with rigid sections. 
	
	\smallskip
	\noindent
	Finally, we turn the above zigzag diagram into a natural zigzag of sectioned $X$-fibrations by applying the functor of lemma \ref{Reedy} to obtain the following diagram, where we have kept the names of the induced morphisms the same for simplicity: 
	\begin{equation} \label{infinityfib2}
		\xymatrix{
			\mathscr{E} \ar[d]^{\pi} \ar[r]^{\eta_r \tilde{\iota}} & \rtE \ar[d]^{\tilde{\pi}} & \rBGx_X(X, \mbox{P}_\pi)   \ar[l]_{\Pi_\pi \quad \quad } \ar[d]^{\tilde{\pi}(X)} \ar[r]^{\quad \mathfrak{J}_\pi} & \rBGx_X(X) \ar[d]^{\pi(X)}   \\
			\mathscr{B} \ar[r]^{\eta_r \tilde{\iota}^{\{0\}}} & \rtB &  \rBGx_X(\mbox{P}_\pi)  \ar[l]_{\Pi_\pi^{\{0\}}  \quad} \ar[r]^{\quad \,  \mathfrak{J}_\pi^{\{0\}}} & \rBGx_X.}
	\end{equation}
	It is clear that the construction of diagram \ref{infinityfib2} is natural with respect to morphisms of $\pi$. This establishes the proof of theorem \ref{infinityfib}. 
\end{proof}

\noindent
We end this section by relating the spectrum $\nabla(\mathscr{B})$ to $\nabla(\mathscr{E})$ for an $X$-fibration $\pi : \mathscr{E} \longrightarrow \mathscr{B}$

\medskip
\begin{thm}\label{stablefib}
	Assume $X$ is a pointed, connected simplicial set and $\pi : \mathscr{E} \longrightarrow \mathscr{B}$ is a sectioned $X$-fibration. Then the canonical map $\nabla(\pi) : \nabla(\mathscr{E}) \longrightarrow \nabla(\mathscr{B})$ is a weak equivalence of spectra. 
\end{thm}
\begin{proof}
	Let us begin by observing that $\nabla$ takes a morphism of $X$-fibrations to a homotopy pullback. This follows from the definition of the functor $\nabla$ (theorem \ref{nabla}) and the fact that $\nabla$ takes special $\Gamma$-spaces to $\Omega$-spectra on positively indexed spaces (theorem \ref{Segal}). Hence, by \ref{infinityfib}, we only need to prove our theorem in the special case of the universal fibration
	\[ \pi(X) : \rBGx_X(X) \longrightarrow \rBGx_X.\]
	In fact, we can do better. Consider the graded topological permutative category $\Sigma(X)$:
	\[ \mbox{Ob}(\Sigma(X)) = \coprod_{i \in \N} |X^i|, \quad  \mbox{Mor}(|X^i|,|X^j|) = \emptyset, \, \,  \mbox{if} \, \,  i \neq j, \quad  \mbox{Mor}(i,i) = |X^i| \times \Sigma_i . \]
	with the monoidal structure on objects induced by the smash product map on $|X^i|$ as in eample \ref{CXX}, and morphisms induced by the canonical right permutation action of $\Sigma_i$ on $|X^i|$. As always, we obtain a special $\Gamma$-space $\BSm(X)$ with rigid sections that comes endowed with an $X$-fibration:
	\[ \pi_\Sigma : \BSm(X) \longrightarrow \BSm. \]
	It is sufficient therefore to prove our theorem for $\pi_\Sigma$ since the fiber of $\nabla(\pi_\Sigma)$ agrees with the fiber of $\nabla(\pi(X))$. Up to equivalence, the monoids  $\BSm(X)({\bf 1})$ and $\BSm({\bf 1})$ are:
	\[ \BSm(X)({\bf 1}) = \coprod_{n \geq 0} \B(|X^n|,\Sigma_n), \quad \quad \BSm({\bf 1}) = \coprod_{n \geq 0} \BSi_n. 
	\]
	By the connectivity of $X$, the monoid of components of both these is $\N$. Furthermore, $\pi_{\Sigma}$ admits a $\{ 0 \}$-section. By \cite{Se}(Section 4), we see that the topological group completion of both these spaces factors through a stabilization procedure that maps the component indexed by the integer $n$, to the component indexed by the integer $(n+1)$. This map is given by multiplication with an element in the first component $\BSm_1 \subset \BSm({\bf 1})$. This element acts on $\BSm(X)({\bf 1})$ through the $\{0\}$-section. Notice that the product of any element of $\B(|X^n|,\Sigma_n)$ with an element in the image of the $\{0 \}$-section of $\BSm_1$ factors through the collapse map to  $\B(|X^n|,\Sigma_n) \longrightarrow \BSi_n$ (since this $\{0\}$-section takes values in the basepoint of $|X|$). Hence the stabilization of $\BSm(X)({\bf 1})$ and $\BSm({\bf 1})$ agree. The result follows. 
\end{proof}


\pagestyle{empty}
\bibliographystyle{amsplain}
\providecommand{\bysame}{\leavevmode\hbox
to3em{\hrulefill}\thinspace}

\end{document}